\documentclass[]{nmeauth}

\usepackage{subfigure}
\usepackage{endfloat}

\newcommand{\real}{\mathbb{R}} 
\newcommand{\integer}{\mathbb{N}}
\newcommand{\vectornorm}[1]{\left\|#1\right\|}

\usepackage[numbers,sort&compress]{natbib}
\usepackage[lined,boxed]{algorithm2e}

\begin{document}
\NME{0}{0}{00}{00}{00}

\runningheads{M. Arnst, R. Ghanem, E. Phipps, and J. Red-Horse}{Stochastic modeling of coupled problems}

\received{4 January 2010}
\norevised{}
\noaccepted{}

\title{Reduced chaos expansions with random coefficients\\in reduced-dimensional stochastic modeling of coupled problems}

\author{M. Arnst\affil{1}$^{,}$\affil{2}\corrauth, R. Ghanem\affil{2}, E. Phipps\affil{3}, and J. Red-Horse\affil{3}}

\address{\affilnum{1} B52/3, Universit\'{e} de Li\`{e}ge, Chemin des Chevreuils 1, B-4000 Li\`{e}ge, Belgium.\\
                 \affilnum{2} 210 KAP Hall, University of Southern California, Los Angeles, CA 90089, USA.\\
                 \affilnum{3} Sandia National Laboratories\footnotemark[2], P.O. Box 5800, Albuquerque, NM 87185, USA.}

\corraddr{B52/3, Universit\'{e} de Li\`{e}ge, Chemin des Chevreuils 1, B-4000 Li\`{e}ge, Belgium.}
\footnotetext[2]{Sandia National Laboratories is a multi-program laboratory managed and operated by Sandia Corporation, a wholly owned subsidiary of Lockheed Martin Corporation, for the U.S. Department of Energy's National Nuclear Security Administration under contract DE-AC04-94AL85000.}

\begin{abstract}
Coupled problems with various combinations of multiple physics, scales, and domains can be found in numerous areas of science and engineering.
A key challenge in the formulation and implementation of corresponding coupled models is to facilitate communication of information across physics, scale, and domain interfaces, as well as between iterations of solvers used for response computations.
In a probabilistic context, any information that is to be communicated between subproblems or iterations should be characterized by an appropriate probabilistic representation. 
In this work, we consider stochastic coupled problems whose subproblems involve only uncertainties that are statistically independent of one another; for these problems, we present a characterization of the exchanged information by using a reduced chaos expansion with random coefficients. 
This expansion provides a reduced-dimensional representation of the exchanged information, while maintaining segregation between sources of uncertainty that stem from different subproblems.
Further, we present a measure transformation that allows stochastic expansion methods to exploit this dimension reduction to obtain an efficient solution of subproblems in a reduced-dimensional space.
We show that owing to the uncertainty source segregation, requisite orthonormal polynomials and quadrature rules can be readily obtained by tensorization.
Finally, the proposed methodology is demonstrated by applying it to a multiphysics problem in nuclear engineering.
\end{abstract}


\section{Introduction}
In the fields of climate research, combustion, and renewable energy and in many other critical areas of science and engineering, models and simulations depend on a common base of mathematical formulations, algorithms, and implementations involving multiphysics, multiscale, and multidomain characteristics.  
Demand for predictive computational results in these areas has motivated the development of uncertainty quantification~(UQ) approaches for coupled problems with various combinations of multiple physics, scales, and domains.

Probability theory provides a rigorous mathematical framework for UQ.
The first step in a probabilistic UQ analysis typically involves the use of mathematical statistics methods~\citep{cramer1946,kullback1968} for the probabilistic characterization of the uncertain features associated with a model as one or more random variables, random fields, random matrices, or random operators.  
The second step is to use the model to map this probabilistic representation of inputs into a probabilistic representation of responses.
This can be achieved in several ways, for example, by using Monte Carlo sampling techniques~\citep{robert2005} and stochastic expansion methods.
The latter typically involve the representation of the predictions as a polynomial chaos~(PC) expansion.
Several approaches are available to calculate the coefficients in this expansion, for example, embedded projection~\citep{ghanem2003,soize2004}, nonintrusive projection~\citep{soize2004}, and collocation~\citep{ghanem1998,ghanem1998b,ghiocel2002,xiu2005,babuska2007}.

A key challenge in the formulation and implementation of coupled models is to facilitate communication of information across physics, scale, and domain interfaces, as well as between iterations of solvers used for response computations.
This information can comprise physical properties, energies, solution patches, and other quantities.
Although the number of sources of uncertainty can be expected to be large in most stochastic coupled problems, we believe that exchanged probabilistic information often resides in a considerably lower dimensional space than the sources themselves.
In stochastic multiphysics and stochastic multiscale problems, the exchanged information can be expected to be of \textit{low effective stochastic dimension} when it consists of a solution field that is smoothed by a forward operator and when it is obtained by summarizing fine-scale quantities into coarse-scale representations, respectively.

In a previous paper~\citep{arnst2011a}, we had thus proposed the use of a \textit{dimension-reduction} technique, namely, the Karhunen-Lo\`{e}ve~(KL) decomposition, to represent exchanged probabilistic information.
In a subsequent paper~\citep{arnst2011b}, we had presented a \textit{measure-transformation} technique that allows stochastic expansion methods to exploit this dimension reduction to obtain a computationally efficient solution of subproblems in a reduced-dimensional space.
In~\citep{arnst2011a} and~\citep{arnst2011b}, we had considered a general setting that allowed the uncertainties in the data of the various subproblems to be statistically dependent on one another.  

In this paper, we present an extension of our previous work.
This extension is applicable to stochastic coupled problems whose subproblems incorporate only sources of uncertainty that are statistically independent of one another; for these stochastic coupled problems, we present a characterization of the exchanged information by using a reduced chaos expansion with random coefficients.
This expansion has recently been proposed in~\citep{soize2009}, and in the current context, it enables us to obtain a reduced-dimensional representation of the exchanged information, a key feature of which is that it can maintain segregation between sources of uncertainty that stem from different subproblems.
Further, we present a corresponding measure-transformation technique that allows stochastic expansion methods to exploit this reduced-dimensional representation to obtain an efficient solution of subproblems in a reduced-dimensional space.
We show that owing to the uncertainty source segregation, requisite orthonormal polynomials and quadrature rules can be readily obtained by tensorization.

The remainder of this paper is organized as follows.
First, in Sec.~\ref{sec:sec2}, we outline the proposed methodology.
Next, in Sec.~\ref{sec:sec3}, we describe the reduced chaos expansion with random coefficients.
In Sec.~\ref{sec:sec6}, we provide details on the implementation.
Finally, in Secs.~\ref{sec:sec7} and~\ref{sec:sec8}, we demonstrate the proposed methodology by considering an illustration problem.

\section{Proposed methodology: Dimension reduction and measure transformation}\label{sec:sec2}

\subsection{Model problem}
This paper is devoted to the determination of the solution to a stochastic coupled model of the following form:
\begin{equation}
\label{eq:coupling12}\begin{aligned}
&\boldsymbol{f}(\boldsymbol{u},\boldsymbol{x},\boldsymbol{\xi})=\boldsymbol{0},&&\boldsymbol{y}=\boldsymbol{h}(\boldsymbol{u},\boldsymbol{\xi}),&&\boldsymbol{f}:\real^{r}\times\real^{s_{0}}\times\real^{m}\rightarrow\real^{r},&&\boldsymbol{h}:\real^{r}\times\real^{m}\rightarrow\real^{r_{0}},\\
&\boldsymbol{g}(\boldsymbol{y},\boldsymbol{v},\boldsymbol{\zeta})=\boldsymbol{0},&&\boldsymbol{x}=\boldsymbol{k}(\boldsymbol{v},\boldsymbol{\zeta}),&&\boldsymbol{g}:\real^{r_{0}}\times\real^{s}\times\real^{n}\rightarrow\real^{s},&&\boldsymbol{k}:\real^{s}\times\real^{n}\rightarrow\real^{s_{0}}.
\end{aligned}
\end{equation}
To avoid certain technicalities involved in infinite-dimensional representations, we assume that these equations are discretized representations of a stochastic model that couples two physical processes, two scales, two domains, or a combination of these subproblems.
For instance, these equations may be obtained by the spatial discretization of a steady-state problem or they may be obtained at a single time step after the spatial and temporal discretization of an evolution problem.
Further, we assume that the data of the first subproblem, which enter this subproblem as coefficients or loadings or both, depend on a finite set of uncertain real parameters denoted as $\xi_{1},\ldots,\xi_{m}$ and that the data of the second subproblem depend on a finite set of uncertain real parameters denoted as~$\zeta_{1},\ldots,\zeta_{n}$.
We collect these sources of uncertainty into vectors $\boldsymbol{\xi}=(\xi_{1},\ldots,\xi_{m})$ and~$\boldsymbol{\zeta}=(\zeta_{1},\ldots,\zeta_{n})$, which we model as random variables that are defined on a probability triple $(\Theta,\mathcal{T},P)$ and which take values in~$\real^{m}$ and $\real^{n}$, respectively.
Throughout this paper, we refer to~$\boldsymbol{\xi}$ and~$\boldsymbol{\zeta}$ as \textit{input random variables}.
Lastly, we assume that the input random variables are statistically independent of one another, that is, we assume that~$\boldsymbol{\xi}$ and~$\boldsymbol{\zeta}$ are statistically independent random variables.

The stochastic coupled model in~(\ref{eq:coupling12}) is a general bidirectionally coupled model. 
The \textit{solution random variable}~$\boldsymbol{u}$ of the first subproblem depends on the solution random variable~$\boldsymbol{v}$ of the second subproblem via the \textit{coupling random variable}~$\boldsymbol{x}$.
Likewise, the solution random variable~$\boldsymbol{v}$ depends on the solution random variable~$\boldsymbol{u}$ via the coupling random variable~$\boldsymbol{y}$.

Thus, to solve this stochastic coupled model, it is necessary to find the random variables~$\boldsymbol{u}$ and~$\boldsymbol{v}$ defined on~$(\Theta,\mathcal{T},P)$ and with values in~$\real^{r}$ and~$\real^{s}$ such that~(\ref{eq:coupling12}) is satisfied under the assumption that the stochastic coupled model is well posed, that is, the model admits a unique and stable solution.

Finally, it should be stressed that the model problem given by~(\ref{eq:coupling12}) is identical to the one given in~\citep{arnst2011a} and~\citep{arnst2011b}, except for the fact that~$\boldsymbol{\xi}$ and~$\boldsymbol{\zeta}$ are required to be statistically independent here; these input random variables were allowed to be statistically dependent in~\citep{arnst2011a} and~\citep{arnst2011b}. 

\subsection{Partitioned iterative solution}
Because a coupled model usually defines its response only in an implicit manner, the numerical solution to a coupled model is typically obtained by using an iterative method. 
Here, we assume that iterative methods and associated solvers already exist for each subproblem.
Therefore, to solve the coupled model, we consider a \textit{partitioned} iterative method that reuses these iterative methods as steps in a global iterative method built around them.
Let us assume that each of the iterative methods is based on the reformulation of the associated subproblem as a fixed-point problem as follows:
\begin{equation}
\label{eq:couplingA12}\begin{aligned}
&\boldsymbol{u}=\boldsymbol{a}(\boldsymbol{u},\boldsymbol{x},\boldsymbol{\xi}),&&\boldsymbol{y}=\boldsymbol{h}(\boldsymbol{u},\boldsymbol{\xi}),&&\boldsymbol{a}:\real^{r}\times\real^{s_{0}}\times\real^{m}\rightarrow\real^{r},&&\boldsymbol{h}:\real^{r}\times\real^{m}\rightarrow\real^{r_{0}},\\
&\boldsymbol{v}=\boldsymbol{b}(\boldsymbol{y},\boldsymbol{v},\boldsymbol{\zeta}),&&\boldsymbol{x}=\boldsymbol{k}(\boldsymbol{v},\boldsymbol{\zeta}),&&\boldsymbol{b}:\real^{r_{0}}\times\real^{s}\times\real^{n}\rightarrow\real^{s},&&\boldsymbol{k}:\real^{s}\times\real^{n}\rightarrow\real^{s_{0}}.
\end{aligned}
\end{equation}
We note that these equations can be obtained by setting~$\boldsymbol{a}(\boldsymbol{u},\boldsymbol{v},\boldsymbol{\xi})=\boldsymbol{u}-\boldsymbol{f}(\boldsymbol{u},\boldsymbol{v},\boldsymbol{\xi})$ and~$\boldsymbol{b}(\boldsymbol{u},\boldsymbol{v},\boldsymbol{\zeta})=\boldsymbol{v}-\boldsymbol{g}(\boldsymbol{u},\boldsymbol{v},\boldsymbol{\zeta})$, but alternative reformulations, such as those involving direct solutions of the subproblems or of their linear approximations, are often better adapted.
We then consider the solution of the stochastic coupled model by a \textit{Gauss-Seidel} iterative method, using suitable initial values~$\boldsymbol{u}^{0}$, $\boldsymbol{v}^{0}$, and~$\boldsymbol{x}^{0}=\boldsymbol{k}(\boldsymbol{v}^{0},\boldsymbol{\zeta})$, as follows:
\begin{equation}
\label{eq:couplingAGS12}\begin{aligned}
&\boldsymbol{u}^{\ell}=\boldsymbol{a}\big(\boldsymbol{u}^{\ell-1},\boldsymbol{x}^{\ell-1},\boldsymbol{\xi}\big),&&\quad\quad\quad\boldsymbol{y}^{\ell}=\boldsymbol{h}(\boldsymbol{u}^{\ell},\boldsymbol{\xi}),\\
&\boldsymbol{v}^{\ell}=\boldsymbol{b}\big(\boldsymbol{y}^{\ell},\boldsymbol{v}^{\ell-1},\boldsymbol{\zeta}\big),&&\quad\quad\quad\boldsymbol{x}^{\ell}=\boldsymbol{k}(\boldsymbol{v}^{\ell},\boldsymbol{\zeta}).
\end{aligned}
\end{equation}
This is not the only partitioned iterative method available; however, for simplicity, we employ only this method in this work.
Although we implement the proposed methodology by using the Gauss-Seidel iterative method, we note that the proposed methodology can be readily used with other iterative methods, such as Jacobi, relaxation, and Newton methods.

\subsection{Dimension reduction}\label{sec:dimensionreduction}
We believe that exchanged information often resides in a considerably lower dimensional space than the sources of uncertainty themselves.
In~\citep{arnst2011a} and~\citep{arnst2011b}, we had thus investigated the effectiveness of dimension-reduction techniques for the representation of the exchanged information; specifically, rather than the exchanging of the coupling random variables $\boldsymbol{x}^{\ell}$ and $\boldsymbol{y}^{\ell}$ and solution random variables $\boldsymbol{u}^{\ell}$ and $\boldsymbol{v}^{\ell}$ in their original form, we had proposed the approximation of these random variables by using a truncated KL decomposition as they pass from subproblem to subproblem and from iteration to iteration.
Here, we present an extension of this previous work.
For stochastic coupled problems whose subproblems have only input random variables that are statistically independent of one another, we present a characterization of the exchanged information by using an alternative dimension-reduction technique, namely, by using a \textit{reduced chaos expansion with random coefficients}.
This expansion has recently been proposed in~\citep{soize2009}, and in the current context, it can be applied to solve the stochastic coupled model as follows:
\begin{equation}
\label{eq:couplingAGSred12}\begin{aligned}
&\hat{\boldsymbol{u}}{}^{\ell}=\boldsymbol{a}\big(\hat{\boldsymbol{u}}{}^{\ell-1,e},\hat{\boldsymbol{x}}{}^{\ell-1,e},\boldsymbol{\xi}\big),&&\quad\quad\quad\hat{\boldsymbol{y}}{}^{\ell}=\boldsymbol{h}(\hat{\boldsymbol{u}}{}^{\ell},\boldsymbol{\xi}),\\
&\hat{\boldsymbol{v}}{}^{\ell}=\boldsymbol{b}\big(\hat{\boldsymbol{y}}{}^{\ell,d},\hat{\boldsymbol{v}}^{\ell-1,d},\boldsymbol{\zeta}\big),&&\quad\quad\quad\hat{\boldsymbol{x}}{}^{\ell}=\boldsymbol{k}(\hat{\boldsymbol{v}}{}^{\ell},\boldsymbol{\zeta}),
\end{aligned}
\end{equation}
where $\boldsymbol{q}^{\ell,d}=(\hat{\boldsymbol{y}}{}^{\ell,d},\hat{\boldsymbol{v}}{}^{\ell-1,d})$ and $\boldsymbol{r}^{\ell,e}=(\hat{\boldsymbol{u}}{}^{\ell-1,e},\hat{\boldsymbol{x}}{}^{\ell-1,e})$ are reduced chaos expansions with random coefficients of $\boldsymbol{q}^{\ell}=(\hat{\boldsymbol{y}}{}^{\ell},\hat{\boldsymbol{v}}{}^{\ell-1})$ and $\boldsymbol{r}^{\ell}=(\hat{\boldsymbol{u}}{}^{\ell-1},\hat{\boldsymbol{x}}{}^{\ell-1})$, respectively, which read as
\begin{equation}
\label{eq:KLqqqq}\begin{aligned}
&\boldsymbol{q}^{\ell,d}(\boldsymbol{\xi},\boldsymbol{\zeta})=\overline{\boldsymbol{q}}{}^{\ell}(\boldsymbol{\zeta})+\sum_{j=1}^{d}\sqrt{\lambda_{j}^{\ell}}\,\eta_{j}^{\ell}(\boldsymbol{\xi})\,\boldsymbol{\phi}^{j,\ell}(\boldsymbol{\zeta}),\\
&\boldsymbol{r}^{\ell,e}(\boldsymbol{\xi},\boldsymbol{\zeta})=\overline{\boldsymbol{r}}{}^{\ell}(\boldsymbol{\xi})+\sum_{j=1}^{e}\sqrt{\kappa_{j}^{\ell}}\,\iota_{j}^{\ell}(\boldsymbol{\zeta})\,\boldsymbol{\psi}^{j,\ell}(\boldsymbol{\xi}).
\end{aligned}
\end{equation}
The reduced chaos expansion with random coefficients is discussed in detail later.
However, here, it is to be noted that these expansions take the form of KL decompositions with random basis vectors.
We say that these expansions maintain \textit{segregation} between the input random variables because each expansion captures one input random variable in the usual manner using the reduced random variables, while capturing the other input random variable, statistically independent from the first, by constructing the basis vectors as random vectors.

Because truncations often result in approximation errors, we use a hat superscript to distinguish between successive approximations determined by~(\ref{eq:couplingAGS12}) and~(\ref{eq:couplingAGSred12}).

We note that~(\ref{eq:KLqqqq}) provides a combined reduced-dimensional representation of~$\hat{\boldsymbol{y}}{}^{\ell}$ and~$\hat{\boldsymbol{v}}{}^{\ell-1}$ in terms of the \textit{reduced random variables}~$\boldsymbol{\eta}^{\ell}=(\eta_{1}^{\ell},\ldots,\eta_{d}^{\ell})$, and a combined reduced-dimensional representation of $\hat{\boldsymbol{u}}{}^{\ell-1}$ and~$\hat{\boldsymbol{x}}{}^{\ell-1}$ in terms of the reduced random variables~$\boldsymbol{\iota}^{\ell}=(\iota_{1}^{\ell},\ldots,\iota_{e}^{\ell})$.
However, these are not the only dimension reductions that can be considered.  
The proposed methodology can be readily used with variants of this dimension reduction, such as variants involving separate reduced-dimensional representations of the coupling and solution random variables, with each representation having its own reduced random variables.  

Further, although our notations do not express the potential dependence of $d$ and $e$ on $\ell$, we note that the reduced dimensions may depend on the number of iterations. 


\subsection{Measure transformation}\label{sec:measuretransformation}
The successive approximations in~(\ref{eq:couplingAGS12}) determined by the iterative method that does not involve dimension reduction can be constructed as random variables of the following form:
\begin{equation}
\label{eq:usefulnessbefore12}
\begin{aligned}
&\boldsymbol{u}{}^{\ell}(\theta)\equiv{\boldsymbol{u}}{}^{\ell}\big(\boldsymbol{\xi}(\theta),\boldsymbol{\zeta}(\theta)\big),\\
&\boldsymbol{v}{}^{\ell}(\theta)\equiv{\boldsymbol{v}}{}^{\ell}\big(\boldsymbol{\xi}(\theta),\boldsymbol{\zeta}(\theta)\big),
\end{aligned}
\end{equation}
that is, $\boldsymbol{u}{}^{\ell}$ and~$\boldsymbol{v}{}^{\ell}$ can be constructed as transformations of~$\boldsymbol{\xi}=(\xi_{1},\ldots,\xi_{m})$ and~$\boldsymbol{\zeta}=(\zeta_{1},\ldots,\zeta_{n})$.
Thus, $\boldsymbol{u}{}^{\ell}$ and~$\boldsymbol{v}{}^{\ell}$ exist in solution spaces of stochastic dimension $m+n$.
Hence, implementations of~(\ref{eq:couplingAGS12}) using stochastic expansion methods typically involve the approximation of the transformations that map~$\boldsymbol{\xi}$ and~$\boldsymbol{\zeta}$ into~$\boldsymbol{u}{}^{\ell}$ and~$\boldsymbol{v}{}^{\ell}$ by finite-dimensional representations.

Depending on whether the nonintrusive projection, embedded projection, or collocation method is chosen, this approximation requires the construction of orthonormal basis functions, quadrature rules, moment tensors, or a combination of these with respect to the joint probability distribution of~$\boldsymbol{\xi}$ and~$\boldsymbol{\zeta}$.
Owing to the assumption that~$\boldsymbol{\xi}$ and~$\boldsymbol{\zeta}$ are statistically independent, the joint probability distribution of~$\boldsymbol{\xi}$ and~$\boldsymbol{\zeta}$ is the product probability distribution~$P_{\boldsymbol{\xi}}\times P_{\boldsymbol{\zeta}}$, where~$P_{\boldsymbol{\xi}}$ and~$P_{\boldsymbol{\zeta}}$ are the probability distributions of~$\boldsymbol{\xi}$ and~$\boldsymbol{\zeta}$, respectively.  
Thus, orthonormal basis functions and quadrature rules with respect to~$P_{\boldsymbol{\xi}}\times P_{\boldsymbol{\zeta}}$ can be constructed by tensorization from orthonormal basis functions and quadrature rules constructed separately with respect to~$P_{\boldsymbol{\xi}}$ and~$P_{\boldsymbol{\zeta}}$; for details, the reader is referred to~\citep{soize2004} and~\citep{reed1980}.
Hence, implementations of~(\ref{eq:couplingAGS12}) typically involve the approximation of $\boldsymbol{u}{}^{\ell}$ and~$\boldsymbol{v}{}^{\ell}$ as
\begin{equation}
\label{eq:usefulnessbefore12b}
\begin{aligned}
&{\boldsymbol{u}}{}^{\ell,p}=\sum_{|\boldsymbol{\alpha}|+|\boldsymbol{\beta}|=0}^{p}\boldsymbol{u}{}^{\ell}_{\boldsymbol{\alpha}\boldsymbol{\beta}}\,\varphi_{\boldsymbol{\alpha}}(\boldsymbol{\xi})\,\psi_{\boldsymbol{\beta}}(\boldsymbol{\zeta}),&&\quad\quad\quad\boldsymbol{u}{}^{\ell}_{\boldsymbol{\alpha}\boldsymbol{\beta}}\in\real^{r},\\
&{\boldsymbol{v}}{}^{\ell,p}=\sum_{|\boldsymbol{\alpha}|+|\boldsymbol{\beta}|=0}^{p}\boldsymbol{v}{}^{\ell}_{\boldsymbol{\alpha}\boldsymbol{\beta}}\,\varphi_{\boldsymbol{\alpha}}(\boldsymbol{\xi})\,\psi_{\boldsymbol{\beta}}(\boldsymbol{\zeta}),&&\quad\quad\quad\boldsymbol{v}{}^{\ell}_{\boldsymbol{\alpha}\boldsymbol{\beta}}\in\real^{s},
\end{aligned}
\end{equation}
where~$\{\varphi_{\boldsymbol{\alpha}},\boldsymbol{\alpha}\in\integer^{m}\}$ and~$\{\psi_{\boldsymbol{\beta}},\boldsymbol{\beta}\in\integer^{n}\}$ are orthonormal bases for the spaces of~$P_{\boldsymbol{\xi}}$-square-integrable functions from~$\real^{m}$ into~$\real$ and of~$P_{\boldsymbol{\zeta}}$-square-integrable functions from~$\real^{n}$ into~$\real$, respectively.
These orthonormal bases are indexed by the usual multi-indices~$\boldsymbol{\alpha}=(\alpha_{1},\ldots,\alpha_{m})$ in~$\integer^{m}$, with~$|\boldsymbol{\alpha}|=\alpha_{1}+\ldots+\alpha_{m}$, and~$\boldsymbol{\beta}=(\beta_{1},\ldots,\beta_{n})$ in~$\integer^{n}$, with~$|\boldsymbol{\beta}|=\beta_{1}+\ldots+\beta_{n}$. 

In contrast, owing to the dimension reduction, the successive approximations determined by the iterative method~(\ref{eq:couplingAGSred12}) can be constructed as random variables of the following form:
\begin{equation}
\label{eq:usefulness12}\begin{aligned}
&\hat{\boldsymbol{u}}{}^{\ell}(\theta)\equiv\hat{\boldsymbol{u}}{}^{\ell}\Big(\boldsymbol{\xi}(\theta),\boldsymbol{\iota}^{\ell}\big(\boldsymbol{\zeta}(\theta)\big)\Big),\\
&\hat{\boldsymbol{v}}{}^{\ell}(\theta)\equiv\hat{\boldsymbol{v}}{}^{\ell}\Big(\boldsymbol{\eta}^{\ell}\big(\boldsymbol{\xi}(\theta)\big),\boldsymbol{\zeta}(\theta)\Big),
\end{aligned}
\end{equation}
that is, $\hat{\boldsymbol{u}}{}^{\ell}$ can be constructed as a transformation of~$\boldsymbol{\xi}=(\xi_{1},\ldots,\xi_{m})$ and~$\boldsymbol{\iota}^{\ell}=(\iota_{1}^{\ell},\ldots,\iota_{e}^{\ell})$ and~$\hat{\boldsymbol{v}}{}^{\ell}$ can be constructed as a transformation of~$\boldsymbol{\eta}^{\ell}=(\eta_{1}^{\ell},\ldots,\eta_{d}^{\ell})$ and~$\boldsymbol{\zeta}=(\zeta_{1},\ldots,\zeta_{n})$.
Thus, $\hat{\boldsymbol{u}}{}^{\ell}$ and $\hat{\boldsymbol{v}}{}^{\ell}$ exist in solution spaces of stochastic dimensions~$m+e$ and~$d+n$, respectively.

Because~$\boldsymbol{\xi}$ and~$\boldsymbol{\zeta}$ are statistically independent and because the reduced chaos expansions with random coefficients given by~(\ref{eq:KLqqqq}) imply that~$\boldsymbol{\iota}^{\ell}$ is a transformation of only~$\boldsymbol{\zeta}$ and $\boldsymbol{\eta}^{\ell}$ is a transformation of only $\boldsymbol{\xi}$, the random variables~$\boldsymbol{\xi}$ and~$\boldsymbol{\iota}^{\ell}$ and, likewise, the random variables~$\boldsymbol{\eta}^{\ell}$ and~$\boldsymbol{\zeta}$ are statistically independent of one another.
Thus, orthonormal basis functions and quadrature rules with respect to the joint probability distribution~$P_{\boldsymbol{\xi}}\times P_{\boldsymbol{\iota}}^{\ell}$ of~$\boldsymbol{\xi}$ and~$\boldsymbol{\iota}^{\ell}$ and those with respect to the joint probability distribution~$P_{\boldsymbol{\eta}}^{\ell}\times P_{\boldsymbol{\zeta}}$ of~$\boldsymbol{\eta}^{\ell}$ and~$\boldsymbol{\zeta}$ can be obtained by tensorization from orthonormal basis functions and quadrature rules constructed separately with respect to the probability distributions~$P_{\boldsymbol{\xi}}$ and $P_{\boldsymbol{\iota}}^{\ell}$ and from those constructed separately with respect to the probability distributions~$P_{\boldsymbol{\zeta}}$ and~$P_{\boldsymbol{\eta}}^{\ell}$, respectively.
Hence, implementations of~(\ref{eq:couplingAGSred12}) can exploit the dimension reduction to approximate $\hat{\boldsymbol{u}}{}^{\ell}$ and~$\hat{\boldsymbol{v}}{}^{\ell}$ as follows:
\begin{equation}
\label{eq:usefulness12b}\begin{aligned}
&\hat{\boldsymbol{u}}{}^{\ell,q}=\sum_{|\boldsymbol{\alpha}|+|\boldsymbol{\varsigma}|=0}^{q}\hat{\boldsymbol{u}}{}^{\ell}_{\boldsymbol{\alpha}\boldsymbol{\varsigma}}\,\varphi_{\boldsymbol{\alpha}}(\boldsymbol{\xi})\,\Upsilon{}_{\boldsymbol{\varsigma}}^{\ell}\big(\boldsymbol{\iota}^{\ell}\big),&&\quad\quad\quad\hat{\boldsymbol{u}}{}^{\ell}_{\boldsymbol{\alpha}\boldsymbol{\varsigma}}\in\real^{r},\\
&\hat{\boldsymbol{v}}{}^{\ell,q}=\sum_{|\boldsymbol{\gamma}|+|\boldsymbol{\beta}|=0}^{q}\hat{\boldsymbol{v}}{}^{\ell}_{\boldsymbol{\gamma}\boldsymbol{\beta}}\,\Gamma{}_{\boldsymbol{\gamma}}^{\ell}\big(\boldsymbol{\eta}^{\ell}\big)\,\psi_{\boldsymbol{\beta}}(\boldsymbol{\zeta}),&&\quad\quad\quad\hat{\boldsymbol{v}}{}^{\ell}_{\boldsymbol{\gamma}\boldsymbol{\beta}}\in\real^{s},
\end{aligned}
\end{equation}
where~$\{\Upsilon{}_{\boldsymbol{\varsigma}}^{\ell},\boldsymbol{\beta}\in\integer^{e}\}$ and~$\{\Gamma{}_{\boldsymbol{\gamma}}^{\ell},\boldsymbol{\gamma}\in\integer^{d}\}$ are orthonormal bases for the spaces of~$P_{\boldsymbol{\iota}}^{\ell}$-square-integrable functions from~$\real^{e}$ into~$\real$ and of~$P_{\boldsymbol{\eta}}^{\ell}$-square-integrable functions from~$\real^{d}$ into~$\real$, respectively, indexed by the usual multi-indices~$\boldsymbol{\varsigma}=(\varsigma_{1},\ldots,\varsigma_{e})$ in~$\integer^{e}$, with~$|\boldsymbol{\varsigma}|=\varsigma_{1}+\ldots+\varsigma_{e}$, and~$\boldsymbol{\gamma}=(\gamma_{1},\ldots,\gamma_{d})$ in~$\integer^{d}$, with~$|\boldsymbol{\gamma}|=\gamma_{1}+\ldots+\gamma_{d}$. 



We say that this approach involves a \textit{measure transformation} because it follows from~(\ref{eq:KLqqqq}) that the probability distributions of the reduced random variables~$\boldsymbol{\iota}^{\ell}$ and~$\boldsymbol{\eta}^{\ell}$ are transformations of the probability distributions of the input random variables~$\boldsymbol{\xi}$ and~$\boldsymbol{\zeta}$, as described later. 

Throughout this work, we employ orthonormal bases that consist of polynomials and we thus refer to them as \textit{polynomial chaos} bases. 

We note that although our notations do not express the potential dependence of $p$ and $q$ on the subproblem or $\ell$, the subsets of the basis functions used to construct the finite-dimensional representations may depend on the subproblem and the number of iterations.


\subsection{Effectiveness of the reduced chaos expansions with random coefficients}
In the reduced chaos expansions with random coefficients, the statistical independence of~$\boldsymbol{\xi}$ and~$\boldsymbol{\iota}^{\ell}$ and that of~$\boldsymbol{\eta}^{\ell}$ and~$\boldsymbol{\zeta}$, which are a consequence of the statistical independence of~$\boldsymbol{\xi}$ and~$\boldsymbol{\zeta}$, are significant because these statistical independences enable polynomial chaos and quadrature rules with respect to the joint probability distribution~$P_{\boldsymbol{\xi}}\times P_{\boldsymbol{\iota}}^{\ell}$ and those with respect to the joint probability distribution~$P_{\boldsymbol{\eta}}^{\ell}\times P_{\boldsymbol{\zeta}}$ to be obtained by tensorization from polynomial chaos and quadrature rules obtained separately with respect to~$P_{\boldsymbol{\xi}}$ and $P_{\boldsymbol{\iota}}^{\ell}$ and from those obtained separately with respect to~$P_{\boldsymbol{\zeta}}$ and~$P_{\boldsymbol{\eta}}^{\ell}$, respectively.

In contrast, even though the methodology presented in~\citep{arnst2011a} and~\citep{arnst2011b} applies to more general stochastic coupled problems whose subproblems need not have only input random variables that are statistically independent of one another, the use of KL decompositions for the representation of the exchanged information does not allow the statistical dependence of~$\boldsymbol{\xi}$ and~$\boldsymbol{\zeta}$ to be exploited for the construction of the requisite polynomial chaos and quadrature rules by using tensorization. 
In general, the use of KL decompositions leads to reduced-dimensional representations of the exchanged information, whose reduced random variables---let us still denote them as~$\boldsymbol{\iota}^{\ell}$ and~$\boldsymbol{\eta}^{\ell}$---are each a transformation of both input random variables~$\boldsymbol{\xi}$ and~$\boldsymbol{\zeta}$, and whose basis vectors are deterministic vectors.
As discussed in~\citep{arnst2011a} and~\citep{arnst2011b}, the implementation of the associated measure transformation leads, in this case, to the construction of polynomial chaos and quadrature rules with respect to the joint probability distributions~$P_{(\boldsymbol{\xi},\boldsymbol{\iota})}^{\ell}$ (of~$\boldsymbol{\xi}$ and~$\boldsymbol{\iota}^{\ell}$) and~$P_{(\boldsymbol{\eta},\boldsymbol{\zeta})}^{\ell}$ (of~$\boldsymbol{\eta}^{\ell}$ and~$\boldsymbol{\zeta}$) because the random variables~$\boldsymbol{\xi}$ and~$\boldsymbol{\iota}^{\ell}$ and the random variables~$\boldsymbol{\eta}^{\ell}$ and~$\boldsymbol{\zeta}$ are then statistically dependent on one another. 

In many applications, the components of the input random variables~$\boldsymbol{\xi}$ and~$\boldsymbol{\zeta}$ are statistically independent ``labeled" random variables; then, polynomial chaos and quadrature rules with respect to~$P_{\boldsymbol{\xi}}$ and~$P_{\boldsymbol{\zeta}}$ can be readily constructed by tensorization of univariate polynomial chaos and quadrature rules read from tables in the literature~\citep{stroud1966,abramowitz1972,reed1980,ghanem2003,xiu2003,soize2004,holtz2010}.
Well-known examples include the use of Hermite and Legendre polynomial chaos and of Gauss-Hermite and Gauss-Legendre quadrature rules with respect to Gaussian and uniform probability distributions. 

However, the reduced random variables~$\boldsymbol{\iota}^{\ell}$ and~$\boldsymbol{\eta}^{\ell}$ associated with reduced chaos expansions with random coefficients or KL decompositions are usually statistically dependent and not ``labeled."
Hence, usually, polynomial chaos and quadrature rules with respect to the probability distributions~$P_{\boldsymbol{\iota}}^{\ell}$ and~$P_{\boldsymbol{\eta}}^{\ell}$ and, likewise, polynomial chaos and quadrature rules with respect to the joint probability distributions~$P_{(\boldsymbol{\xi},\boldsymbol{\iota})}^{\ell}$ and~$P_{(\boldsymbol{\eta},\boldsymbol{\zeta})}^{\ell}$ cannot be read from tables in the literature; they should be computationally constructed.

Now, the computational cost associated with the computation of polynomial chaos and quadrature rules with respect to a given probability distribution can be expected to increase with an increase in the dimension of the space on which this probability distribution is defined.
Thus, for stochastic coupled problems whose subproblems have only input random variables that are statistically independent of one another, the use of reduced chaos expansions with random coefficients for the representation of the exchanged information has the advantage of reducing the computational effort associated with the implementation of the measure transformation since only the computation of polynomial chaos and quadrature rules with respect to the probability distributions of the reduced random variables is necessary.

\subsection{Effectiveness of the proposed methodology}
The main feature of the proposed methodology is that it provides a solution of the subproblems in a reduced-dimensional space when the exchanged information has a low effective stochastic dimension; specifically, a solution in a reduced-dimensional space is obtained when the reduced dimensions can be selected such that $d<m$ and $e<n$, while sufficient accuracy is maintained (refer to~(\ref{eq:usefulnessbefore12}) and~(\ref{eq:usefulness12})).
This is a significant benefit since stochastic expansion methods suffer from a curse of dimensionality: their computational cost increases quickly with an increase in the stochastic dimension.
The proposed methodology breaks the curse of dimensionality by inhibiting the increase in stochastic dimension when information is exchanged.
     
Finally, we note that the proposed methodology can be adapted to meet various requirements of specific applications. 
For instance, one could only use a dimension-reduction technique to represent exchanged random variables that are of low effective stochastic dimension or one could only implement a measure transformation to solve subproblems whose computational cost would thereupon be lowered. 

\section{Reduced chaos expansion with random coefficients}\label{sec:sec3}

\subsection{Problem setting}
Let~$\boldsymbol{\xi}$, $\boldsymbol{\zeta}$, and~$\boldsymbol{q}$ be three random variables defined on the probability triple~$(\Theta,\mathcal{T},P)$ with values in the Euclidean spaces~$\real^{m}$, $\real^{n}$, and~$\real^{w}$, respectively; let their probability distributions be~$P_{\boldsymbol{\xi}}$, $P_{\boldsymbol{\zeta}}$, and~$P_{\boldsymbol{q}}$.
Further, let the random variables~$\boldsymbol{\xi}$ and~$\boldsymbol{\zeta}$ be statistically independent; hence, the joint probability distribution of~$\boldsymbol{\xi}$ and~$\boldsymbol{\zeta}$ is the product probability distribution~$P_{\boldsymbol{\xi}}\times P_{\boldsymbol{\zeta}}$.
Lastly, let the random variable~$\boldsymbol{q}$ be a transformation of~$\boldsymbol{\xi}$ and~$\boldsymbol{\zeta}$ given by
\begin{equation}
\boldsymbol{q}(\theta)\equiv\boldsymbol{q}\big(\boldsymbol{\xi}(\theta),\boldsymbol{\zeta}(\theta)\big)
\end{equation}
and let it be of the second order:
\begin{equation}
\int_{\real^{m}}\int_{\real^{n}}\vectornorm{\boldsymbol{q}(\boldsymbol{\xi},\boldsymbol{\zeta})}^{2}dP_{\boldsymbol{\xi}}dP_{\boldsymbol{\zeta}}<+\infty,
\end{equation}
where~$\vectornorm{\cdot}$ denotes the Euclidean norm.
Under these assumptions, we construct below a reduced chaos expansion with random coefficients of~$\boldsymbol{q}$.

\subsection{Interpretation}
In the context of the proposed methodology, we consider~$\boldsymbol{q}$ as a random variable that passes between subproblems or iterations during the determination of the solution to a stochastic coupled model.
With reference to~(\ref{eq:KLqqqq}), at a specific iteration, $\boldsymbol{q}$ could collect the components of $\hat{\boldsymbol{y}}{}^{\ell}$ and $\hat{\boldsymbol{v}}{}^{\ell-1}$ with $w=r_{0}+s$ or those of $\hat{\boldsymbol{u}}{}^{\ell-1}$ and $\hat{\boldsymbol{x}}{}^{\ell-1}$ with $w=r+s_{0}$.
 
\subsection{Chaos expansion with random coefficients}
Let~$\{\varphi_{\boldsymbol{\alpha}},\boldsymbol{\alpha}\in\integer^{m}\}$ and~$\{\psi_{\boldsymbol{\beta}},\boldsymbol{\beta}\in\integer^{n}\}$ be polynomial chaos bases for the spaces of~$P_{\boldsymbol{\xi}}$-square-integrable functions from~$\real^{m}$ into~$\real$ and of~$P_{\boldsymbol{\zeta}}$-square-integrable functions from~$\real^{n}$ into~$\real$, respectively, indexed by the usual multi-indices~$\boldsymbol{\alpha}=(\alpha_{1},\ldots,\alpha_{m})$ in~$\integer^{m}$ and~$\boldsymbol{\beta}=(\beta_{1},\ldots,\beta_{n})$ in~$\integer^{n}$.
Throughout this paper, we assume that~$\varphi_{\boldsymbol{0}}=1$ and~$\psi_{\boldsymbol{0}}=1$; hence, we have
\begin{equation}
\int_{\real^{m}}\varphi_{\boldsymbol{\alpha}}(\boldsymbol{\xi})dP_{\boldsymbol{\xi}}=0,\quad\quad\boldsymbol{\alpha}\neq\boldsymbol{0},\quad\quad\text{and}\quad\quad\int_{\real^{n}}\psi_{\boldsymbol{\beta}}(\boldsymbol{\zeta})dP_{\boldsymbol{\zeta}}=0,\quad\quad\boldsymbol{\beta}\neq\boldsymbol{0}.\label{eq:mean1}
\end{equation}

Owing to the statistical independence of~$\boldsymbol{\xi}$ and~$\boldsymbol{\zeta}$ and, thus, to the fact that their joint probability distribution is the product distribution~$P_{\boldsymbol{\xi}}\times P_{\boldsymbol{\zeta}}$, the random variable~$\boldsymbol{q}$ has the following chaos expansion related to the tensor product of the aforementioned spaces:
\begin{equation}
\boldsymbol{q}(\boldsymbol{\xi},\boldsymbol{\zeta})=\sum_{\boldsymbol{\alpha}\in\integer^{m}}\sum_{\boldsymbol{\beta}\in\integer^{n}}\boldsymbol{q}_{\boldsymbol{\alpha}\boldsymbol{\beta}}\,\varphi_{\boldsymbol{\alpha}}(\boldsymbol{\xi})\,\psi_{\boldsymbol{\beta}}(\boldsymbol{\zeta}),\quad\quad\boldsymbol{q}_{\boldsymbol{\alpha}\boldsymbol{\beta}}=\int_{\real^{m}}\int_{\real^{n}}\boldsymbol{q}(\boldsymbol{\xi},\boldsymbol{\zeta})\varphi_{\boldsymbol{\alpha}}(\boldsymbol{\xi})\psi_{\boldsymbol{\beta}}(\boldsymbol{\zeta})dP_{\boldsymbol{\xi}}dP_{\boldsymbol{\zeta}};
\end{equation}
here, both sets of sources of uncertainty, that is, $\boldsymbol{\xi}$ and~$\boldsymbol{\zeta}$, are captured in the usual manner using polynomial chaos, and the coefficients~$\boldsymbol{q}_{\boldsymbol{\alpha}\boldsymbol{\beta}}$ are constructed as deterministic vectors in~$\real^{w}$.

Following the definitions given in~\citep{soize2009}, the random variable~$\boldsymbol{q}$ has the following two corresponding \textit{chaos expansions with random coefficients}:
\begin{equation}
\label{eq:chaosrandom12}
\begin{aligned}
\boldsymbol{q}(\boldsymbol{\xi},\boldsymbol{\zeta})&=\sum_{\boldsymbol{\alpha}\in\integer^{m}}\boldsymbol{q}_{\boldsymbol{\alpha}}(\boldsymbol{\zeta})\,\varphi_{\boldsymbol{\alpha}}(\boldsymbol{\xi}),&&\quad\quad\quad\boldsymbol{q}_{\boldsymbol{\alpha}}(\boldsymbol{\zeta})=\sum_{\boldsymbol{\beta}\in\integer^{n}}\boldsymbol{q}_{\boldsymbol{\alpha}\boldsymbol{\beta}}\,\psi_{\boldsymbol{\beta}}(\boldsymbol{\zeta}),\\
\boldsymbol{q}(\boldsymbol{\xi},\boldsymbol{\zeta})&=\sum_{\boldsymbol{\beta}\in\integer^{n}}\boldsymbol{q}_{\boldsymbol{\beta}}(\boldsymbol{\xi})\,\psi_{\boldsymbol{\beta}}(\boldsymbol{\zeta}),&&\quad\quad\quad\boldsymbol{q}_{\boldsymbol{\beta}}(\boldsymbol{\xi})=\sum_{\boldsymbol{\alpha}\in\integer^{m}}\boldsymbol{q}_{\boldsymbol{\alpha}\boldsymbol{\beta}}\,\varphi_{\boldsymbol{\alpha}}(\boldsymbol{\xi});
\end{aligned}
\end{equation}
here, each chaos expansion with random coefficients captures one set of sources of uncertainty in the usual manner using polynomial chaos and the other set by constructing the coefficients as random variables.
It should be emphasized that because~$\boldsymbol{\xi}$ and~$\boldsymbol{\zeta}$ are statistically independent and because the expansions given by~(\ref{eq:chaosrandom12}) imply that the random coefficients~$\boldsymbol{q}_{\boldsymbol{\alpha}}$ are transformations of only~$\boldsymbol{\zeta}$ and the random coefficients~$\boldsymbol{q}_{\boldsymbol{\beta}}$ are transformations of only~$\boldsymbol{\xi}$, the random coefficients~$\boldsymbol{q}_{\boldsymbol{\alpha}}$ are statistically independent of~$\boldsymbol{\xi}$ and, likewise, the random coefficients~$\boldsymbol{q}_{\boldsymbol{\beta}}$ are statistically independent of~$\boldsymbol{\zeta}$. 

Below, we reduce the second of the two chaos expansions with random coefficients in~(\ref{eq:chaosrandom12}) to obtain the corresponding reduced chaos expansion with random coefficients.
We note that the first of the two chaos expansions with random coefficients in~(\ref{eq:chaosrandom12}) can be readily reduced in a similar manner; however, for the sake of brevity, this reduction is not shown.

\subsection{Reduced chaos expansion with random coefficients}
Let~$\boldsymbol{q}$ be approximated by a chaos expansion with random coefficients truncated at total degree~$p$ as follows:
\begin{equation}
\boldsymbol{q}^{p}(\boldsymbol{\xi},\boldsymbol{\zeta})=\sum_{|\boldsymbol{\beta}|=0}^{p}\boldsymbol{q}_{\boldsymbol{\beta}}(\boldsymbol{\xi})\,\psi_{\boldsymbol{\beta}}(\boldsymbol{\zeta}),\label{eq:qp}
\end{equation}
wherein the random coefficients~$\boldsymbol{q}_{\boldsymbol{\beta}}$ are defined by~(\ref{eq:chaosrandom12}).
Then, according to~\citep{soize2009}, an associated reduced chaos expansion with random coefficients is obtained by approximating the random coefficients in this expansion by using a truncated KL decomposition.

Let~$\mu=\text{card}(\{\psi_{\boldsymbol{\beta}},|\boldsymbol{\beta}|\leq p\})=(n+p)!/n!/p!$ denote the number of terms in the chaos expansion with random coefficients given by~(\ref{eq:qp}).
Then, let the random coefficients~$\{\boldsymbol{q}_{\boldsymbol{\beta}},|\boldsymbol{\beta}|\leq p\}$ of~$\boldsymbol{q}^{p}$ be collected in a random variable with values in~$\real^{\mu\times w}$:
\begin{equation}
\begin{bmatrix}
\boldsymbol{q}_{\boldsymbol{\beta}_{1}}(\boldsymbol{\xi})\\
\vdots\\
\boldsymbol{q}_{\boldsymbol{\beta}_{\mu}}(\boldsymbol{\xi})\\
\end{bmatrix}
=
\begin{bmatrix}
\sum_{\boldsymbol{\alpha}\in\integer^{m}}\boldsymbol{q}_{\boldsymbol{\alpha}\boldsymbol{\beta}_{1}}\varphi_{\boldsymbol{\alpha}}(\boldsymbol{\xi})\\
\vdots\\
\sum_{\boldsymbol{\alpha}\in\integer^{m}}\boldsymbol{q}_{\boldsymbol{\alpha}\boldsymbol{\beta}_{\mu}}\varphi_{\boldsymbol{\alpha}}(\boldsymbol{\xi})
\end{bmatrix},\label{eq:vectorq}
\end{equation}
where~$\boldsymbol{\beta}_{1},\ldots,\boldsymbol{\beta}_{\mu}$ are the multi-indices of the random coefficients contained in~$\{\boldsymbol{q}_{\boldsymbol{\beta}},|\boldsymbol{\beta}|\leq p\}$ arranged in a sequence. 
With reference to~(\ref{eq:mean1}), the mean vector of this random variable is the $(\mu\times w)$-dimensional vector that collects the blocks~$\overline{\boldsymbol{q}}{}_{\boldsymbol{\beta}}=\int_{\real^{m}}\boldsymbol{q}_{\boldsymbol{\beta}}(\boldsymbol{\xi})dP_{\boldsymbol{\xi}}$, and it can be expressed as a function of the coefficients~$\boldsymbol{q}_{\boldsymbol{\alpha}\boldsymbol{\beta}}$ as follows:
\begin{equation}
\begin{bmatrix}
\overline{\boldsymbol{q}}{}_{\boldsymbol{\beta}_{1}}\\
\vdots\\
\overline{\boldsymbol{q}}{}_{\boldsymbol{\beta}_{\mu}}\\
\end{bmatrix}
=
\begin{bmatrix}
\boldsymbol{q}_{\boldsymbol{0}\boldsymbol{\beta}_{1}}\\
\vdots\\
\boldsymbol{q}_{\boldsymbol{0}\boldsymbol{\beta}_{\mu}}
\end{bmatrix}.
\end{equation}
Likewise, the covariance matrix is the $(\mu\times w)$-dimensional square matrix that collects the blocks~$\boldsymbol{C}_{\boldsymbol{q}_{\boldsymbol{\beta}}\boldsymbol{q}_{\tilde{\boldsymbol{\beta}}}}=\int_{\real^{m}}\big(\boldsymbol{q}_{\boldsymbol{\beta}}(\boldsymbol{\xi})-\overline{\boldsymbol{q}}{}_{\boldsymbol{\beta}}\big)\big(\boldsymbol{q}_{\tilde{\boldsymbol{\beta}}}(\boldsymbol{\xi})-\overline{\boldsymbol{q}}{}_{\tilde{\boldsymbol{\beta}}}\big)^{\mathrm{T}}dP_{\boldsymbol{\xi}}$, and it can be expressed as a function of the coefficients~$\boldsymbol{q}_{\boldsymbol{\alpha}\boldsymbol{\beta}}$ as follows:
\begin{equation}
\begin{bmatrix}
\boldsymbol{C}_{\boldsymbol{q}_{\boldsymbol{\beta}_{1}}\boldsymbol{q}_{\boldsymbol{\beta}_{1}}} & \ldots & \boldsymbol{C}_{\boldsymbol{q}_{\boldsymbol{\beta}_{1}}\boldsymbol{q}_{\boldsymbol{\beta}_{\mu}}}\\
\vdots & & \vdots\\
\boldsymbol{C}_{\boldsymbol{q}_{\boldsymbol{\beta}_{\mu}}\boldsymbol{q}_{\boldsymbol{\beta}_{1}}} & \ldots & \boldsymbol{C}_{\boldsymbol{q}_{\boldsymbol{\beta}_{\mu}}\boldsymbol{q}_{\boldsymbol{\beta}_{\mu}}}\\
\end{bmatrix}
\hspace{-1mm}=\hspace{-1mm}
\begin{bmatrix}
\sum_{\boldsymbol{\alpha}\in\integer^{m},\boldsymbol{\alpha}\neq\boldsymbol{0}}\boldsymbol{q}_{\boldsymbol{\alpha}\boldsymbol{\beta}_{1}}\boldsymbol{q}_{\boldsymbol{\alpha}\boldsymbol{\beta}_{1}}^{\mathrm{T}} & \ldots & \sum_{\boldsymbol{\alpha}\in\integer^{m},\boldsymbol{\alpha}\neq\boldsymbol{0}}\boldsymbol{q}_{\boldsymbol{\alpha}\boldsymbol{\beta}_{1}}\boldsymbol{q}_{\boldsymbol{\alpha}\boldsymbol{\beta}_{\mu}}^{\mathrm{T}}\\
\vdots & & \vdots\\
\sum_{\boldsymbol{\alpha}\in\integer^{m},\boldsymbol{\alpha}\neq\boldsymbol{0}}\boldsymbol{q}_{\boldsymbol{\alpha}\boldsymbol{\beta}_{\mu}}\boldsymbol{q}_{\boldsymbol{\alpha}\boldsymbol{\beta}_{1}}^{\mathrm{T}} & \ldots & \sum_{\boldsymbol{\alpha}\in\integer^{m},\boldsymbol{\alpha}\neq\boldsymbol{0}}\boldsymbol{q}_{\boldsymbol{\alpha}\boldsymbol{\beta}_{\mu}}\boldsymbol{q}_{\boldsymbol{\alpha}\boldsymbol{\beta}_{\mu}}^{\mathrm{T}}\\
\end{bmatrix}.
\end{equation}

Let $\boldsymbol{W}$ be a $w$-dimensional, square, symmetric, and positive definite matrix, termed the \textit{weighting matrix}.
Because the abovementioned covariance matrix is symmetric and positive semidefinite, the solution of the generalized eigenproblem
\begin{equation}
\begin{bmatrix}
\boldsymbol{W}^{\mathrm{T}}\boldsymbol{C}_{\boldsymbol{q}_{\boldsymbol{\beta}_{1}}\boldsymbol{q}_{\boldsymbol{\beta}_{1}}}\boldsymbol{W} & \ldots & \boldsymbol{W}^{\mathrm{T}}\boldsymbol{C}_{\boldsymbol{q}_{\boldsymbol{\beta}_{1}}\boldsymbol{q}_{\boldsymbol{\beta}_{\mu}}}\boldsymbol{W}\\
\vdots & & \vdots\\
\boldsymbol{W}^{\mathrm{T}}\boldsymbol{C}_{\boldsymbol{q}_{\boldsymbol{\beta}_{\mu}}\boldsymbol{q}_{\boldsymbol{\beta}_{1}}}\boldsymbol{W} & \ldots & \boldsymbol{W}^{\mathrm{T}}\boldsymbol{C}_{\boldsymbol{q}_{\boldsymbol{\beta}_{\mu}}\boldsymbol{q}_{\boldsymbol{\beta}_{\mu}}}\boldsymbol{W}\\
\end{bmatrix}
\begin{bmatrix}
\boldsymbol{\phi}^{j}_{\boldsymbol{\beta}_{1}}\\
\vdots\\
\boldsymbol{\phi}^{j}_{\boldsymbol{\beta}_{\mu}}
\end{bmatrix}
=
\lambda_{j}
\begin{bmatrix}
\boldsymbol{W} \\
& \ddots &\\
& & \boldsymbol{W}
\end{bmatrix}
\begin{bmatrix}
\boldsymbol{\phi}^{j}_{\boldsymbol{\beta}_{1}}\\
\vdots\\
\boldsymbol{\phi}^{j}_{\boldsymbol{\beta}_{\mu}}
\end{bmatrix}
\end{equation}
provides~$\mu\times w$ eigenvalues, $\lambda_{1}\geq\lambda_{2}\geq\ldots\geq\lambda_{\mu\times w}\geq 0$, and~$\mu\times w$ eigenvectors, $\boldsymbol{\phi}^{1},\ldots,\boldsymbol{\phi}^{\mu\times w}$.
These eigenvectors form a $\boldsymbol{W}$-weighted orthonormal basis of $\real^{\mu\times w}$, that is,
\begin{equation}
\begin{bmatrix}
\big(\boldsymbol{\phi}^{i}_{\boldsymbol{\beta}_{1}}\big)^{\mathrm{T}} & \ldots & \big(\boldsymbol{\phi}^{i}_{\boldsymbol{\beta}_{\mu}}\big)^{\mathrm{T}}
\end{bmatrix}
\begin{bmatrix}
\boldsymbol{W} \\
& \ddots &\\
& & \boldsymbol{W}
\end{bmatrix}
\begin{bmatrix}
\boldsymbol{\phi}^{j}_{\boldsymbol{\beta}_{1}}\\
\vdots\\
\boldsymbol{\phi}^{j}_{\boldsymbol{\beta}_{\mu}}
\end{bmatrix}
=\delta_{ij},\label{eq:normalitykl1}
\end{equation}
where $\delta_{ij}$ is the Kronecker delta; $\delta_{ij}$ is equal to 1 if $i=j$, and it is 0 otherwise.
Then, the KL decomposition of the vector given by~(\ref{eq:vectorq}) is obtained as
\begin{equation}
\begin{bmatrix}
\boldsymbol{q}_{\boldsymbol{\beta}_{1}}(\boldsymbol{\xi})\\
\vdots\\
\boldsymbol{q}_{\boldsymbol{\beta}_{\mu}}(\boldsymbol{\xi})\\
\end{bmatrix}
=
\begin{bmatrix}
\overline{\boldsymbol{q}}{}_{\boldsymbol{\beta}_{1}}\\
\vdots\\
\overline{\boldsymbol{q}}{}_{\boldsymbol{\beta}_{\mu}}\\
\end{bmatrix}
+\sum_{j=1}^{\mu\times w}
\sqrt{\lambda_{j}}\eta_{j}(\boldsymbol{\xi})
\begin{bmatrix}
\boldsymbol{\phi}^{j}_{\boldsymbol{\beta}_{1}}\\
\vdots\\
\boldsymbol{\phi}^{j}_{\boldsymbol{\beta}_{\mu}}
\end{bmatrix},\label{eq:KLvectorq}
\end{equation}
where the reduced random variables~$\eta_{j}$ are random variables with values in~$\real$ such that
\begin{equation}
\eta_{j}(\boldsymbol{\xi})=\frac{1}{\sqrt{\lambda_{j}}}
\begin{bmatrix}
\big(\boldsymbol{q}_{\boldsymbol{\beta}_{1}}(\boldsymbol{\xi})-\overline{\boldsymbol{q}}{}_{\boldsymbol{\beta}_{1}}\big)^{\mathrm{T}} & \ldots & \big(\boldsymbol{q}_{\boldsymbol{\beta}_{\mu}}(\boldsymbol{\xi})-\overline{\boldsymbol{q}}{}_{\boldsymbol{\beta}_{\mu}}\big)^{\mathrm{T}}
\end{bmatrix}
\begin{bmatrix}
\boldsymbol{W} \\
& \ddots &\\
& & \boldsymbol{W}
\end{bmatrix}
\begin{bmatrix}
\boldsymbol{\phi}^{j}_{\boldsymbol{\beta}_{1}}\\
\vdots\\
\boldsymbol{\phi}^{j}_{\boldsymbol{\beta}_{\mu}}
\end{bmatrix}.
\end{equation}
These reduced random variables~$\eta_{j}$ are zero-mean and uncorrelated:
\begin{align}
\label{eq:eta1}&\int_{\real^{m}}\eta_{j}(\boldsymbol{\xi})dP_{\boldsymbol{\xi}}=0,\\
\label{eq:eta2}&\int_{\real^{m}}\eta_{i}(\boldsymbol{\xi})\eta_{j}(\boldsymbol{\xi})dP_{\boldsymbol{\xi}}=\delta_{ij}.
\end{align}
By substituting~(\ref{eq:KLvectorq}) in~(\ref{eq:qp}), the following representation of~$\boldsymbol{q}^{p}$ is obtained:
\begin{equation}
\boldsymbol{q}^{p}(\boldsymbol{\xi},\boldsymbol{\zeta})=\overline{\boldsymbol{q}}{}^{p}(\boldsymbol{\zeta})+\sum_{j=1}^{\mu\times w}\sqrt{\lambda_{j}}\eta_{j}(\boldsymbol{\xi})\boldsymbol{\phi}^{j,p}(\boldsymbol{\zeta}).\label{eq:qp2}
\end{equation}
This representation takes the form of a KL decomposition that captures one input random variable, $\boldsymbol{\xi}$, in the usual manner using the reduced random variables and the other input random variable, $\boldsymbol{\zeta}$, by constructing the basis vectors as random vectors, that is,
\begin{align}
\overline{\boldsymbol{q}}{}^{p}(\boldsymbol{\zeta})&=\sum_{|\boldsymbol{\beta}|=0}^{p}\boldsymbol{q}_{\boldsymbol{0}\boldsymbol{\beta}}\psi_{\boldsymbol{\beta}}(\boldsymbol{\zeta}),\\
\boldsymbol{\phi}^{j,p}(\boldsymbol{\zeta})&=\sum_{|\boldsymbol{\beta}|=0}^{p}\boldsymbol{\phi}_{\boldsymbol{\beta}}^{j}\psi_{\boldsymbol{\beta}}(\boldsymbol{\zeta}).
\end{align}
By truncating~(\ref{eq:qp2}) after~$d$ terms, we obtain a \textit{reduced chaos expansion with random coefficients}~$\boldsymbol{q}^{p,d}$ of~$\boldsymbol{q}$ as follows:
\begin{equation}
\boldsymbol{q}^{p,d}(\boldsymbol{\xi},\boldsymbol{\zeta})=\overline{\boldsymbol{q}}{}^{p}(\boldsymbol{\zeta})+\sum_{j=1}^{d}\sqrt{\lambda_{j}}\eta_{j}(\boldsymbol{\xi})\boldsymbol{\phi}^{j,p}(\boldsymbol{\zeta}).\label{eq:qpd}
\end{equation}
Because of the orthonormality properties~(\ref{eq:normalitykl1}) and~(\ref{eq:eta2}), the truncation error satisfies
\begin{equation}
\int_{\real^{m}}\int_{\real^{n}}\vectornorm{\boldsymbol{q}^{p}(\boldsymbol{\xi},\boldsymbol{\zeta})-\boldsymbol{q}^{p,d}(\boldsymbol{\xi},\boldsymbol{\zeta})}_{\boldsymbol{W}}^{2}dP_{\boldsymbol{\xi}}dP_{\boldsymbol{\zeta}}=\sum_{j=d+1}^{\mu\times w}\lambda_{j},\label{eq:KLaccuracy}
\end{equation}
where~$\vectornorm{\cdot}_{\boldsymbol{W}}$ is the~$\boldsymbol{W}$-weighted norm on~$\real^{w}$ such that $\vectornorm{\boldsymbol{p}}_{\boldsymbol{W}}=\sqrt{\boldsymbol{p}^{\mathrm{T}}\boldsymbol{W}\boldsymbol{p}}$ for any~$\boldsymbol{p}$ in~$\real^{w}$.
Thus, the accuracy of the reduced chaos expansion with random coefficients in~(\ref{eq:qpd}) can be improved systematically by increasing the number of terms that are retained.

\subsection{Optimality}

%

Let~$\mathcal{P}^{p}_{n}\otimes\real^{w}$ be the space of random variables with values in~$\real^{w}$ that are representable as a chaos expansion in~$\boldsymbol{\zeta}$ of maximum total degree~$p$.
Then, for any basis~$\{\boldsymbol{e}^{1,p},\ldots,\boldsymbol{e}^{\mu\times w,p}\}$ of~$\mathcal{P}^{p}_{n}\otimes\real^{w}$ with~$\boldsymbol{e}^{j,p}=\sum_{|\boldsymbol{\beta}|=0}^{p}\boldsymbol{e}^{j}_{\boldsymbol{\beta}}\psi_{\boldsymbol{\beta}}$ that is~$\boldsymbol{W}$-weighted orthonormal in that
\begin{equation}
\sum_{|\boldsymbol{\beta}|=0}^{p}\big(\boldsymbol{e}^{i}_{\boldsymbol{\beta}}\big)^{\mathrm{T}}\boldsymbol{W}\boldsymbol{e}^{j}_{\boldsymbol{\beta}}=\delta_{ij},
\end{equation}
the random variable~$\boldsymbol{q}^{p}$ can be expanded as follows:
\begin{equation}
\boldsymbol{q}^{p}(\boldsymbol{\xi},\boldsymbol{\zeta})=\overline{\boldsymbol{q}}{}^{p}(\boldsymbol{\zeta})+\sum_{j=1}^{\mu\times w}\bigg(\sum_{|\boldsymbol{\beta}|=0}^{p}\big(\boldsymbol{q}_{\boldsymbol{\beta}}(\boldsymbol{\xi})-\overline{\boldsymbol{q}}{}_{\boldsymbol{\beta}}\big)^{\mathrm{T}}\boldsymbol{W}\boldsymbol{e}^{j}_{\boldsymbol{\beta}}\bigg)\boldsymbol{e}^{j,p}(\boldsymbol{\zeta}).
\end{equation}
The approximation error~$\epsilon_{d}=\sum_{j=1}^{\mu\times w}\sum_{|\boldsymbol{\beta}|=0}^{p}(\boldsymbol{q}_{\boldsymbol{\beta}}-\overline{\boldsymbol{q}}{}_{\boldsymbol{\beta}})^{\mathrm{T}}\boldsymbol{W}\boldsymbol{e}^{j}_{\boldsymbol{\beta}}\boldsymbol{e}^{j,p}$ introduced because of the truncation of this expansion after~$d$ terms can be readily shown to satisfy
\begin{equation}
\int_{\real^{m}}\int_{\real^{n}}\vectornorm{\epsilon_{d}(\boldsymbol{\xi},\boldsymbol{\zeta})}_{\boldsymbol{W}}^{2}dP_{\boldsymbol{\xi}}dP_{\boldsymbol{\zeta}}=\sum_{j=d+1}^{\mu\times w}\sum_{|\boldsymbol{\beta}|=0}^{p}\sum_{|\tilde{\boldsymbol{\beta}}|=0}^{p}\big(\boldsymbol{e}^{j}_{\boldsymbol{\beta}}\big)^{\mathrm{T}}\boldsymbol{W}^{\mathrm{T}}\boldsymbol{C}_{\boldsymbol{q}_{\boldsymbol{\beta}}\boldsymbol{q}_{\tilde{\boldsymbol{\beta}}}}\boldsymbol{W}\boldsymbol{e}^{j}_{\tilde{\boldsymbol{\beta}}}.\label{eq:optimalitya}
\end{equation}
Expression~(\ref{eq:optimalitya}) indicates that the reduced chaos expansion with random coefficients is optimal in that from among all expansions with~$d$ random basis vectors that are~$\boldsymbol{W}$-weighted orthonormal, it minimizes the~$\boldsymbol{W}$-weighted mean-square norm of the approximation error:
\begin{equation}
\int_{\real^{m}}\int_{\real^{n}}\vectornorm{\boldsymbol{q}^{p}-\boldsymbol{q}^{p,d}}_{\boldsymbol{W}}^{2}dP_{\boldsymbol{\xi}}dP_{\boldsymbol{\zeta}}=\min_{\substack{\{\boldsymbol{e}^{1,p},\ldots,\boldsymbol{e}^{d,p}\}\in\mathcal{P}^{p}_{n}\otimes\real^{w}\\\sum_{|\boldsymbol{\beta}|=0}^{p}(\boldsymbol{e}^{i}_{\boldsymbol{\beta}})^{\mathrm{T}}\boldsymbol{W}\boldsymbol{e}^{j}_{\boldsymbol{\beta}}=\delta_{ij}}}\int_{\real^{m}}\int_{\real^{n}}\vectornorm{\epsilon_{d}}_{\boldsymbol{W}}^{2}dP_{\boldsymbol{\xi}}dP_{\boldsymbol{\zeta}}.\label{eq:optimality}
\end{equation}
Expression~(\ref{eq:optimality}) indicates that the reduced chaos expansion with random coefficients has an optimality property that is similar to that of the KL decomposition; for details on the optimality property of the KL decomposition, the reader is referred to~\citep{ghanem2003}.
However, we note that whereas the KL decomposition achieves optimality among all expansions of a given finite length that capture all sources of uncertainty using the reduced random variables and that construct the basis vectors as deterministic vectors, the reduced chaos expansion with random coefficients achieves optimality among all expansions of a given finite length that capture one set of sources of uncertainty using the reduced random variables and the other set by constructing the basis vectors as random vectors.




\subsection{Concluding remarks}
The reduced chaos expansion with random coefficients described here can be viewed as an adaptation of the one introduced in~\citep{soize2009}; the construction presented here features a weighting matrix, whereas the one given in~\citep{soize2009} does not.
This weighting matrix is particularly useful for the construction of a reduced-dimensional representation of a random variable that solves a space-time discretized stochastic model.
Then, by choosing the weighting matrix as the Gram matrix of the discretization basis, we obtain a reduced-dimensional representation that is consistent with Hilbertian projections in the Hilbert space of random variables to which that random variable belonged prior to the discretization of the stochastic model.
Therefore, the representation satisfies the optimality condition~(\ref{eq:optimality}) in a weighted norm that is consistent with the norm that was relevant to that random variable prior to the discretization of the stochastic model; for more details, the reader is referred to~\citep{arnst2011a}. 
Nevertheless, the reduced chaos expansion with random coefficients introduced in~\citep{soize2009} is recovered by setting the weighting matrix equal to the identity matrix.




\section{Implementation}\label{sec:sec6}

\subsection{Reduced chaos expansion with random coefficients}
In this section, we describe the implementation of the reduced chaos expansion with random coefficients of a random variable that is represented by a chaos expansion involving only polynomial chaos whose total degree does not exceed a given value.
We show how this implementation in turn naturally provides a representation of the reduced random variables by a chaos expansion involving only polynomial chaos whose total degree does not exceed the given value.
Specifically, we adopt the notations used in Secs.~\ref{sec:sec2} and~\ref{sec:sec3}, and we construct a reduced chaos expansion with random coefficients of a second-order random variable~$\boldsymbol{q}^{p}$ with values in~$\real^{w}$ that is represented by the following chaos expansion:
\begin{equation}
\boldsymbol{q}^{p}(\boldsymbol{\xi},\boldsymbol{\zeta})=\sum_{|\boldsymbol{\alpha}|+|\boldsymbol{\beta}|=0}^{p}\boldsymbol{q}_{\boldsymbol{\alpha}\boldsymbol{\beta}}\,\varphi_{\boldsymbol{\alpha}}(\boldsymbol{\xi})\,\psi_{\boldsymbol{\beta}}(\boldsymbol{\zeta}),\quad\quad\quad\boldsymbol{q}_{\boldsymbol{\alpha}\boldsymbol{\beta}}\in \real^{w}.\label{eq:upcekl}
\end{equation}
Then, $\boldsymbol{q}^{p}$ can be expressed by the following chaos expansion with random coefficients:
\begin{equation}
\boldsymbol{q}^{p}(\boldsymbol{\xi},\boldsymbol{\zeta})=\sum_{|\boldsymbol{\beta}|=0}^{p}\boldsymbol{q}_{\boldsymbol{\beta}}^{p-|\boldsymbol{\beta}|}(\boldsymbol{\xi})\,\psi_{\boldsymbol{\beta}}(\boldsymbol{\zeta}),\quad\quad\quad\boldsymbol{q}_{\boldsymbol{\beta}}^{p-|\boldsymbol{\beta}|}(\boldsymbol{\xi})=\sum_{|\boldsymbol{\alpha}|=0}^{p-|\boldsymbol{\beta}|}\boldsymbol{q}_{\boldsymbol{\alpha}\boldsymbol{\beta}}\,\varphi_{\boldsymbol{\alpha}}(\boldsymbol{\xi}).\label{eq:upcekl2}
\end{equation}
Owing to the truncation of the chaos expansion of~$\boldsymbol{q}^{p}$ given by~(\ref{eq:upcekl}) at total degree~$p$, we obtain a representation of each random coefficient~$\boldsymbol{q}_{\boldsymbol{\beta}}^{p-|\boldsymbol{\beta}|}$ that has the form of a chaos expansion involving only polynomial chaos up to total degree~$p-|\boldsymbol{\beta}|$.
Because of the orthonormality of the polynomial chaos~$\varphi_{\boldsymbol{\alpha}}$, the mean vectors and cross-covariance matrices of the random coefficients are immediately obtained as follows:
\begin{equation}
\overline{\boldsymbol{q}}_{\boldsymbol{\beta}}=\boldsymbol{q}_{\boldsymbol{0}\boldsymbol{\beta}}\quad\text{and}\quad\boldsymbol{C}_{\boldsymbol{q}_{\boldsymbol{\beta}}\boldsymbol{q}_{\tilde{\boldsymbol{\beta}}}}=\sum_{|\boldsymbol{\alpha}|=1}^{p-\max(|\boldsymbol{\beta}|,|\tilde{\boldsymbol{\beta}}|)}\boldsymbol{q}_{\boldsymbol{\alpha}\boldsymbol{\beta}}\boldsymbol{q}_{\boldsymbol{\alpha}\tilde{\boldsymbol{\beta}}}^{\mathrm{T}},
\end{equation}   
where owing to the truncation of the chaos expansion of~$\boldsymbol{q}^{p}$ given by~(\ref{eq:upcekl}) at total degree~$p$, cross-covariance matrices~$\boldsymbol{C}_{\boldsymbol{q}_{\boldsymbol{\beta}}\boldsymbol{q}_{\tilde{\boldsymbol{\beta}}}}$ for which~$|\boldsymbol{\beta}|>p-1$, or~$|\tilde{\boldsymbol{\beta}}|>p-1$, or both, vanish.
Then, the solution of the generalized eigenproblem~$\sum_{|\tilde{\boldsymbol{\beta}}|=0}^{p-1}\boldsymbol{W}^{\mathrm{T}}\boldsymbol{C}_{\boldsymbol{q}_{\boldsymbol{\beta}}\boldsymbol{q}_{\tilde{\boldsymbol{\beta}}}}\boldsymbol{W}\boldsymbol{\phi}^{j}_{\tilde{\boldsymbol{\beta}}}=\lambda_{j}\boldsymbol{W}\boldsymbol{\phi}^{j}_{\boldsymbol{\beta}},\;|\boldsymbol{\beta}|\leq p-1$, provides the eigenvalues and associated eigenmodes required to construct a reduced chaos expansion with random coefficients~$\boldsymbol{q}^{p,d}$ of~$\boldsymbol{q}^{p}$ as follows:
\begin{equation}
\boldsymbol{q}^{p,d}(\boldsymbol{\xi},\boldsymbol{\zeta})=\overline{\boldsymbol{q}}{}^{p}(\boldsymbol{\zeta})+\sum_{j=1}^{d}\sqrt{\lambda_{j}}\eta_{j}^{p}(\boldsymbol{\xi})\boldsymbol{\phi}^{j,p-1}(\boldsymbol{\zeta}),
\end{equation}
where owing to the truncation of the chaos expansion of~$\boldsymbol{q}^{p}$ in~(\ref{eq:upcekl}) at total degree~$p$, the basis vectors~$\overline{\boldsymbol{q}}{}^{p}$ and~$\boldsymbol{\phi}^{j,p-1}$ are represented by chaos expansions~$\overline{\boldsymbol{q}}{}^{p}=\sum_{\boldsymbol{\beta}=0}^{p}\overline{\boldsymbol{q}}_{\boldsymbol{\beta}}\psi_{\boldsymbol{\beta}}$ and~$\boldsymbol{\phi}^{j,p-1}=\sum_{|\boldsymbol{\beta}|=0}^{p-1}\boldsymbol{\phi}^{j}_{\boldsymbol{\beta}}\psi_{\boldsymbol{\beta}}$ involving polynomial chaos up to total degrees~$p$ and~$p-1$.
The reduced random variables~$\eta_{j}^{p}$ are random variables with values in~$\real$ such that
\begin{equation}
\eta_{j}^{p}(\boldsymbol{\xi})=\frac{1}{\sqrt{\lambda_{j}}}\sum_{|\boldsymbol{\beta}|=0}^{p-1}\big(\boldsymbol{q}_{\boldsymbol{\beta}}^{p-|\boldsymbol{\beta}|}(\boldsymbol{\xi})-\overline{\boldsymbol{q}}_{\boldsymbol{\beta}}\big)^{\mathrm{T}}\boldsymbol{W}\boldsymbol{\phi}_{\boldsymbol{\beta}}^{j},\label{eq:etaj}
\end{equation}
and they are zero-mean and uncorrelated.
By substituting~(\ref{eq:upcekl2}) in~(\ref{eq:etaj}), we obtain the representation of each reduced random variable as a chaos expansion:
\begin{equation}
\eta_{j}^{p}(\boldsymbol{\xi})=\sum_{|\boldsymbol{\alpha}|=1}^{p}\eta_{j,\boldsymbol{\alpha}}\varphi_{\boldsymbol{\alpha}}(\boldsymbol{\xi})\quad\text{with}\quad\eta_{j,\boldsymbol{\alpha}}=\frac{1}{\sqrt{\lambda_{j}}}\sum_{|\boldsymbol{\beta}|=0}^{p-1}\boldsymbol{q}_{\boldsymbol{\alpha}\boldsymbol{\beta}}^{\mathrm{T}}\boldsymbol{W}\boldsymbol{\phi}_{\boldsymbol{\beta}}^{j}.\label{eq:etajpce}
\end{equation}
Expression~(\ref{eq:etajpce}) indicates that the reduced chaos expansion with random coefficients of a random variable that is represented by a chaos expansion involving only polynomial chaos whose total degree does not exceed a given value naturally provides a complete probabilistic characterization of the reduced random variables by a chaos expansion involving only polynomial chaos whose total degree does not exceed the given value.

We note that while this section provides details on the implementation of the reduced chaos expansion of a random variable that is represented by a chaos expansion that is truncated at a given total degree, the proposed implementation can be readily extended to random variables represented by chaos expansions involving other subsets of polynomial chaos, such as those obtained by full tensorization; for examples, the reader is referred to~\citep{back2011}.

Finally, we note that the methodology proposed in Sec.~\ref{sec:sec2} for obtaining the solution to stochastic coupled problems provides a representation of the solution and coupling variables associated with the subproblems as chaos expansions in terms of combinations of the input random variables~$\boldsymbol{\xi}$ and~$\boldsymbol{\zeta}$ and the reduced random variables~$\boldsymbol{\iota}^{\ell}$ and~$\boldsymbol{\eta}^{\ell}$ of the reduced-dimensional representations of the exchanged information (refer to~(\ref{eq:usefulness12})).
Thus, when these reduced random variables are represented by chaos expansions in terms of the input random variables~$\boldsymbol{\xi}$ and~$\boldsymbol{\zeta}$ themselves (refer to~(\ref{eq:etajpce})), the proposed methodology requires the construction of reduced chaos expansions with random coefficients of random variables that are represented by compositions of chaos expansions.
The implementation of such reduced chaos expansions with random coefficients falls within the scope of the implementation mentioned above because the composition of two chaos expansions can always be written equivalently as a chaos expansion of the form of~(\ref{eq:upcekl}), albeit as one that is truncated at a higher total degree.

\subsection{Construction of the requisite polynomial chaos}
Depending on the stochastic expansion method that is chosen, the approximation of the solutions to the subproblems by chaos expansions in terms of combinations of the input random variables~$\boldsymbol{\xi}$ and~$\boldsymbol{\zeta}$ and the reduced random variables~$\boldsymbol{\iota}^{\ell}$ and~$\boldsymbol{\eta}^{\ell}$, proposed in Sec.~\ref{sec:sec2}, requires the construction of polynomial chaos (discussed in this section) and the construction of moment tensors, quadrature rules, or both (discussed in the next section). 
Here, we construct polynomial chaos with respect to the product probability distribution~$P_{\boldsymbol{\eta}}^{\ell}\times P_{\boldsymbol{\zeta}}$, which is defined on the Euclidean space~$\real^{d+n}$ and which involves the probability distribution~$P_{\boldsymbol{\eta}}^{\ell}$ of the reduced random variable~$\boldsymbol{\eta}^{\ell}$ with values in~$\real^{d}$ and the probability distribution~$P_{\boldsymbol{\zeta}}$ of the input random variable~$\boldsymbol{\zeta}$ with values in~$\real^{n}$.
Clearly, polynomial chaos with respect to~$P_{\boldsymbol{\xi}}\times P_{\boldsymbol{\iota}}^{\ell}$ can be constructed in a similar manner; however, this construction is omitted here for brevity.

\subsubsection{Polynomial chaos with respect to $P_{\boldsymbol{\zeta}}$.}\label{sec:pcPzeta}
In many applications, the components of the input random variables are statistically independent ``labeled" random variables.
Thus, let us assume that the probability distribution of the input random variable~$\boldsymbol{\zeta}$ is the product probability distribution~$P_{\boldsymbol{\zeta}}=P_{\zeta_{1}}\times\ldots\times P_{\zeta_{n}}$ consisting of ``labeled'' univariate probability distributions~$P_{\zeta_{1}},\ldots,P_{\zeta_{n}}$.
Then, for each dimension~$j=1,\ldots,n$, a basis~$\{\psi_{\beta_{j}}^{j},\beta_{j}\in\integer\}$ of polynomial chaos of increasing degree for the space of~$P_{\zeta_{j}}$-square-integrable functions from~$\real$ into~$\real$ can be readily read from tables in the literature, and a corresponding basis~$\{\psi_{\boldsymbol{\beta}},\boldsymbol{\beta}\in\integer^{n}\}$ of polynomial chaos of increasing total degree for the space of~$P_{\boldsymbol{\zeta}}$-square-integrable functions from~$\real^{n}$ into~$\real$, indexed by the multi-indices~$\boldsymbol{\beta}$ in~$\integer^{n}$, can be readily obtained by tensorization by setting~$\psi_{\boldsymbol{\beta}}(\boldsymbol{\zeta})=\psi_{\beta_{1}}^{1}(\zeta_{1})\times\ldots\times\psi_{\beta_{n}}^{n}(\zeta_{n})$; for details, the reader is referred to~\citep{reed1980,ghanem2003,xiu2003,soize2004}.

\subsubsection{Polynomial chaos with respect to $P_{\boldsymbol{\eta}}^{\ell}$.}
The reduced random variables associated with a reduced chaos expansion with random coefficients are usually statistically dependent and not ``labeled."
Hence, polynomial chaos with respect to the probability distribution~$P_{\boldsymbol{\eta}}^{\ell}$ of the reduced random variable~$\boldsymbol{\eta}^{\ell}$ usually cannot be read from tables in the literature, and they should be computationally constructed.
In this paper, we use the procedure presented in~\citep{arnst2011b}, wherein we first arrange the multivariate monomials~$\boldsymbol{\eta}^{\boldsymbol{\gamma}}=\eta_{1}^{\gamma_{1}}\times\ldots\times\eta_{d}^{\gamma_{d}}$ in a sequence of increasing total degree and then orthonormalize this sequence by means of the Gram-Schmidt method; we use the inner product associated with~$P_{\boldsymbol{\eta}}^{\ell}$ to obtain the requisite polynomial chaos.
Let~$\{\Gamma^{\ell}_{\boldsymbol{\gamma}},\boldsymbol{\gamma}\in\integer^{d}\}$ denote the polynomial chaos basis thus obtained for the space of~$P_{\boldsymbol{\eta}}^{\ell}$-square-integrable functions from~$\real^{d}$ into~$\real$, indexed by the multi-indices~$\boldsymbol{\gamma}$ in~$\integer^{d}$.

\subsubsection{Polynomial chaos with respect to $P_{\boldsymbol{\eta}}^{\ell}\times P_{\boldsymbol{\zeta}}$.}
Once the polynomial chaos bases~$\{\Gamma^{\ell}_{\boldsymbol{\gamma}},\boldsymbol{\gamma}\in\integer^{d}\}$ and~$\{\psi_{\boldsymbol{\beta}},\boldsymbol{\beta}\in\integer^{n}\}$ for the spaces of~$P_{\boldsymbol{\eta}}^{\ell}$-square-integrable functions from~$\real^{d}$ into~$\real$ and of~$P_{\boldsymbol{\zeta}}$-square-integrable functions from~$\real^{n}$ into~$\real$, respectively, are available, a corresponding polynomial chaos basis for the space of~$P_{\boldsymbol{\eta}}^{\ell}\times P_{\boldsymbol{\zeta}}$-square-integrable functions from~$\real^{d+n}$ into~$\real$ is readily obtained by tensorization as~$\{\Gamma^{\ell}_{\boldsymbol{\gamma}}\psi_{\boldsymbol{\beta}},\boldsymbol{\gamma}\in\integer^{d},\boldsymbol{\beta}\in\integer^{n}\}$, as discussed in Sec.~\ref{sec:pcPzeta}.
Hence, any~$(P_{\boldsymbol{\eta}}^{\ell}\times P_{\boldsymbol{\zeta}})$-square-integrable function~$f$ from~$\real^{d+n}$ into~$\real$ can be expanded as follows:
\begin{equation}
f(\boldsymbol{\eta}^{\ell},\boldsymbol{\zeta})=\sum_{\boldsymbol{\gamma}\in\integer^{d}}\sum_{\boldsymbol{\beta}\in\integer^{n}}f_{\boldsymbol{\gamma}\boldsymbol{\beta}}\Gamma^{\ell}_{\boldsymbol{\gamma}}(\boldsymbol{\eta}^{\ell})\psi_{\boldsymbol{\beta}}(\boldsymbol{\zeta}),\quad f_{\boldsymbol{\gamma}\boldsymbol{\beta}}=\int_{\real^{d}}\int_{\real^{n}}f(\boldsymbol{\eta}^{\ell},\boldsymbol{\zeta})\Gamma^{\ell}_{\boldsymbol{\gamma}}(\boldsymbol{\eta}^{\ell})\psi_{\boldsymbol{\beta}}(\boldsymbol{\zeta})dP_{\boldsymbol{\eta}}^{\ell}dP_{\boldsymbol{\zeta}}.
\end{equation}
Then, various finite subsets of polynomial chaos can be selected to construct chaos expansions of finite length; for examples, the reader is referred to~\citep{back2011}.
In this work, we use finite subsets of the form~$\{\Gamma^{\ell}_{\boldsymbol{\gamma}}\psi_{\boldsymbol{\beta}},|\boldsymbol{\gamma}|+|\boldsymbol{\beta}|\leq q\}$ to construct chaos expansions of the form of~(\ref{eq:usefulness12b}), involving only polynomials whose total degree does not exceed a given value~$q$.

\subsection{Construction of the requisite quadrature rules}
In this section, we consider the construction of a family of quadrature rules for integration with respect to the product probability distribution~$P_{\boldsymbol{\eta}}^{\ell}\times P_{\boldsymbol{\zeta}}$, which is defined on the Euclidean space~$\real^{d+n}$ and which consists of the probability distribution~$P_{\boldsymbol{\eta}}^{\ell}$ of the reduced random variable~$\boldsymbol{\eta}^{\ell}$ with values in~$\real^{d}$ and the probability distribution~$P_{\boldsymbol{\zeta}}$ of the input random variable~$\boldsymbol{\zeta}$ with values in~$\real^{n}$.
Clearly, a family of quadrature rules with respect to~$P_{\boldsymbol{\xi}}\times P_{\boldsymbol{\iota}}^{\ell}$ can be constructed in a similar manner; however, this construction is omitted here for brevity.

\subsubsection{Quadrature rules with respect to $P_{\boldsymbol{\zeta}}$.}\label{sec:quadrulePzeta}
As in the previous section, let us assume that~$P_{\boldsymbol{\zeta}}$ is the product probability distribution given by~$P_{\boldsymbol{\zeta}}=P_{\zeta_{1}}\times\ldots\times P_{\zeta_{n}}$ and consisting of the ``labeled'' univariate probability distributions~$P_{\zeta_{1}},\ldots,P_{\zeta_{n}}$.
Then, for each dimension~$j=1,\ldots,n$, families of quadrature rules of increasing accuracy for integration with respect to~$P_{\zeta_{j}}$ can be readily constructed and often read from tables in the literature.
While quadrature rules for integration with respect to univariate probability distributions can be obtained following various approaches~\citep{stroud1966,abramowitz1972,holtz2010}, we have used, in this study, Gaussian quadrature rules for integration with respect to univariate probability distributions.
Thus, for each dimension~$j=1,\ldots,n$, let~$Q^{\lambda}_{\zeta_{j}}$ be the level-$\lambda$ Gaussian quadrature rule that allows the integral of any continuous function~$f$ from~$\real$ into~$\real$ with respect to~$P_{\zeta_{j}}$ to be approximated by a weighted sum of integrand evaluations as follows:
\begin{equation}
Q^{\lambda}_{\zeta_{j}}(f)=\sum_{k=1}^{\lambda}f\big(\zeta_{j,k}^{\lambda}\big)v_{j,k}^{\lambda}.
\end{equation}
It is well known that a level-$\lambda$ Gaussian quadrature rule allows integrals of univariate polynomials up to degree~$2\lambda-1$ to be calculated exactly.

Because the components of~$\boldsymbol{\zeta}=(\zeta_{1},\ldots,\zeta_{n})$ are statistically independent, quadrature rules for integration with respect to the multivariate probability distribution~$P_{\boldsymbol{\zeta}}$ can be readily synthesized from the aforementioned quadrature rules for integration with respect to the univariate probability distributions~$P_{\zeta_{1}},\ldots, P_{\zeta_{n}}$ by tensorization.
While various approaches are available to obtain such tensorized quadrature rules \citep{holtz2010}, in this work, we use fully tensorized quadrature rules if the stochastic dimension is low, say~$n\leq 5$, and sparse-grid quadrature rules otherwise.
A level-$\lambda$ fully tensorized quadrature rule that allows the integral of any continuous function~$f$ from~$\real^{n}$ into~$\real$ with respect to~$P_{\boldsymbol{\zeta}}$ to be approximated by a weighted sum of integrand evaluations is obtained as follows:
\begin{equation}
\big(Q_{\zeta_{1}}^{\lambda}\otimes\ldots\otimes Q_{\zeta_{n}}^{\lambda}\big)(f)=\sum_{k_{1}=1}^{\lambda}\ldots\sum_{k_{n}=1}^{\lambda}f(\zeta_{1,k_{1}}^{\lambda},\ldots,\zeta_{n,k_{n}}^{\lambda})v_{1,k_{1}}^{\lambda}\times\ldots\times v_{1,k_{n}}^{\lambda}.\label{eq:fullytensorized}
\end{equation}
Further, a level-$\lambda$ sparse-grid quadrature rule is obtained as follows:
\begin{equation}
Q^{\lambda}_{\boldsymbol{\zeta}}(f)=\sum_{\lambda\leq|\boldsymbol{l}|\leq\lambda+n-1}(-1)^{\lambda+(n-1)-|\boldsymbol{l}|}\begin{pmatrix}n-1\\\lambda+(n-1)-|\boldsymbol{l}|\end{pmatrix}\big(Q^{\lambda_{l_{1}}}_{\zeta_{1}}\otimes\ldots\otimes Q^{\lambda_{l_{n}}}_{\zeta_{n}}\big)(f),\label{eq:smolyak}
\end{equation}
where~$Q^{\lambda_{l_{1}}}_{\zeta_{1}}\otimes\ldots\otimes Q^{\lambda_{l_{n}}}_{\zeta_{n}})(f)=\sum_{k_{1}=1}^{\lambda_{l_{1}}}\ldots\sum_{k_{n}=1}^{\lambda_{l_{n}}}f(\zeta_{1,k_{1}}^{\lambda_{l_{1}}},\ldots,\zeta_{n,k_{n}}^{\lambda_{l_{n}}})v_{1,k_{1}}^{\lambda_{l_{1}}}\times\ldots\times v_{1,k_{n}}^{\lambda_{l_{n}}}$.
Various approaches are available to select~$\lambda_{l_{j}}$ as a function of~$l_{j}$ in~(\ref{eq:smolyak}), for example, the classical Smolyak approach with~$\lambda_{l_{j}}=2^{l_{j}}-1$ and the so-called slowly increasing approach with~$\lambda_{l_{j}}=l_{j}$.
In this work, we use the Smolyak approach.

By the argument given in~\citep{novak1999}, it can be shown that if each univariate quadrature rule of level~$\lambda_{l_{j}}$ allows integrals of univariate polynomials of at least up to a degree of~$2l_{j}-1$ to be calculated exactly, then the sparse-grid quadrature rule of level~$\lambda$ given by~(\ref{eq:smolyak}) allows integrals of multivariate polynomials of at least up to a total degree of~$2\lambda-1$ to be calculated exactly.

\subsubsection{Quadrature rules with respect to $P_{\boldsymbol{\eta}}^{\ell}$.}
The reduced random variables associated with a reduced chaos expansion with random coefficients are usually statistically dependent and not ``labeled."
Hence, usually, quadrature rules for integration with respect to the probability distribution~$P_{\boldsymbol{\eta}}^{\ell}$ of the reduced random variable~$\boldsymbol{\eta}^{\ell}$ cannot be read from tables in the literature, and they should be computationally constructed.

In this work, we use the procedure proposed and described in detail in~\citep{arnst2011b}.
This procedure is applicable, provided that the reduced random variable is itself characterized by a chaos expansion in another random variable whose probability distribution is such that a quadrature rule for integration with respect to this probability distribution can be easily constructed, for example,~$\boldsymbol{\eta}^{\ell}\equiv\boldsymbol{\eta}^{\ell,p}(\boldsymbol{\xi})=\sum_{|\boldsymbol{\alpha}|=1}^{p}\boldsymbol{\eta}^{\ell}_{\boldsymbol{\alpha}}\varphi_{\boldsymbol{\alpha}}(\boldsymbol{\xi})$ (refer to~(\ref{eq:etajpce})).
In this procedure, first, a quadrature rule~$Q_{\boldsymbol{\xi}}$ for integration with respect to~$P_{\boldsymbol{\xi}}$ is constructed.
Next, a change of variables is carried out to obtain a quadrature rule~$Q_{\boldsymbol{\xi}}(\cdot\circ\boldsymbol{\eta}^{\ell,p})$ for integration with respect to~$P_{\boldsymbol{\eta}}^{\ell}$.
Finally, a subset selection algorithm is used to obtain an \textit{embedded quadrature rule},
\begin{equation}
Q_{\boldsymbol{\eta}}^{\ell,\lambda}(f)=\sum_{k=1}^{\nu_{\lambda}^{\ell}}f(\boldsymbol{\eta}_{k}^{\ell,\lambda})w_{k}^{\ell,\lambda},
\end{equation}
that uses a subset of the nodes of the quadrature rule~$Q_{\boldsymbol{\xi}}(\cdot\circ\boldsymbol{\eta}^{\ell,p})$ as nodes and that has level~$\lambda$ in that it allows integrals of multivariate polynomials of at least up to a total degree of~$2\lambda-1$ to be calculated exactly.


\subsubsection{Quadrature rules with respect to $P_{\boldsymbol{\eta}}^{\ell}\times P_{\boldsymbol{\zeta}}$.}
Once~$\lambda$-indexed families of quadrature rules~$Q_{\boldsymbol{\eta}}^{\ell,\lambda}$ and~$Q_{\boldsymbol{\zeta}}^{\lambda}$ for integration with respect to~$P_{\boldsymbol{\eta}}^{\ell}$ and~$P_{\boldsymbol{\zeta}}$ are available, corresponding quadrature rules for integration with respect to~$P_{\boldsymbol{\eta}}^{\ell}\times P_{\boldsymbol{\zeta}}$ can be synthesized by tensorization.
As mentioned in Sec.~\ref{sec:quadrulePzeta}, various approaches are available to obtain tensorized quadrature rules.
In this work, we use the sparse-grid construction.
A level-$\lambda$ sparse-grid quadrature rule that allows the integral of any continuous function~$f$ from~$\real^{d+n}$ into~$\real$ with respect to~$P_{\boldsymbol{\eta}}^{\ell}\times P_{\boldsymbol{\zeta}}$ to be approximated by a weighted sum of integrand evaluations is obtained as follows:
\begin{equation}
Q_{(\boldsymbol{\xi},\boldsymbol{\eta})}^{\ell,\lambda}(f)=\sum_{\lambda\leq k+l\leq\lambda+1}(-1)^{\lambda+1-(k+l)}\big(Q_{\boldsymbol{\eta}}^{\ell,\lambda_{k}}\otimes Q_{\boldsymbol{\zeta}}^{\lambda_{l}}\big)(f),\label{eq:sparseqrulee}
\end{equation}
where~$(Q_{\boldsymbol{\eta}}^{\ell,\lambda_{k}}\otimes Q_{\boldsymbol{\zeta}}^{\lambda_{l}})(f)=\sum_{\tilde{k}=1}^{\nu_{\lambda_{k}}^{\ell}}\sum_{\tilde{l}=1}^{\nu_{\lambda_{l}}}f(\boldsymbol{\eta}_{\tilde{k}}^{\ell,\lambda_{k}},\boldsymbol{\zeta}_{\tilde{l}}^{\lambda_{l}})w_{\tilde{k}}^{\ell,\lambda_{k}}\times v_{\tilde{l}}^{\lambda_{l}}$; here, $\{(\boldsymbol{\eta}_{\tilde{k}}^{\ell,\lambda_{k}},w_{\tilde{k}}^{\ell,\lambda_{k}}),1\leq\tilde{k}\leq\nu_{\lambda_{k}}^{\ell}\}$ and~$\{(\boldsymbol{\zeta}_{\tilde{l}}^{\lambda_{l}},v_{\tilde{l}}^{\lambda_{l}}),1\leq\tilde{l}\leq\nu_{\lambda_{l}}\}$ are the nodes and weights of~$Q_{\boldsymbol{\eta}}^{\ell,\lambda_{k}}$ and~$Q_{\boldsymbol{\zeta}}^{\lambda_{l}}$, respectively.
As mentioned in Sec.~\ref{sec:quadrulePzeta}, various approaches are available to choose~$\lambda_{k}$ and~$\lambda_{l}$ as a function of~$k$ and~$l$.
In this work, we use here the slowly increasing approach with~$\lambda_{k}=k$ and~$\lambda_{l}=l$.

Using the argument given in~\citep{novak1999}, it can be shown that if the component quadrature rules of levels~$\lambda_{k}$ and~$\lambda_{l}$ allow integrals of multivariate polynomials of at least up to total degrees of~$2\lambda_{k}-1$ and~$2\lambda_{l}-1$, respectively, to be calculated exactly, then the sparse-grid quadrature rule of level~$\lambda$ given by~(\ref{eq:sparseqrulee}) allows integrals of multivariate polynomials of at least up to a total degree of~$2\lambda-1$ to be calculated exactly.

\section{Realization for a multiphysics problem}\label{sec:sec7}

\subsection{Problem formulation}
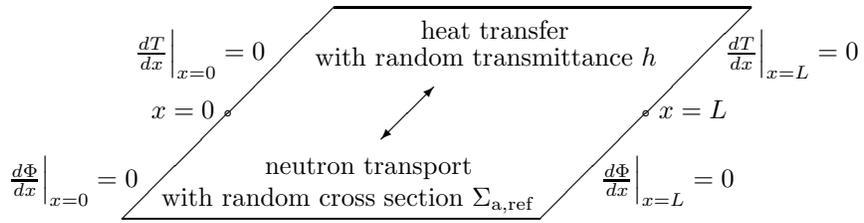
\begin{figure}[htp]
  \begin{center}
    \begin{picture}(250,100)(0,0)
      \put(50,50){\line(1,1){40}}
      \put(50,50){\line(-1,-1){40}}
      \put(208,50){\line(1,1){40}}
      \put(208,50){\line(-1,-1){40}}
      \put(10,10){\line(1,0){158}}
      \put(90,90){\line(1,0){158}}
      \put(123,78){\makebox{heat transfer}}
      \put(53,27){\makebox{$\;\;\;\;$neutron transport}}
      \put(83,68){\makebox{with random transmittance $h$}}
      \put(25,15){\makebox{with random cross section $\Sigma_{\text{a,ref}}$}}
      \put(118,50){\vector(1,1){10}}
      \put(118,50){\vector(-1,-1){10}}
      \put(50,50){\circle{2}}
      \put(208,50){\circle{2}}
      \put(21,48){\makebox{$x=0$}}
      \put(213,48){\makebox{$x=L$}}
      \put(15,70){\makebox{$\frac{dT}{dx}\Big|_{x=0}=0$}}
      \put(238,70){\makebox{$\frac{dT}{dx}\Big|_{x=L}=0$}}
      \put(-33,22){\makebox{$\frac{d\Phi}{dx}\Big|_{x=0}=0$}}
      \put(191,22){\makebox{$\frac{d\Phi}{dx}\Big|_{x=L}=0$}}
    \end{picture}
  \end{center}
  \caption{Schematic representation of the problem.}\label{fig:figure0}
\end{figure}
We consider the stationary transport of neutrons in a one-dimensional reactor with temperature feedback~\citep{lamarsh2002}.
Let the reactor occupy an open interval~$]0,L[$ (Fig.~\ref{fig:figure0}).
The problem involves finding the temperature~$T$ and neutron flux~$\Phi$ such that
\begin{equation}
\label{eq:neutron1c}
\begin{aligned}
&\frac{d}{dx}\left(k\frac{dT}{dx}\right)-h(T-T_{\infty})=-E_{\text{f}}\Sigma_{\text{f}}(T)\Phi,\\
&\frac{d}{dx}\left(D(T)\frac{d\Phi}{dx}\right)-\Big(\Sigma_{\text{a}}(T)-\nu\Sigma_{\text{f}}(T)\Big)\Phi=-s
\end{aligned}
\end{equation}
under homogeneous Neumann boundary conditions.
The first term on the left-hand side of the heat subproblem represents heat conduction, and the second term represents the transmission of heat to the surroundings; further, the right-hand side represents a distributed heat source proportional to the neutron flux.  
The first term on the left-hand side of the neutronics subproblem represents neutron diffusion, and the second term represents the net effect of the absorption and generation of neutrons; further, the right-hand side represents a distributed neutron source.  
The coefficients~$k$ and~$h$ are the heat conductivity and heat transmittance, respectively; the temperature $T_{\infty}$ is the ambient temperature; and $\nu$ and~$E_{\text{f}}$ are the number of neutrons and the energy released per fission reaction, respectively.
The coefficients~$D$, $\Sigma_{\text{a}}$, and~$\Sigma_{\text{f}}$ are the neutron diffusion constant, fission cross section, and absorption cross section, respectively; the dependence of these coefficients on the reactor temperature is as follows:
\begin{equation}
D\big(T(x)\big)=D_{\text{ref}}\sqrt{\frac{T(x)}{T_{\text{ref}}}},\quad\Sigma_{\text{a}}\big(T(x)\big)=\Sigma_{\text{a,ref}}\sqrt{\frac{T_{\text{ref}}}{T(x)}},\quad\Sigma_{\text{f}}\big(T(x)\big)=\Sigma_{\text{f,ref}}\sqrt{\frac{T_{\text{ref}}}{T(x)}}.\label{eq:couplingmechanism}
\end{equation}

\subsection{Deterministic weak formulation}
Let~$H=H^{1}(]0,L[)$ be the space of functions that are sufficiently regular to describe the solutions to the heat and neutronics subproblems.
Then, the weak formulation involves finding~$T$ and~$\Phi$ in~$H$ such that
\begin{equation}
\label{eq:neutron2}\begin{aligned}
&\int_{0}^{L}k\frac{dT}{dx}\frac{dS}{dx}dx+\int_{0}^{L}h(T-T_{\infty})Sdx=\int_{0}^{L}E_{\text{f}}\Sigma_{\text{f}}(T)\Phi Sdx,&&\forall S\in H,\\
&\int_{0}^{L}D(T)\frac{d\Phi}{dx}\frac{d\Psi}{dx}dx+\int_{0}^{L}\Big(\Sigma_{\text{a}}(T)-\nu\Sigma_{\text{f}}(T)\Big)\Phi\Psi dx=\int_{0}^{L}s\Psi dx,&&\forall \Psi\in H.\\
\end{aligned}
\end{equation}

\subsection{Random thermal transmittance and absorption cross section}
We incorporate uncertainties by modeling the thermal transmittance and absorption cross section as random fields~$\{h(x,\boldsymbol{\xi}),0\leq x\leq L\}$ and~$\{\Sigma_{\text{a,ref}}(x,\boldsymbol{\zeta}),0\leq x\leq L\}$ such that
\begin{align}
\label{eq:hN}h(x,\boldsymbol{\xi})&=\overline{h}\bigg(1+\delta_{h}\sum_{j=1}^{m}\sqrt{\lambda_{j}}\,\sqrt{3}\xi_{j}\,\phi^{j}(x)\bigg),\\
\label{eq:Sigmaan}\Sigma_{\text{a,ref}}(x,\boldsymbol{\zeta})&=\overline{\Sigma}_{\text{a,ref}}\bigg(1+\delta_{\Sigma}\sum_{j=1}^{n}\sqrt{\tilde{\lambda}_{j}}\,\sqrt{3}\zeta_{j}\,\tilde{\phi}{}^{j}(x)\bigg),
\end{align}
where the random variables~$\xi_{j}$ and~$\zeta_{j}$ are statistically independent uniform random variables defined on the probability triple~$(\Theta,\mathcal{T},P)$ and with values in $[-1,1]$; thus, the random variables~$\sqrt{3}\xi_{j}$ and~$\sqrt{3}\zeta_{j}$ are uniform random variables with unit standard deviation.
Further, $\lambda_{j}$ and~$\phi^{j}$ ($\tilde{\lambda}_{j}$ and~$\tilde{\phi}{}^{j}$) are the eigenvalues and eigenmodes of the eigenproblem~$\mathcal{C}(\phi^{j})=\lambda_{j}\phi^{j}$ ($\widetilde{\mathcal{C}}(\tilde{\phi}{}^{j})=\tilde{\lambda}_{j}\tilde{\phi}{}^{j}$); $\mathcal{C}$ and~$\widetilde{\mathcal{C}}$ are covariance integral operators with the following kernels:
\begin{align}
\label{eq:covara}C(x,y)&=\frac{4a_{h}^{2}}{\pi^{2}(x-y)^{2}}\sin^{2}\left(\frac{\pi(x-y)}{2a_{h}}\right),\\
\label{eq:covarb}\widetilde{C}(x,y)&=\frac{4a_{\Sigma}^{2}}{\pi^{2}(x-y)^{2}}\sin^{2}\left(\frac{\pi(x-y)}{2a_{\Sigma}}\right).
\end{align}
Here, the parameters~$a_{h}$ and~$a_{\Sigma}$ are the spatial correlation lengths of the thermal transmittance and absorption cross section random fields, respectively.
Clearly, the random fields thus obtained are such that the random variables~$h(x,\boldsymbol{\xi})$ and~$\Sigma_{\text{a,ref}}(x,\boldsymbol{\zeta})$ have mean values~$\overline{h}$ and~$\overline{\Sigma}_{\text{a,ref}}$ and coefficients of variation~$\delta_{h}$ and~$\delta_{\Sigma}$, respectively, at every position~$x$, at least when the approximation errors introduced because of the truncation of the expansions after~$m$ and~$n$ terms, respectively, are not taken into account.  

\subsection{Stochastic weak formulation} 
The weak formulation of the stochastic problem involves finding random variables~$T$ and~$\Phi$ on~$(\Theta,\mathcal{T},P)$ with values in~$H$ such that with $h(\boldsymbol{\xi})\equiv h(\cdot,\boldsymbol{\xi})$ and~$\Sigma_{\text{a}}(T,\boldsymbol{\zeta})\equiv\Sigma_{\text{a,ref}}(\cdot,\boldsymbol{\zeta})\sqrt{{T_{\text{ref}}}/{T}}$,
\begin{equation}
\label{eq:neutron2b}
\begin{aligned}
&\int_{0}^{L}k\frac{dT}{dx}\frac{dS}{dx}dx+\int_{0}^{L}h(\boldsymbol{\xi})(T-T_{\infty})Sdx=\int_{0}^{L}E_{\text{f}}\Sigma_{\text{f}}(T)\Phi Sdx,&&\forall S\in H,\\
&\int_{0}^{L}D(T)\frac{d\Phi}{dx}\frac{d\Psi}{dx}dx+\int_{0}^{L}\Big(\Sigma_{\text{a}}(T,\boldsymbol{\zeta})-\nu\Sigma_{\text{f}}(T)\Big)\Phi\Psi dx=\int_{0}^{L}s\Psi dx,&&\forall \Psi\in H.
\end{aligned}
\end{equation}

\subsection{Discretization of space}
The finite element~(FE) method is used for the discretization of space.
The domain~$[0,L]$ is meshed using $r-1$ elements of equal length.  
Let~$N_{1},\ldots,N_{r}$ be a basis of element-wise linear shape functions such that~$N_{j}$ is equal to~$1$ at the~$j$-th node and~0 at all other nodes. 
Using this basis, the random temperature~$T$ and random neutron flux~$\Phi$ are approximated as follows:
\begin{equation}
\label{eq:Vh1}
\begin{aligned}
T^{r}(x)=\sum_{j=1}^{r}T_{j}N_{j}(x),\quad\quad\quad T_{j}\in\real,\\
\Phi^{r}(x)=\sum_{j=1}^{r}\Phi_{j}N_{j}(x),\quad\quad\quad \Phi_{j}\in\real.\\
\end{aligned}
\end{equation} 
Then, the FE discretization of the stochastic weak formulation in~(\ref{eq:neutron2b}) involves the determination of the random vectors~$\boldsymbol{T}=(T_{1},\ldots,T_{r})$ and~$\boldsymbol{\Phi}=(\Phi_{1},\ldots,\Phi_{r})$ defined on~$(\Theta,\mathcal{T},P)$, with values in~$\real^{r}$, which collect the nodal values of the random temperature and random neutron flux such that
\begin{equation}
\label{eq:discrfE3b}
\begin{aligned}
&[\boldsymbol{K}+\boldsymbol{H}(\boldsymbol{\xi})]\boldsymbol{T}=\boldsymbol{q}(\boldsymbol{\Phi},\boldsymbol{T}),\\
&[\boldsymbol{D}(\boldsymbol{T})+\boldsymbol{M}(\boldsymbol{T},\boldsymbol{\zeta})]\boldsymbol{\Phi}=\boldsymbol{s}.
\end{aligned}
\end{equation}
Here, $\boldsymbol{K}$, $\boldsymbol{H}(\boldsymbol{\xi})$, $\boldsymbol{D}(\boldsymbol{T})$, and $\boldsymbol{M}(\boldsymbol{T},\boldsymbol{\zeta})$ are $r$-dimensional matrices, and $\boldsymbol{q}(\boldsymbol{\Phi},\boldsymbol{T})$ and $\boldsymbol{s}$ are $r$-dimensional vectors such that
\begin{align}
\boldsymbol{S}_{1}^{\mathrm{T}}\boldsymbol{K}\boldsymbol{S}_{2}&=\int_{0}^{L}k\frac{dS^{r}_{1}}{dx}\frac{dS^{r}_{2}}{dx}dx,\\
\boldsymbol{S}_{1}^{\mathrm{T}}\boldsymbol{H}(\boldsymbol{\xi})\boldsymbol{S}_{2}&=\int_{0}^{L}h(\boldsymbol{\xi})S^{r}_{1}S^{r}_{2}dx,\\
\boldsymbol{\Psi}_{1}^{\mathrm{T}}\boldsymbol{D}(\boldsymbol{T})\boldsymbol{\Psi}_{2}&=\int_{0}^{L}D\big(T^{r}\big)\frac{d\Psi^{r}_{1}}{dx}\frac{d\Psi^{r}_{2}}{dx}dx,\\
\boldsymbol{\Psi}_{1}^{\mathrm{T}}\boldsymbol{M}(\boldsymbol{T},\boldsymbol{\zeta})\boldsymbol{\Psi}_{2}&=\int_{0}^{L}\Big(\Sigma_{\text{a}}\big(T^{r},\boldsymbol{\zeta}\big)-\nu\Sigma_{\text{f}}\big(T^{r}\big)\Big)\Psi^{r}_{1}\Psi^{r}_{2}dx,\\
\boldsymbol{S}^{\mathrm{T}}\boldsymbol{q}(\boldsymbol{T},\boldsymbol{\Phi})&=\int_{0}^{L}E_{\text{f}}\Sigma_{\text{f}}\big(T^{r}\big)\Phi^{r}S^{r}dx+\int_{0}^{L}hT_{\infty}S^{r}dx,\\
\boldsymbol{S}^{\mathrm{T}}\boldsymbol{s}&=\int_{0}^{L}s\Psi^{r}dx.
\end{align}

\subsection{Reformulation as a realization of the model problem}
The aforementioned illustration problem can be reformulated as a particular realization of the general model problem introduced in Sec.~\ref{sec:sec2}:
\begin{equation}
\begin{aligned}
&\boldsymbol{T}=\boldsymbol{a}(\boldsymbol{T},\boldsymbol{\Phi},\boldsymbol{\xi}),&&\quad\quad\quad\boldsymbol{a}:\real^{r}\times\real^{r}\times\real^{m}\rightarrow\real^{r},\\
&\boldsymbol{\Phi}=\boldsymbol{b}(\boldsymbol{T},\boldsymbol{\zeta}),&&\quad\quad\quad\boldsymbol{b}:\real^{r}\times\real^{n}\rightarrow\real^{r},
\end{aligned}
\end{equation}
where $\boldsymbol{a}(\boldsymbol{T},\boldsymbol{\Phi},\boldsymbol{\xi})=[\boldsymbol{K}+\boldsymbol{H}(\boldsymbol{\xi})]^{-1}\boldsymbol{q}(\boldsymbol{\Phi},\boldsymbol{T})$ and $\boldsymbol{b}(\boldsymbol{T},\boldsymbol{\zeta})=[\boldsymbol{D}(\boldsymbol{T})+\boldsymbol{M}(\boldsymbol{T},\boldsymbol{\zeta})]^{-1}\boldsymbol{s}$.
This reformulation indicates that the illustration problem is a simplified realization of the model problem.
This is inferred from two features.
First, the neutronics subproblem admits a direct solution that does not require iteration.
Second, the neutronics and heat subproblems are coupled directly through their solution variables rather than intermediate coupling variables.

\subsection{Dimension reduction}
Now, we will demonstrate the proposed methodology by approximating the random temperature by a reduced chaos expansion with random coefficients as it is communicated from the heat subproblem to the neutronics subproblem.
At iteration~$\ell$, let the random temperature be represented by the following chaos expansion:
\begin{equation}
\widehat{\boldsymbol{T}}{}^{\ell,p}(\boldsymbol{\xi},\boldsymbol{\zeta})=\sum_{|\boldsymbol{\alpha}|+|\boldsymbol{\beta}|=0}^{p}\widehat{\boldsymbol{T}}{}^{\ell}_{\boldsymbol{\alpha}\boldsymbol{\beta}}\,\varphi_{\boldsymbol{\alpha}}(\boldsymbol{\xi})\,\psi_{\boldsymbol{\beta}}(\boldsymbol{\zeta}),\quad\quad\quad\widehat{\boldsymbol{T}}{}^{\ell}_{\boldsymbol{\alpha}\boldsymbol{\beta}}\in \real^{r}.\label{eq:upceklb}
\end{equation}
Then, the random temperature has the following chaos expansion with random coefficients:
\begin{equation}
\widehat{\boldsymbol{T}}{}^{\ell,p}(\boldsymbol{\xi},\boldsymbol{\zeta})=\sum_{|\boldsymbol{\beta}|=0}^{p}\widehat{\boldsymbol{T}}{}_{\boldsymbol{\beta}}^{\ell,p-|\boldsymbol{\beta}|}(\boldsymbol{\xi})\,\psi_{\boldsymbol{\beta}}(\boldsymbol{\zeta}),\quad\quad\quad\widehat{\boldsymbol{T}}{}_{\boldsymbol{\beta}}^{\ell,p-|\boldsymbol{\beta}|}(\boldsymbol{\xi})=\sum_{|\boldsymbol{\alpha}|=0}^{p-|\boldsymbol{\beta}|}\widehat{\boldsymbol{T}}{}^{\ell}_{\boldsymbol{\alpha}\boldsymbol{\beta}}\,\varphi_{\boldsymbol{\alpha}}(\boldsymbol{\xi}).\label{eq:upceklb2}
\end{equation}
This expansion involves the representation of each random coefficient as a chaos expansion involving only polynomial chaos up to total degree~$p-|\boldsymbol{\beta}|$.
The mean vectors and cross-covariance matrices of the random coefficients are then given by
\begin{equation}
\overline{\boldsymbol{T}}{}^{\ell}_{\boldsymbol{\beta}}=\widehat{\boldsymbol{T}}{}^{\ell}_{\boldsymbol{0}\boldsymbol{\beta}}\quad\text{and}\quad\boldsymbol{C}{}^{\ell}_{\widehat{\boldsymbol{T}}_{\boldsymbol{\beta}}\widehat{\boldsymbol{T}}_{\tilde{\boldsymbol{\beta}}}}=\sum_{|\boldsymbol{\alpha}|=1}^{p-\max(|\boldsymbol{\beta}|,|\tilde{\boldsymbol{\beta}}|)}\widehat{\boldsymbol{T}}{}^{\ell}_{\boldsymbol{\alpha}\boldsymbol{\beta}}(\widehat{\boldsymbol{T}}{}^{\ell}_{\boldsymbol{\alpha}\tilde{\boldsymbol{\beta}}})^{\mathrm{T}},
\end{equation}   
where cross-covariance matrices~$\boldsymbol{C}_{\boldsymbol{q}_{\boldsymbol{\beta}}\boldsymbol{q}_{\tilde{\boldsymbol{\beta}}}}$ for which~$|\boldsymbol{\beta}|>p-1$, or~$|\tilde{\boldsymbol{\beta}}|>p-1$, or both, vanish.
Let the $r$-dimensional square matrix~$\boldsymbol{W}$ be the Gram matrix of the FE basis, that is,
\begin{equation}
\boldsymbol{W}=\begin{bmatrix}
\langle N_{1},N_{1}\rangle_{H} & \ldots & \langle N_{1},N_{r}\rangle_{H}\\
\vdots & & \vdots\\
\langle N_{r},N_{1}\rangle_{H} & \ldots & \langle N_{r},N_{r}\rangle_{H}
\end{bmatrix},
\end{equation}
where the inner product~$\langle\cdot,\cdot\rangle_{H}$ is such that~$\langle S_{1},S_{2}\rangle_{H}=\int_{0}^{L}S_{1}S_{2}dx+\int_{0}^{L}(dS_{1}/dx)(dS_{2}/dx)dx$ for any pair of functions~$S_{1}$ and~$S_{2}$ in~$H$. 
Then, the solution of the generalized eigenproblem~$\sum_{|\tilde{\boldsymbol{\beta}}|=0}^{p-1}\boldsymbol{W}^{\mathrm{T}}\boldsymbol{C}{}^{\ell}_{\widehat{\boldsymbol{T}}_{\boldsymbol{\beta}}\widehat{\boldsymbol{T}}_{\tilde{\boldsymbol{\beta}}}}\boldsymbol{W}\boldsymbol{\phi}^{j,\ell}_{\tilde{\boldsymbol{\beta}}}=\lambda_{j}^{\ell}\boldsymbol{W}\boldsymbol{\phi}^{j,\ell}_{\boldsymbol{\beta}},\;|\boldsymbol{\beta}|\leq p-1$, provides the eigenvalues and eigenmodes required to construct a reduced chaos expansion with random coefficients,
\begin{equation}
\widehat{\boldsymbol{T}}{}^{\ell,p,d}(\boldsymbol{\xi},\boldsymbol{\zeta})=\overline{\boldsymbol{T}}{}^{\ell,p}(\boldsymbol{\zeta})+\sum_{j=1}^{d}\sqrt{\lambda_{j}^{\ell}}\,\eta_{j}^{\ell,p}(\boldsymbol{\xi})\,\boldsymbol{\phi}^{j,\ell,p-1}(\boldsymbol{\zeta}),\label{eq:Tpdkl}
\end{equation}
where the basis vectors~$\overline{\boldsymbol{T}}{}^{\ell,p}$ and~$\boldsymbol{\phi}^{j,\ell,p-1}$ are represented by the chaos expansions~$\overline{\boldsymbol{T}}{}^{\ell,p}=\sum_{\boldsymbol{\beta}=0}^{p}\overline{\boldsymbol{T}}{}^{\ell}_{\boldsymbol{\beta}}\psi_{\boldsymbol{\beta}}$ and~$\boldsymbol{\phi}^{j,\ell,p-1}=\sum_{|\boldsymbol{\beta}|=0}^{p-1}\boldsymbol{\phi}^{j,\ell}_{\boldsymbol{\beta}}\psi_{\boldsymbol{\beta}}$ involving polynomial chaos up to total degrees~$p$ and~$p-1$.
The reduced random variables~$\eta_{j}^{\ell,p}$ are random variables with values in~$\real$, such that
\begin{equation}
\eta_{j}^{\ell,p}(\boldsymbol{\xi})=\frac{1}{\sqrt{\lambda_{j}^{\ell}}}\sum_{|\boldsymbol{\beta}|=0}^{p-1}\big(\widehat{\boldsymbol{T}}{}_{\boldsymbol{\beta}}^{\ell,p-|\boldsymbol{\beta}|}(\boldsymbol{\xi})-\overline{\boldsymbol{T}}{}^{\ell}_{\boldsymbol{\beta}}\big)^{\mathrm{T}}\boldsymbol{W}\boldsymbol{\phi}_{\boldsymbol{\beta}}^{j,\ell},\label{eq:etajb}
\end{equation}
and they are zero-mean and uncorrelated.
By substituting~(\ref{eq:upceklb2}) in~(\ref{eq:etajb}), we obtain the representation of each reduced random variable as a chaos expansion:
\begin{equation}
\eta_{j}^{\ell,p}(\boldsymbol{\xi})=\sum_{|\boldsymbol{\alpha}|=1}^{p}\eta^{\ell}_{j,\boldsymbol{\alpha}}\varphi_{\boldsymbol{\alpha}}(\boldsymbol{\xi})\quad\text{with}\quad\eta^{\ell}_{j,\boldsymbol{\alpha}}=\frac{1}{\sqrt{\lambda_{j}^{\ell}}}\sum_{|\boldsymbol{\beta}|=0}^{p-1}(\widehat{\boldsymbol{T}}{}^{\ell}_{\boldsymbol{\alpha}\boldsymbol{\beta}})^{\mathrm{T}}\boldsymbol{W}\boldsymbol{\phi}_{\boldsymbol{\beta}}^{j,\ell},\label{eq:etajpce2}
\end{equation}
thus completely characterizing the reduced random variables as a chaos expansion.

We note that the random neutron flux, in principle, could also be represented by a reduced chaos expansion with random coefficients as it passes from the neutronics subproblem to the heat subproblem.
However, this extension is omitted here.

The present implementation can be expected to be well adapted to problems wherein the stochastic dimension of the random thermal transmittance field is moderate or high, say $m> 5$, but the stochastic dimension of the random absorption cross section field is low, say $n\leq 5$.
Then, the reduction of the random temperature has the potential to lower the number of sources of uncertainty that enter the neutronics subproblem and to thus provide a more efficient solution of the neutronics subproblem in a reduced-dimensional space.
In contrast, a reduction of the random neutron flux would not have the potential to significantly lower the number of sources of uncertainty that enter the heat subproblem; therefore, it is not implemented.  

\subsection{Measure transformation}\label{sec:measuretransformationillustration}
While~$\boldsymbol{\xi}$ and~$\boldsymbol{\zeta}$ are necessarily the sources of uncertainty that enter the heat subproblem, the aforementioned representation of the random temperature by a reduced chaos expansion with random coefficients allows~$\boldsymbol{\eta}^{\ell,p}$ and~$\boldsymbol{\zeta}$ to be construed as the sources of uncertainty that enter the neutronics subproblem.
Then, the proposed methodology leads to the approximation of the random temperature and neutron flux by chaos expansions as follows:
\begin{equation}
\label{eq:PCTphi}\begin{aligned}
\widehat{\boldsymbol{T}}{}^{\ell,p}(\boldsymbol{\xi},\boldsymbol{\zeta})&=\sum_{|\boldsymbol{\alpha}|+|\boldsymbol{\beta}|=0}^{p}\widehat{\boldsymbol{T}}{}^{\ell}_{\boldsymbol{\alpha}\boldsymbol{\beta}}\varphi_{\boldsymbol{\alpha}}(\boldsymbol{\xi})\psi_{\boldsymbol{\beta}}(\boldsymbol{\zeta}),&&\quad\quad\quad\widehat{\boldsymbol{T}}{}^{\ell}_{\boldsymbol{\alpha}\boldsymbol{\beta}}\in\real^{r},\\
\widehat{\boldsymbol{\Phi}}{}^{\ell,q}(\boldsymbol{\eta}^{\ell,p},\boldsymbol{\zeta})&=\sum_{|\boldsymbol{\gamma}|+|\boldsymbol{\beta}|=0}^{q}\widehat{\boldsymbol{\Phi}}{}^{\ell}_{\boldsymbol{\gamma}\boldsymbol{\beta}}\Gamma_{\boldsymbol{\gamma}}^{\ell}(\boldsymbol{\eta}^{\ell,p})\psi_{\boldsymbol{\beta}}(\boldsymbol{\zeta}),&&\quad\quad\quad\widehat{\boldsymbol{\Phi}}{}^{\ell}_{\boldsymbol{\beta}\boldsymbol{\gamma}}\in\real^{r},
\end{aligned}
\end{equation}
that is, we approximate the random temperature by a chaos expansion in~$\boldsymbol{\xi}$ and~$\boldsymbol{\zeta}$ and the random neutron flux by a chaos expansion in~$\boldsymbol{\eta}^{\ell,p}$ and~$\boldsymbol{\zeta}$.

We select~$\{\varphi_{\boldsymbol{\alpha}},\boldsymbol{\alpha}\in\integer^{m}\}$ and~$\{\psi_{\boldsymbol{\beta}},\boldsymbol{\beta}\in\integer^{m}\}$ as normalized Legendre polynomials.  
Following the approach given in Sec.~\ref{sec:sec6}, we construct the polynomial chaos~$\{\Gamma^{\ell}_{\boldsymbol{\gamma}},\boldsymbol{\gamma}\in\integer^{d}\}$ at each iteration using the method given in~\citep{arnst2011b}.
Further, we select the quadrature rule for integration with respect to the joint probability distribution of~$\boldsymbol{\xi}$ and~$\boldsymbol{\zeta}$, whose nodes and weights we denote as~$\{(\boldsymbol{\xi}_{k},\boldsymbol{\zeta}_{k},v_{k}),\;1\leq k\leq N\}$, to be a sparse-grid Gauss-Legendre quadrature rule of dimension~$m+n$ and level $p+1$.
Following the approach given in Sec.~\ref{sec:sec6}, we construct the quadrature rule for integration with respect to the joint probability distribution of~$\boldsymbol{\eta}^{\ell,p}$ and~$\boldsymbol{\zeta}$, whose nodes and weights we denote as $\{(\boldsymbol{\eta}_{k}^{\ell},\boldsymbol{\zeta}_{k}^{\ell},w_{k}^{\ell}),\;1\leq k\leq \nu^{\ell}\}$, at each iteration as a sparse-grid quadrature rule of dimension~$d+n$ and level~$q+2$.
This rule is synthesized from the family of fully tensorized Gauss-Legendre quadrature rules for integration with respect to the probability distribution of $\boldsymbol{\zeta}$ and the family of embedded quadrature rules for integration with respect to the probability distribution of~$\boldsymbol{\eta}^{\ell,p}$ obtained using the method given in~\citep{arnst2011b}. 

\subsection{Selection of the reduced dimension and polynomial degree}
At each iteration, we select the number of terms retained in~(\ref{eq:Tpdkl}) as the smallest dimension $d$ that satisfies
\small{
\begin{equation}
\sqrt{\int_{\real^{m}}\int_{\real^{n}}\hspace{-1mm}\vectornorm{\widehat{\boldsymbol{T}}{}^{\ell,p}(\boldsymbol{\xi},\boldsymbol{\zeta})-\widehat{\boldsymbol{T}}{}^{\ell,p,d}(\boldsymbol{\xi},\boldsymbol{\zeta})}^{2}_{\boldsymbol{W}}\hspace{-1mm}dP_{\boldsymbol{\xi}}dP_{\boldsymbol{\zeta}}}\leq\epsilon_{1}\sqrt{\int_{\real^{m}}\int_{\real^{n}}\hspace{-1mm}\vectornorm{\widehat{\boldsymbol{T}}{}^{\ell,p}(\boldsymbol{\xi},\boldsymbol{\zeta})}^{2}_{\boldsymbol{W}}\hspace{-1mm}dP_{\boldsymbol{\xi}}dP_{\boldsymbol{\zeta}}},\quad\forall\ell\in\integer,\label{eq:criterion1}
\end{equation}}\normalsize
where~$\epsilon_{1}$ is a prescribed tolerance level.
Further, at each iteration, we truncate the chaos expansion of the random neutron flux in (\ref{eq:PCTphi}) at the smallest total degree $q$ that satisfies

\vspace{-3mm}
\small{
\begin{equation}
\sqrt{\int_{\real^{d}}\int_{\real^{n}}\hspace{-1mm}\vectornorm{\widehat{\boldsymbol{\Phi}}{}^{\ell,q}(\boldsymbol{\eta},\boldsymbol{\zeta})-\widehat{\boldsymbol{\Phi}}{}^{\ell,q-1}(\boldsymbol{\eta},\boldsymbol{\zeta})}^{2}_{\boldsymbol{W}}\hspace{-1mm}dP_{\boldsymbol{\eta}}^{\ell}dP_{\boldsymbol{\zeta}}}\leq\epsilon_{2}\sqrt{\int_{\real^{d}}\int_{\real^{n}}\hspace{-1mm}\vectornorm{\widehat{\boldsymbol{\Phi}}{}^{\ell,q}(\boldsymbol{\eta},\boldsymbol{\zeta})}^{2}_{\boldsymbol{W}}\hspace{-1mm}dP_{\boldsymbol{\eta}}^{\ell}dP_{\boldsymbol{\zeta}}},\quad\forall\ell\in\integer,\label{eq:criterion2}
\end{equation}}\normalsize
where~$\epsilon_{2}$ is a prescribed tolerance level.
Clearly, these criteria may result in the dependence of~$d$ and~$q$ on the number of iterations~$\ell$.

\subsection{Concluding remarks}
Algorithm~\ref{algo:algo6} summarizes an implementation of the problem in which the nonintrusive projection method is used for the solution of the subproblems.
Although this algorithm uses the nonintrusive projection method, we note that the proposed methodology can be readily adapted for use with other methods, such as embedded projection and collocation.

The main feature of the proposed implementation is that it provides a solution of the neutronics subproblem in a reduced-dimensional space when the reduced chaos expansion with random coefficients can extract a low-dimensional representation of the random temperature ($d<m$), while maintaining accuracy.
The solution in a reduced-dimensional space can be expected to reduce the number of terms required in the chaos expansion of the random neutron flux to achieve sufficient accuracy.
Further, it can be expected to reduce the number of quadrature nodes required for the nonintrusive projection method to achieve sufficient accuracy in the coefficients of the chaos expansion of the random neutron flux.
Hence, the solution of the neutronics subproblem in a reduced-dimensional space reduces the number of times a sample of the neutronics subproblem must be solved, thus in turn lowering the computational cost.

\begin{algorithm}[htp]
\SetKwInOut{Input}{Input}\SetKwInOut{Output}{Output}
\SetKw{KwAnd}{and}
\SetKwBlock{FirstMonoProblem}{neutronics subproblem}{end}
\SetKwBlock{SecondMonoProblem}{heat subproblem}{end}
\SetKwBlock{NexusRegion}{dimension reduction}{end}
\SetKwBlock{NexusRegionn}{measure transformation}{end}
\SetKwBlock{Initialization}{initialization}{end}
\Input{Error tolerance levels $\small{\epsilon_{1}}\normalsize$ and $\small{\epsilon_{2}}\normalsize$;\\
Basis $\small{\big\{\varphi_{\boldsymbol{\alpha}}\psi_{\boldsymbol{\beta}},\;0\leq|\boldsymbol{\alpha}|+|\boldsymbol{\beta}|\leq p\big\}}\normalsize$ up to total degree $\small{p}\normalsize$;\\
Quadrature rule $\small{\{(\boldsymbol{\xi}_{k},\boldsymbol{\zeta}_{k},v_{k}),\;1\leq k\leq N\}}\normalsize$ of level $\small{p+1}\normalsize$\;}
\vspace{-1mm}
$\small{\ell=1}\normalsize$\;
\Repeat{$($convergence$)$}{
\SecondMonoProblem{
\For{$k=1$ \KwTo $\,N$}{
Solve $\small{\big[\boldsymbol{K}+\boldsymbol{H}(\boldsymbol{\xi}_{k})\big]\widehat{\boldsymbol{T}}{}^{\ell}\big(\boldsymbol{\xi}_{k},\boldsymbol{\zeta}_{k}\big)=\boldsymbol{q}\Big(\widehat{\boldsymbol{T}}{}^{\ell-1,p}(\boldsymbol{\xi}_{k},\boldsymbol{\zeta}_{k}),\widehat{\boldsymbol{\Phi}}{}^{\ell-1,q}\big(\boldsymbol{\eta}^{\ell-1,p}(\boldsymbol{\xi}_{k}),\boldsymbol{\zeta}_{k}\big)\Big)}\normalsize$\;}
Compute chaos coordinates of $\small{\widehat{\boldsymbol{T}}{}^{\ell,p}}\normalsize$ using $\small{\widehat{\boldsymbol{T}}{}^{\ell}_{\boldsymbol{\alpha}\boldsymbol{\beta}}=\sum_{k=1}^{N}\widehat{\boldsymbol{T}}{}^{\ell}\big(\boldsymbol{\xi}_{k},\boldsymbol{\zeta}_{k}\big)\varphi_{\boldsymbol{\alpha}}(\boldsymbol{\xi}_{k})\psi_{\boldsymbol{\beta}}(\boldsymbol{\zeta}_{k})v_{k}}\normalsize$\;
}
\NexusRegion{  
Compute $\small{\overline{\boldsymbol{T}}{}^{\ell}_{\boldsymbol{\beta}}=\widehat{\boldsymbol{T}}{}^{\ell}_{\boldsymbol{0}\boldsymbol{\beta}}}\normalsize$ and $\small{\boldsymbol{C}{}^{\ell}_{\widehat{\boldsymbol{T}}_{\boldsymbol{\beta}}\widehat{\boldsymbol{T}}_{\tilde{\boldsymbol{\beta}}}}=\sum_{|\boldsymbol{\alpha}|=1}^{p-\max(|\boldsymbol{\beta}|,|\tilde{\boldsymbol{\beta}}|)}\widehat{\boldsymbol{T}}{}^{\ell}_{\boldsymbol{\alpha}\boldsymbol{\beta}}(\widehat{\boldsymbol{T}}{}^{\ell}_{\boldsymbol{\alpha}\tilde{\boldsymbol{\beta}}})^{\mathrm{T}}}\normalsize$\;
Solve the eigenproblem $\small{\sum_{|\tilde{\boldsymbol{\beta}}|=0}^{p-1}\boldsymbol{W}^{\mathrm{T}}\boldsymbol{C}{}^{\ell}_{\widehat{\boldsymbol{T}}_{\boldsymbol{\beta}}\widehat{\boldsymbol{T}}_{\tilde{\boldsymbol{\beta}}}}\boldsymbol{W}\boldsymbol{\phi}^{j,\ell}_{\tilde{\boldsymbol{\beta}}}=\lambda_{j}^{\ell}\boldsymbol{W}\boldsymbol{\phi}^{j,\ell}_{\boldsymbol{\beta}},\;0\leq|\boldsymbol{\beta}|\leq p-1}\normalsize$\;
Choose $\small{d}\normalsize$ such that $\small{\sqrt{\sum_{j=d+1}^{r}\lambda_{j}^{\ell}}\leq\epsilon_{1}\sqrt{\sum_{|\boldsymbol{\alpha}|+|\boldsymbol{\beta}|=0}^{p}(\widehat{\boldsymbol{T}}{}_{\boldsymbol{\alpha}\boldsymbol{\beta}}^{\ell})^{\mathrm{T}}\boldsymbol{W}\widehat{\boldsymbol{T}}{}_{\boldsymbol{\alpha}\boldsymbol{\beta}}^{\ell}}}\normalsize$\;
Compute coordinates of $\small{\eta_{j}^{\ell,p}}\normalsize$ using $\small{\eta_{j,\boldsymbol{\alpha}}^{\ell}=\sum_{|\boldsymbol{\beta}|=0}^{p-1}(\widehat{\boldsymbol{T}}{}^{\ell}_{\boldsymbol{\alpha}\boldsymbol{\beta}})^{\mathrm{T}}\boldsymbol{W}\boldsymbol{\phi}_{\boldsymbol{\beta}}^{j,\ell}/\sqrt{\lambda_{j}^{\ell}}}\normalsize$\;
}
\FirstMonoProblem{
$\small{q=0}\normalsize$\;
\Repeat{$\small{\Big(\small{\sqrt{\sum_{|\boldsymbol{\gamma}|+|\boldsymbol{\beta}|=q}\|\widehat{\boldsymbol{\Phi}}{}^{\ell}_{\boldsymbol{\gamma}\boldsymbol{\beta}}\|_{\boldsymbol{W}}^{2}}\leq\epsilon_{2}\sqrt{\sum_{|\boldsymbol{\gamma}|+|\boldsymbol{\beta}|=0}^{q}\|\widehat{\boldsymbol{\Phi}}{}^{\ell}_{\boldsymbol{\gamma}\boldsymbol{\beta}}\|_{\boldsymbol{W}}^{2}}}\normalsize\;\Big)}\normalsize$}{
\NexusRegionn{
Compute basis $\small{\{\Gamma_{\boldsymbol{\gamma}}^{\ell}\psi_{\boldsymbol{\beta}},0\leq|\boldsymbol{\gamma}|+|\boldsymbol{\beta}|\leq q\}}\normalsize$ up to total degree $\small{q}\normalsize$\;
Compute quadrature rule $\small{\{(\boldsymbol{\eta}_{k}^{\ell},\boldsymbol{\zeta}_{k}^{\ell},w_{k}^{\ell}),1\leq k\leq \nu^{\ell}\}}\normalsize$ of level $\small{q+2}\normalsize$\;
}  
\For{$k=1$ \KwTo $\,\nu^{\ell}$}{
Solve $\small{\big[\boldsymbol{D}\big(\widehat{\boldsymbol{T}}{}^{\ell,p,d}(\boldsymbol{\eta}^{\ell}_{k},\boldsymbol{\zeta}^{\ell}_{k})\big)+\boldsymbol{M}\big(\widehat{\boldsymbol{T}}{}^{\ell,p,d}(\boldsymbol{\eta}^{\ell}_{k},\boldsymbol{\zeta}^{\ell}_{k}),\boldsymbol{\zeta}^{\ell}_{k}\big)\big]\widehat{\boldsymbol{\Phi}}{}^{\ell}(\boldsymbol{\eta}^{\ell}_{k},\boldsymbol{\zeta}^{\ell}_{k})=\boldsymbol{s},}\normalsize$\\
with $\small{\widehat{\boldsymbol{T}}{}^{\ell,p,d}(\boldsymbol{\eta}^{\ell}_{k},\boldsymbol{\zeta}^{\ell}_{k})=\overline{\boldsymbol{T}}{}^{\ell,p}(\boldsymbol{\zeta}^{\ell}_{k})+\sum_{j=1}^{d}\sqrt{\lambda_{j}^{\ell}}\eta_{j,k}^{\ell}\boldsymbol{\phi}^{j,\ell,p-1}(\boldsymbol{\zeta}^{\ell}_{k})}\normalsize$\;}
Compute coordinates of $\small{\widehat{\boldsymbol{\Phi}}{}^{\ell,q}}\normalsize$ using $\small{\widehat{\boldsymbol{\Phi}}{}^{\ell}_{\boldsymbol{\gamma}\boldsymbol{\beta}}=\sum_{k=1}^{\nu^{\ell}}\widehat{\boldsymbol{\Phi}}{}^{\ell}(\boldsymbol{\eta}^{\ell}_{k},\boldsymbol{\zeta}^{\ell}_{k})\Gamma_{\boldsymbol{\gamma}}^{\ell}(\boldsymbol{\eta}^{\ell}_{k})\psi_{\boldsymbol{\beta}}(\boldsymbol{\zeta}^{\ell}_{k})w^{\ell}_{k}}\normalsize$\;
$\small{q=q+1}\normalsize$\;
\vspace{1mm}
}
}
$\small{\ell=\ell+1}\normalsize$\;}
\caption{Implementation of the illustration problem.}\label{algo:algo6}
\end{algorithm}

\section{Numerical results}\label{sec:sec8}
We obtained numerical results by considering the following parameter values.  
We assumed the reactor to have a length~$L$ of $100\,[\text{cm}]$.
Further, we assumed a deterministic and position-independent heat conductivity~$k$ of~$100\,\text{$[\text{J/(K\,cm\,s)}]$}$, ambient temperature~$T_{\infty}$ of $390\,[\text{K}]$, a fission energy $E_{\text{f}}$ of $3.0E\text{-}11\,[\text{J/neutrons}]$, a fission cross section $\Sigma_{\text{a,ref}}$ of $0.0075\,[\text{cm}^{-1}]$, a neutron-diffusion constant $D_{\text{ref}}$ of $2.2\,[\text{cm}]$, a multiplication factor $\nu$ of $2.2$, a neutron source strength $s$ of $5.0E11\,[\text{neutrons/(s\,cm$^{3}$)}$].
In addition, the following temperatures were assumed: $T_{\text{ref}}=390\,[\text{K}]$, $T_{\text{min}}=390\,[\text{K}]$, and $T_{\text{max}}=1000\,[\text{K}]$.

\begin{figure}[htp]
  \begin{center}
    \subfigure[Samples.]{\includegraphics[width=0.8\textwidth]{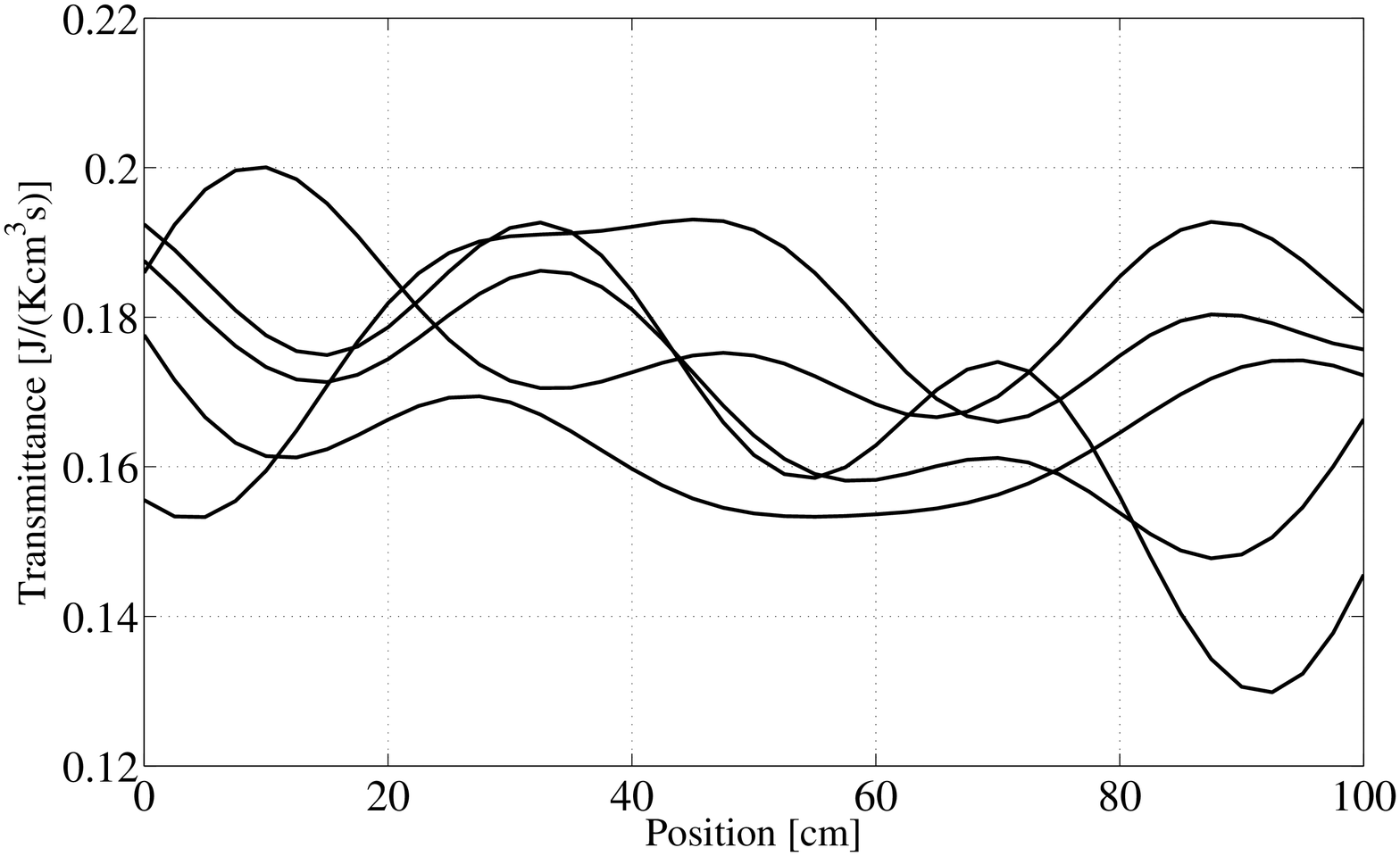}}
    \hfill
    \subfigure[Eigenvalues.]{\includegraphics[width=0.8\textwidth]{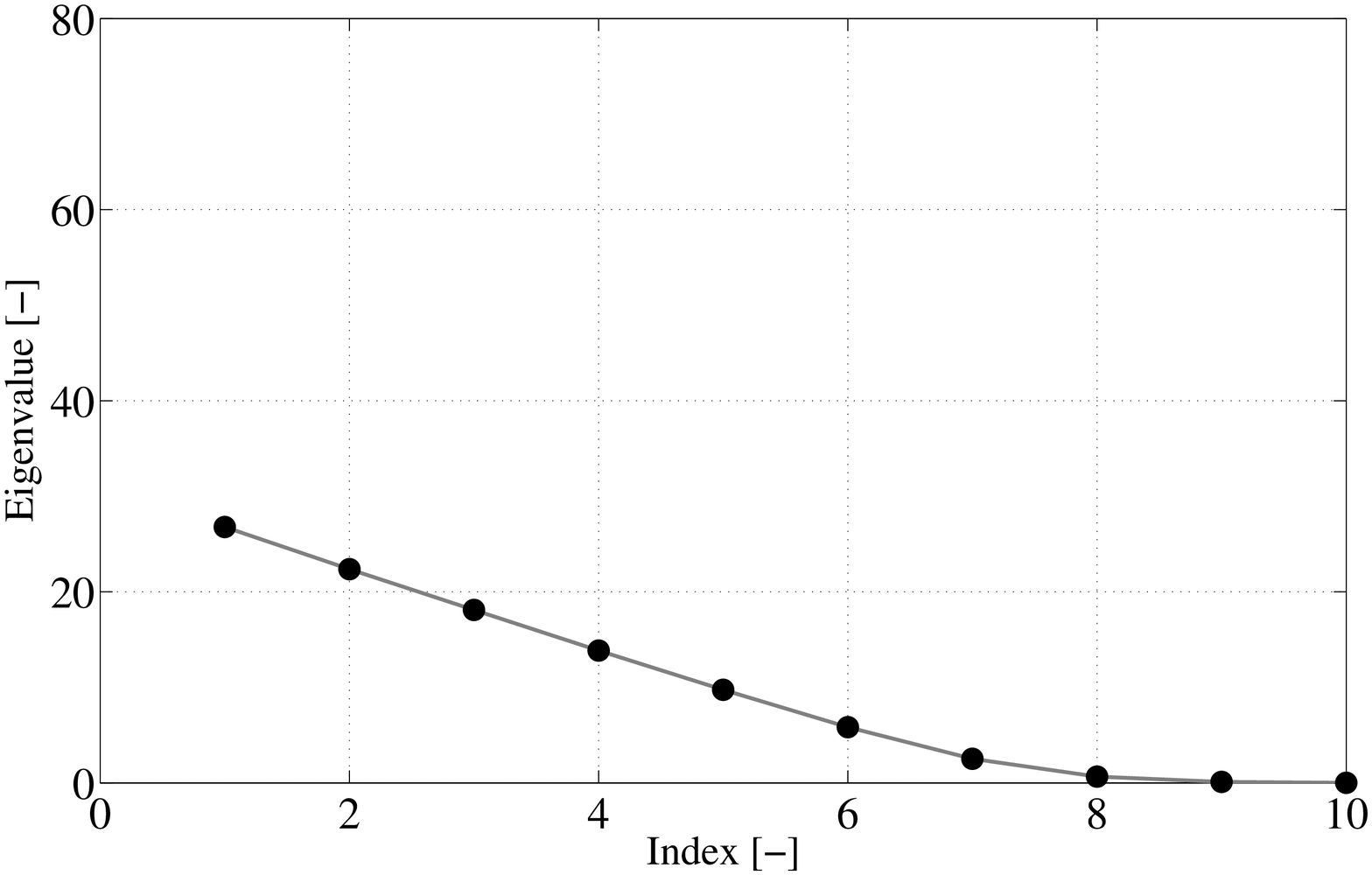}}
    \caption{Thermal transmittance random field: (a)~five samples and~(b)~ten largest eigenvalues of the covariance integral operator whose kernel is given by~(\ref{eq:covara}).}\label{fig:figure1ab}
  \end{center}
\end{figure}

\begin{figure}[htp]
  \begin{center}
    \subfigure[Samples.]{\includegraphics[width=0.8\textwidth]{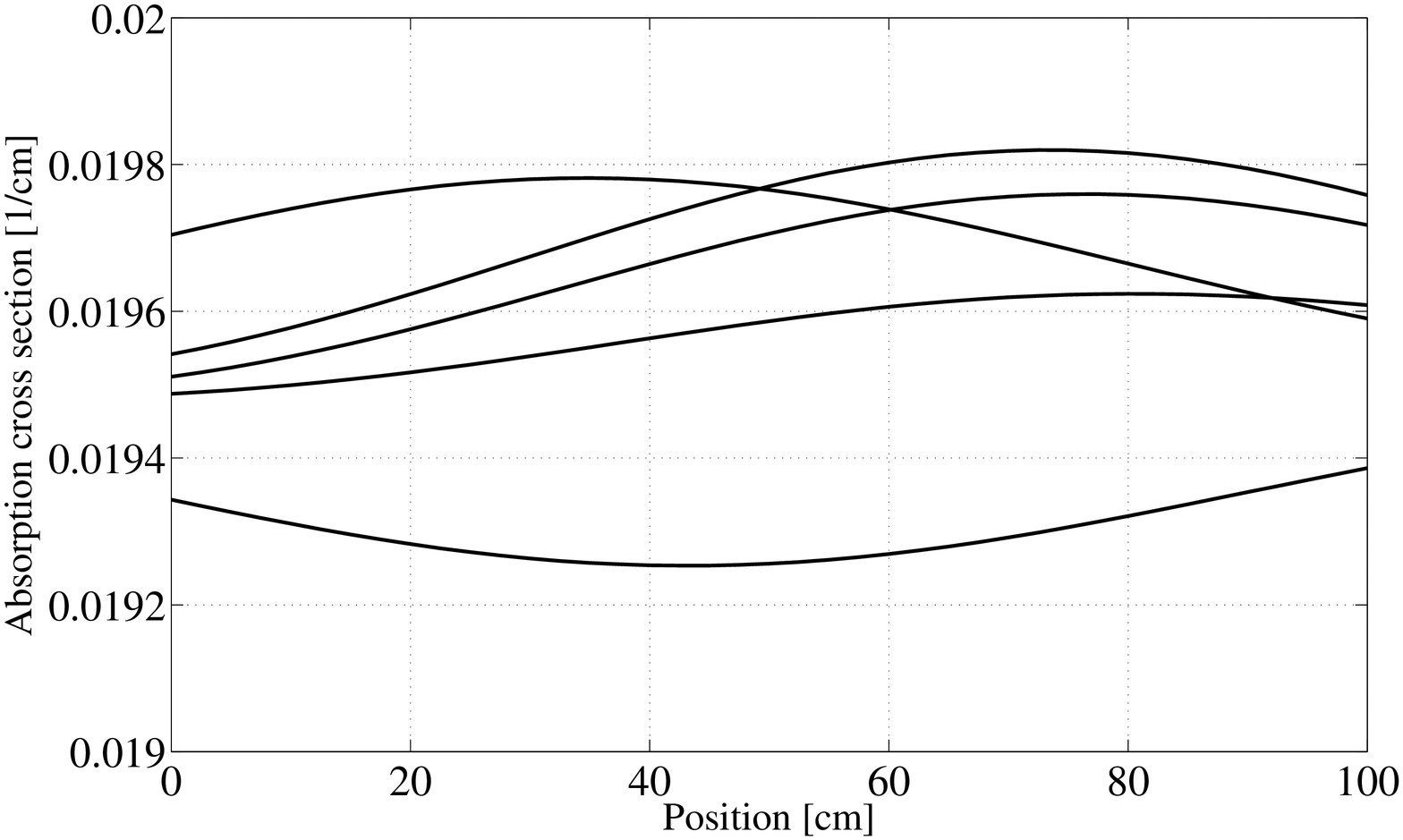}}
    \hfill
    \subfigure[Eigenvalues.]{\includegraphics[width=0.8\textwidth]{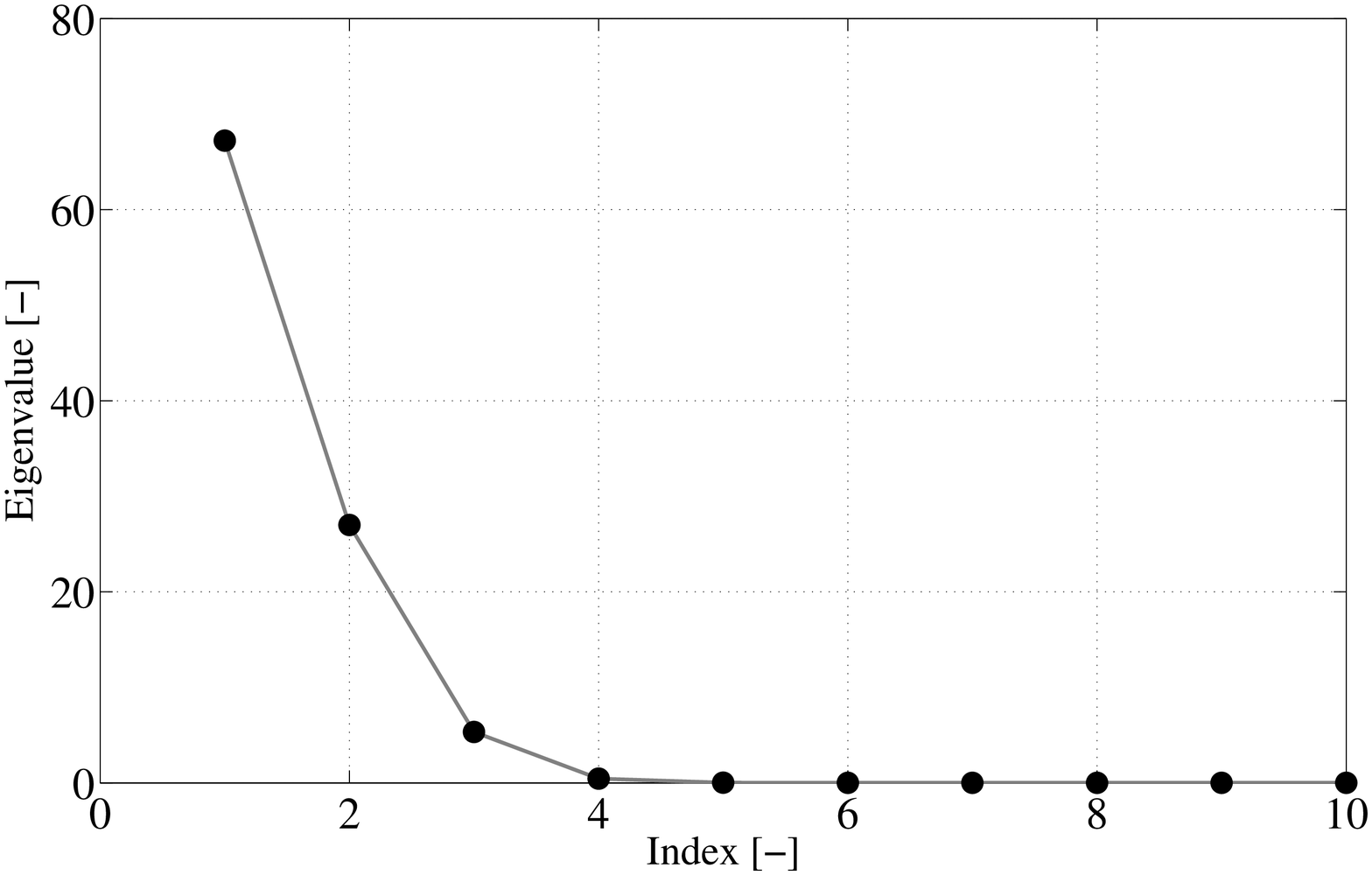}}
    \caption{Absorption cross section random field: (a)~five samples and~(b)~ten largest eigenvalues of the covariance integral operator whose kernel is given by~(\ref{eq:covarb}).}\label{fig:figure1cd}
  \end{center}
\end{figure}

We used thermal transmittance and absorption cross section random fields with position-independent mean values~$\overline{h}=0.17\,[\text{J/(K\,cm$^{3}$\,s)}]$ and~$\Sigma_{\text{a,ref}}=0.0195\,[\text{cm}^{-1}]$, spatial correlation lengths~$a_{h}=15\,[\text{cm}]$ and~$a_{\Sigma}=50\,[\text{cm}]$, and coefficients of variation~$\delta_{h}=10\,\%$ and~$\delta_{\Sigma}=10\,\%$.
We retained~$m=10$ and~$n=2$ terms in expansions~(\ref{eq:hN}) and~(\ref{eq:Sigmaan}), respectively.

Figures~\ref{fig:figure1ab}(a) and~\ref{fig:figure1cd}(a) show a few sample paths of the random fields thus obtained.
We can observe that the samples of the thermal transmittance random field with~$a_{h}=15\,[\text{cm}]$ are less smooth than those of the absorption cross section random field with~$a_{\Sigma}=50\,[\text{cm}]$, that is, the samples of the former random field exhibit more rapid oscillations with respect to the position in the reactor than those of the latter random field.
Figures~\ref{fig:figure1ab}(b) and~\ref{fig:figure1cd}(b) show the 10 largest eigenvalues of the covariance integral operators whose kernels are given by~(\ref{eq:covara}) and~(\ref{eq:covarb}).  
We can observe that the eigenvalues obtained for the thermal transmittance random field decay at a lower rate than those obtained for the absorption cross section random field, indicating the adequacy of the truncation of expansions~(\ref{eq:hN}) and~(\ref{eq:Sigmaan}) after~$m=10$ and~$n=2$ terms.

\subsection{Monte Carlo sampling implementation}
\begin{figure}[htp]
  \begin{center}
    \subfigure[Temperature.]{\includegraphics[width=0.8\textwidth]{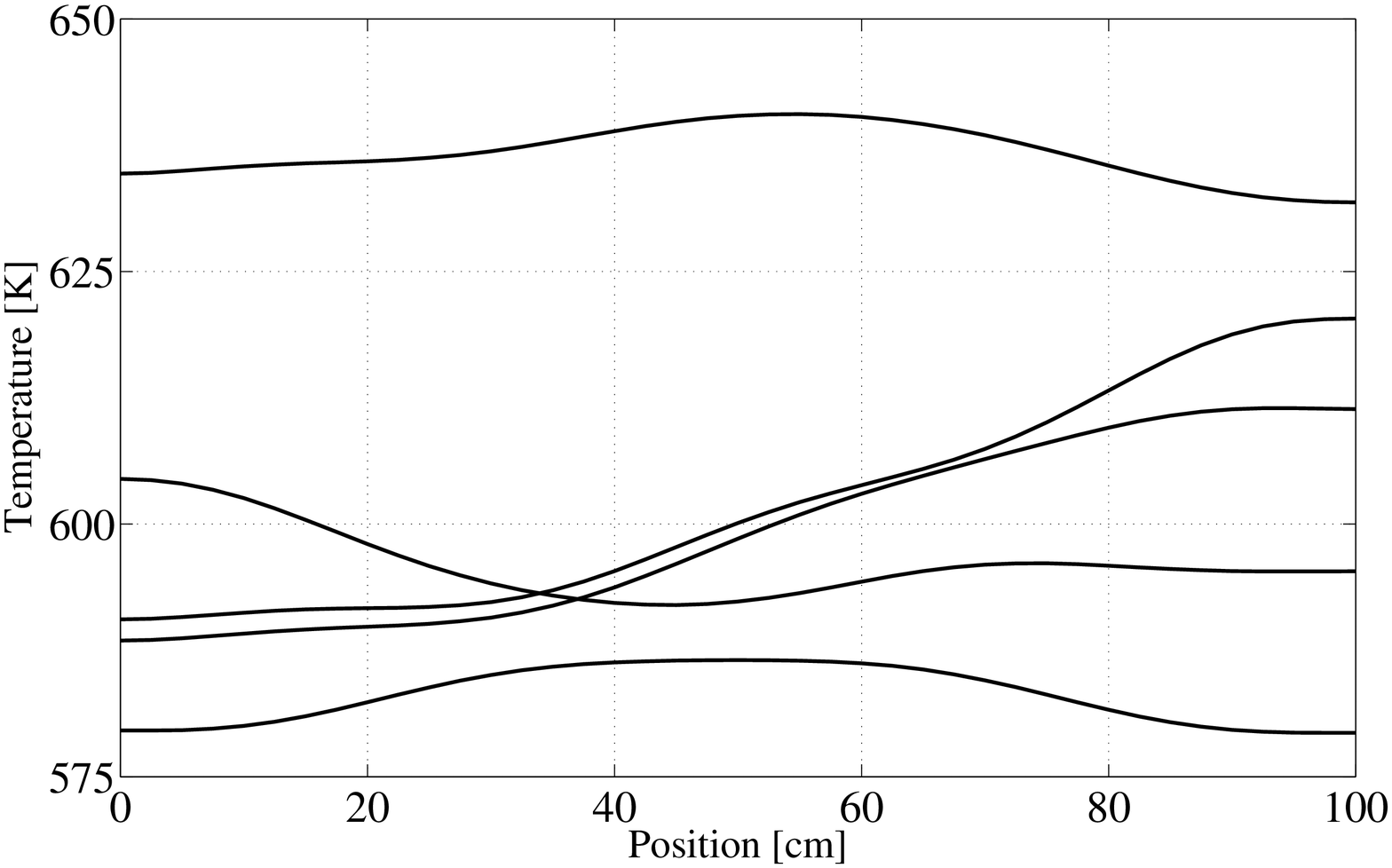}}
    \hfill
    \subfigure[Neutron flux.]{\includegraphics[width=0.8\textwidth]{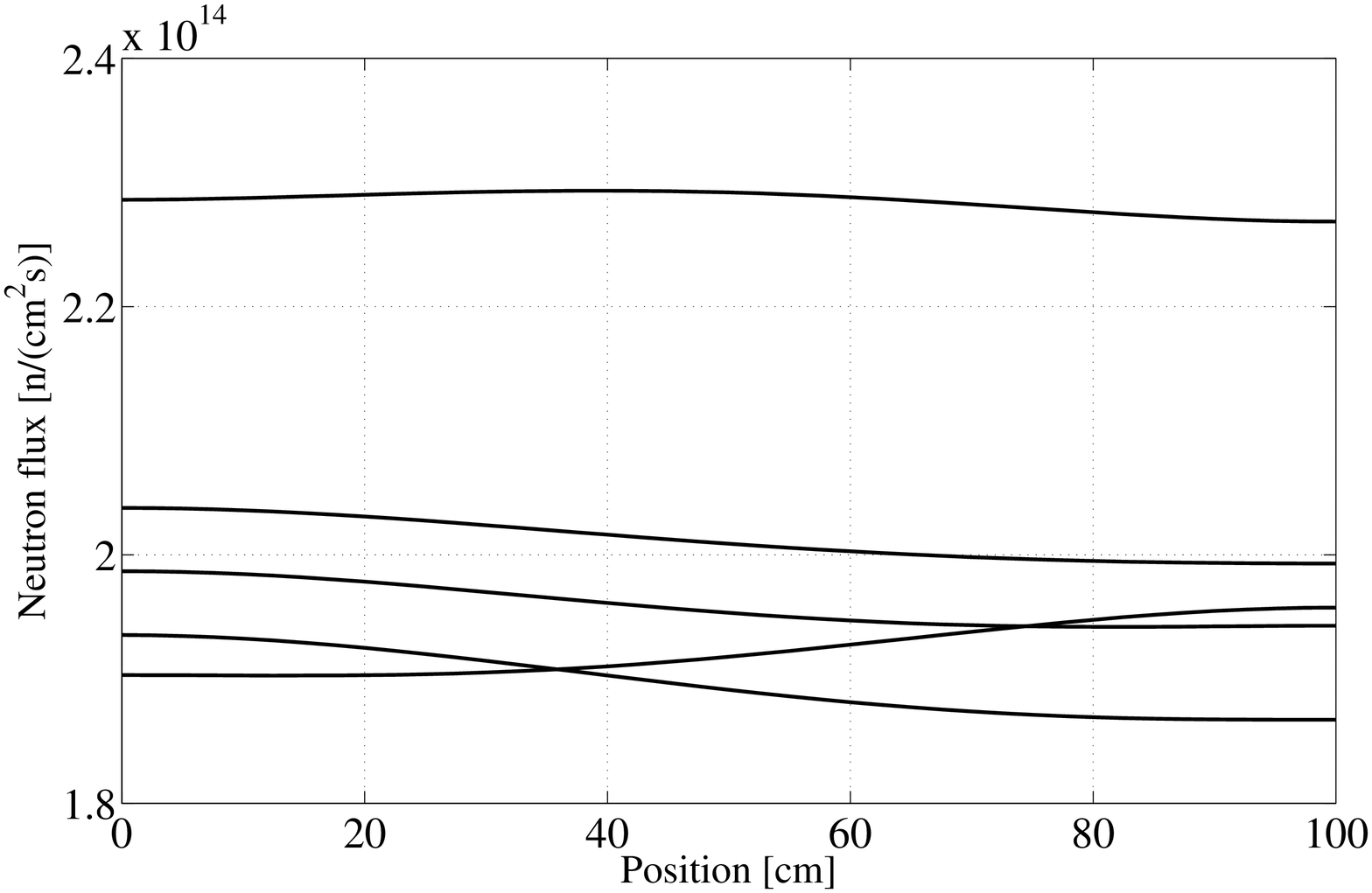}}
    \caption{Monte Carlo simulation: five samples of the solution.}\label{fig:figure2}
  \end{center}
\end{figure}
First, we carried out a Monte Carlo simulation.
We generated $MC=100,000$ pairs of sample paths of the thermal transmittance and absorption cross section random fields.
Then, for each pair of sample paths, we constructed the associated deterministic multiphysics model.
We solved each model by using the FE method for spatial discretization and Gauss-Seidel iteration.
We systematically obtained converged results for~$40$ finite elements and~$20$ iterations.  

Figure~\ref{fig:figure2} shows a few samples of the random temperature and neutron flux thus obtained.
The samples of the random temperature (Fig.~\ref{fig:figure2}(a)) are smoother than those of the thermal transmittance random field~(Fig.~\ref{fig:figure1ab}(a)), that is, the former samples exhibit less rapid oscillations with respect to the position in the reactor than the latter samples.
In~\citep{arnst2011a}, we had shown that this behavior can be attributed to the large magnitude of the diffusion term of the heat subproblem, which reduces the nonuniformity of the samples of the random temperature.

\subsection{PC-based implementation involving dimension reduction and measure transformation}
Next, we implemented the proposed polynomial-chaos-based iterative method involving dimension reduction and measure transformation.
This implementation corresponded exactly to Algorithm~\ref{algo:algo6}.
We obtained results by setting the total degree $p$ of the chaos expansion of the random temperature to $4$ and, with reference to~(\ref{eq:criterion1}) and (\ref{eq:criterion2}), by using a range of values for the error tolerance levels $\epsilon_{1}$ and~$\epsilon_{2}$ to determine the reduced dimension and the total degree of the chaos expansion of the random neutron flux at each iteration.
We discuss the convergence of the results as a function of these error tolerance levels later.
Now, we present detailed results obtained for~$\epsilon_{1}=0.01$ and~$\epsilon_{2}=0.01$.

\begin{figure}[htp]
  \begin{center}
    \includegraphics[width=0.8\textwidth]{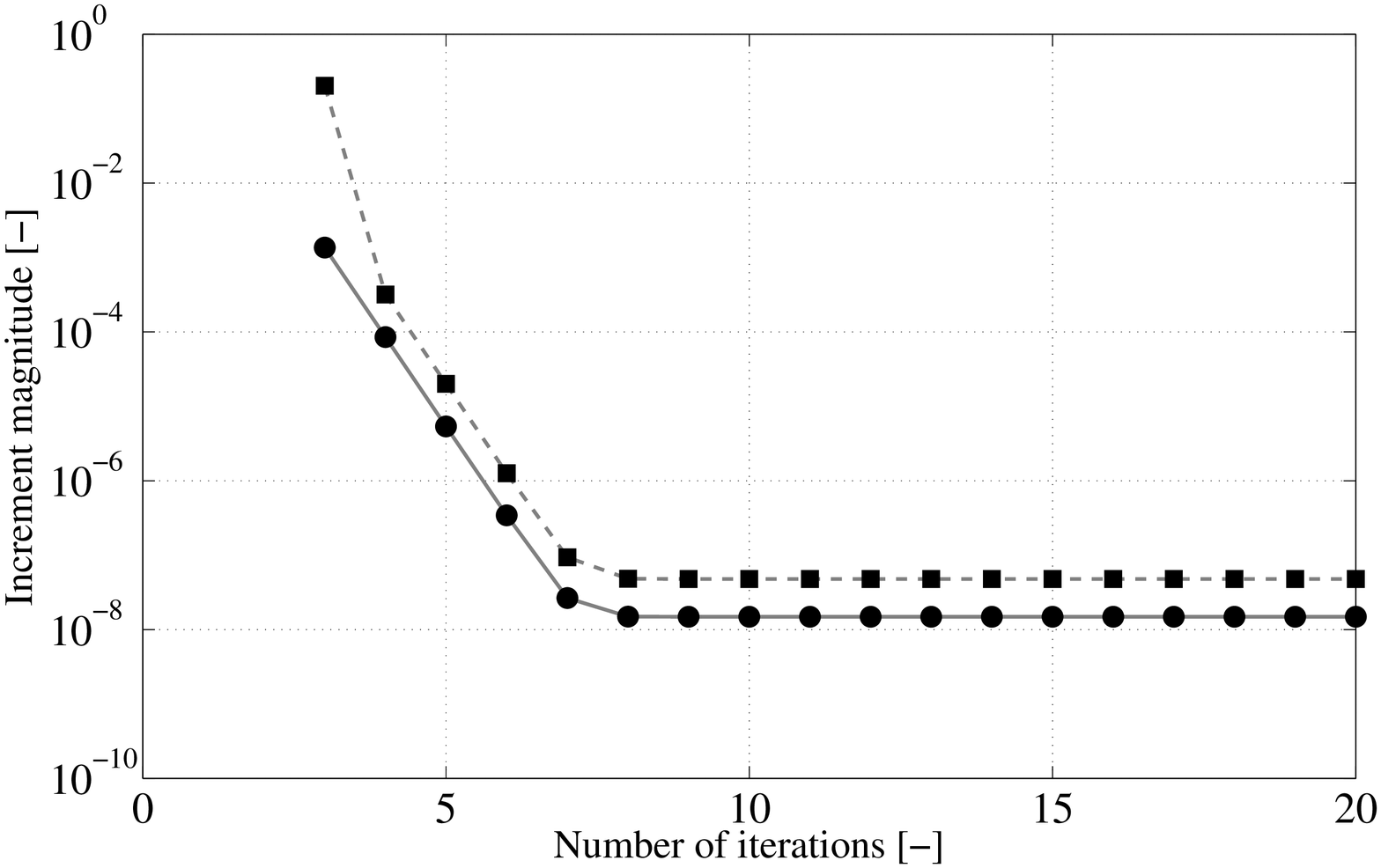}
    \centerline{$\small{\ell\mapsto\sqrt{\frac{1}{MC}\sum_{k=1}^{MC}\|\widehat{\boldsymbol{T}}{}^{\ell,p}(\boldsymbol{\xi}_{k},\boldsymbol{\zeta}_{k})-\widehat{\boldsymbol{T}}{}^{\ell-1,p}(\boldsymbol{\xi}_{k},\boldsymbol{\zeta}_{k})\|_{\boldsymbol{W}}^{2}}\Big/\sqrt{\frac{1}{MC}\sum_{k=1}^{MC}\|{\boldsymbol{T}}{}^{\infty}(\boldsymbol{\xi}_{k},\boldsymbol{\zeta}_{k})\|_{\boldsymbol{W}}^{2}}}\normalsize$ (circles).} 
    \centerline{$\small{\ell\mapsto\sqrt{\frac{1}{MC}\sum_{k=1}^{MC}\|\widehat{\boldsymbol{\Phi}}{}^{\ell,q}(\boldsymbol{\eta}^{\ell,p}(\boldsymbol{\xi}_{k}),\boldsymbol{\zeta}_{k})-\widehat{\boldsymbol{\Phi}}{}^{\ell-1,q}(\boldsymbol{\eta}^{\ell-1,p}(\boldsymbol{\xi}_{k}),\boldsymbol{\zeta}_{k})\|_{\boldsymbol{W}}^{2}}\Big/\sqrt{\frac{1}{MC}\sum_{k=1}^{MC}\|{\boldsymbol{\Phi}}{}^{\infty}(\boldsymbol{\xi}_{k},\boldsymbol{\zeta}_{k})\|_{\boldsymbol{W}}^{2}}}\normalsize$ (squares).}
    \caption{PC-based simulation: convergence with respect to the number of iterations.}\label{fig:figure3}
  \end{center}
\end{figure}

Figure~\ref{fig:figure3} shows the convergence of the iterative method as a function of the number of iterations; note that the superscript $\infty$ is used in the figure captions to indicate convergence with respect to the number of iterations.
The iterative method converged linearly up to approximately iteration~$\ell=7$, after which linear-solver tolerances became dominant and prevented further convergence. 
All results to follow were obtained at iteration~$\ell=20$ and can thus be considered to have converged with respect to the number of iterations.

\begin{figure}[htp]
  \begin{center}
    \subfigure[Mean (thick solid) and first (thin solid), second (dashed), third (dash-dotted), and fourth (dotted) order coefficients $\overline{\boldsymbol{T}}{}^{\infty}_{\boldsymbol{\beta}}$ of the chaos expansion $\overline{\boldsymbol{T}}{}^{\infty,p}(\boldsymbol{\zeta})=\sum_{|\boldsymbol{\beta}|=0}^{p}\overline{\boldsymbol{T}}{}^{\infty}_{\boldsymbol{\beta}}\psi_{\boldsymbol{\beta}}(\boldsymbol{\zeta})$.]{\includegraphics[width=0.47\textwidth]{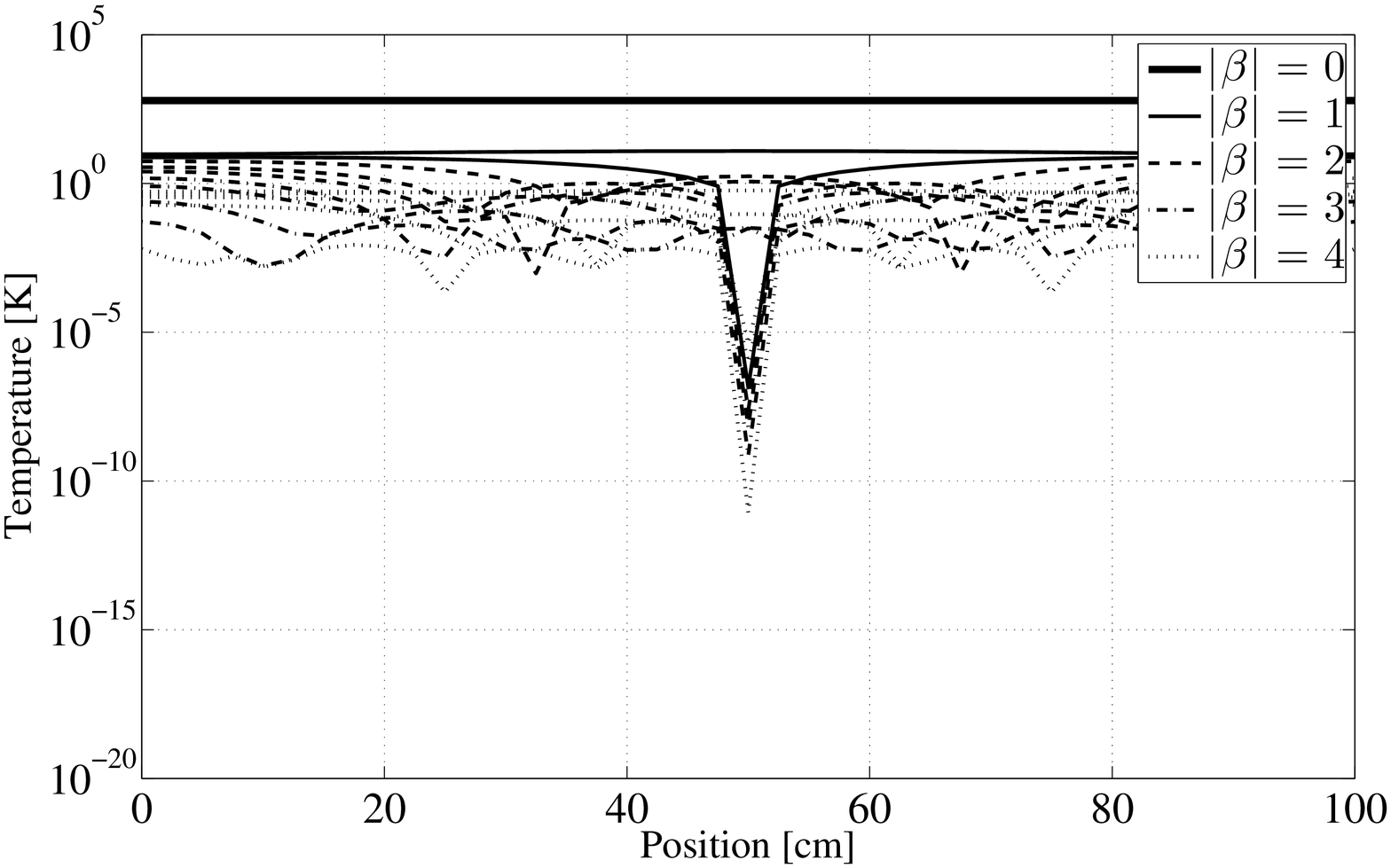}}
    \hfill
    \subfigure[Eigenvalues $\{\lambda_{j}^{\infty},1\leq j\leq 10\}$.]{\includegraphics[width= 0.47\textwidth]{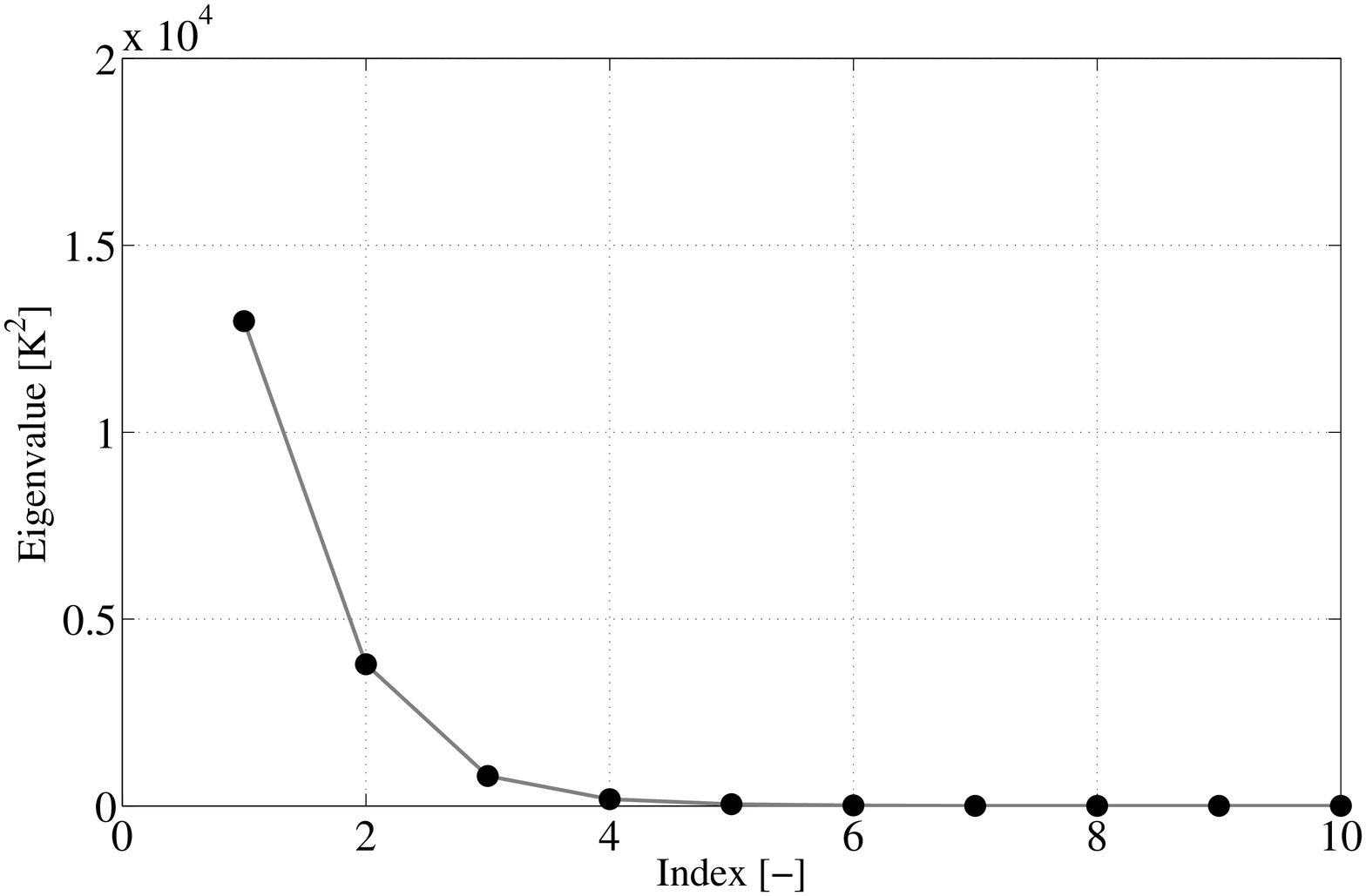}}
    \subfigure[Mean (thick solid) and first (thin solid), second (dashed), and third (dash-dotted) order coefficients $\boldsymbol{\phi}^{1,\infty}_{\boldsymbol{\beta}}$ of the chaos expansion $\boldsymbol{\phi}^{1,\infty,p-1}(\boldsymbol{\zeta})=\sum_{|\boldsymbol{\beta}|=0}^{p-1}\boldsymbol{\phi}^{1,\infty}_{\boldsymbol{\beta}}\psi_{\boldsymbol{\beta}}(\boldsymbol{\zeta})$.]{\includegraphics[width= 0.47\textwidth]{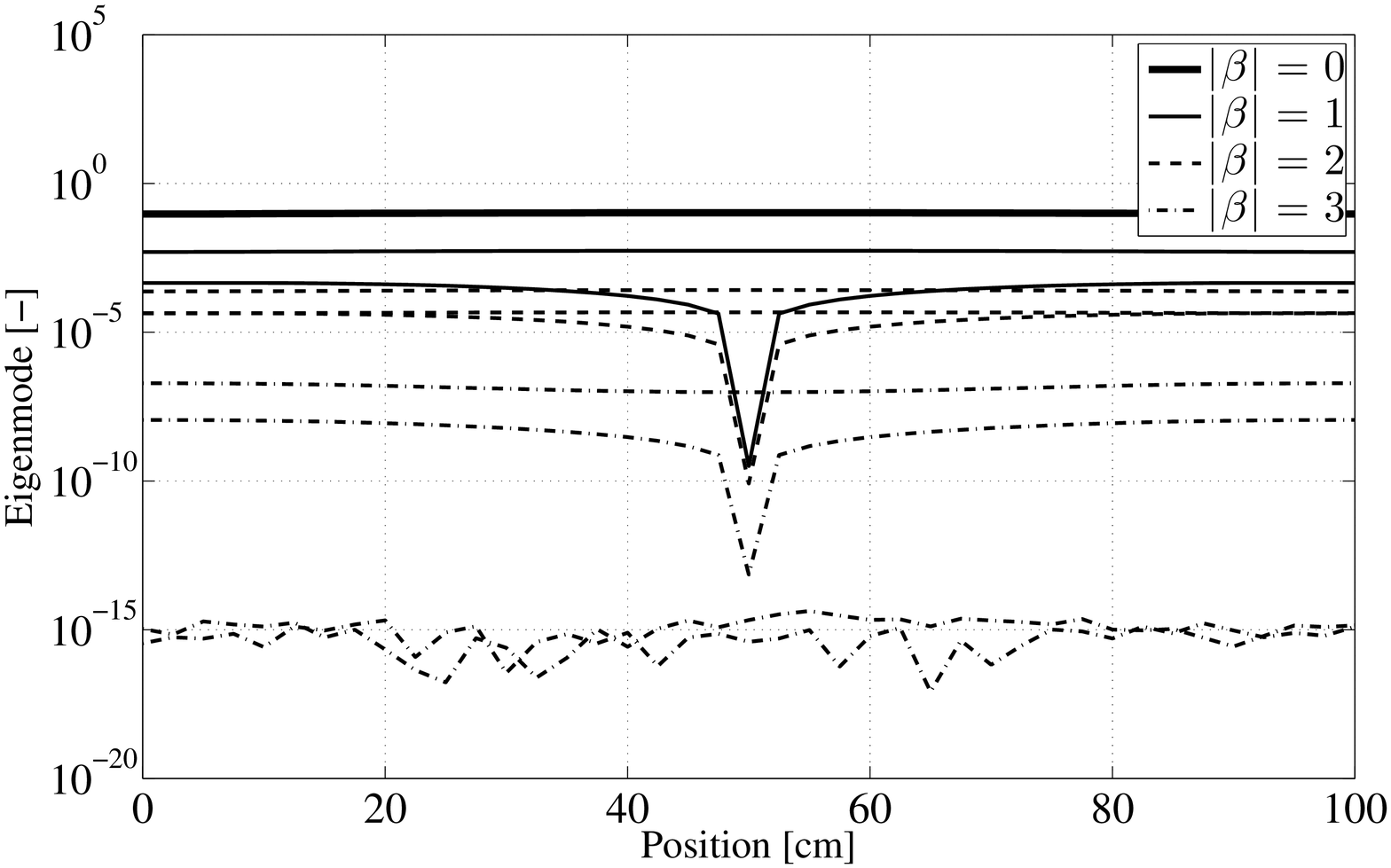}}
    \hfill
    \subfigure[Mean (thick solid) and first (thin solid), second (dashed), and third (dash-dotted) order coefficients $\boldsymbol{\phi}^{2,\infty}_{\boldsymbol{\beta}}$ of the chaos expansion $\boldsymbol{\phi}^{2,\infty,p-1}(\boldsymbol{\zeta})=\sum_{|\boldsymbol{\beta}|=0}^{p-1}\boldsymbol{\phi}^{2,\infty}_{\boldsymbol{\beta}}\psi_{\boldsymbol{\beta}}(\boldsymbol{\zeta})$.]{\includegraphics[width= 0.47\textwidth]{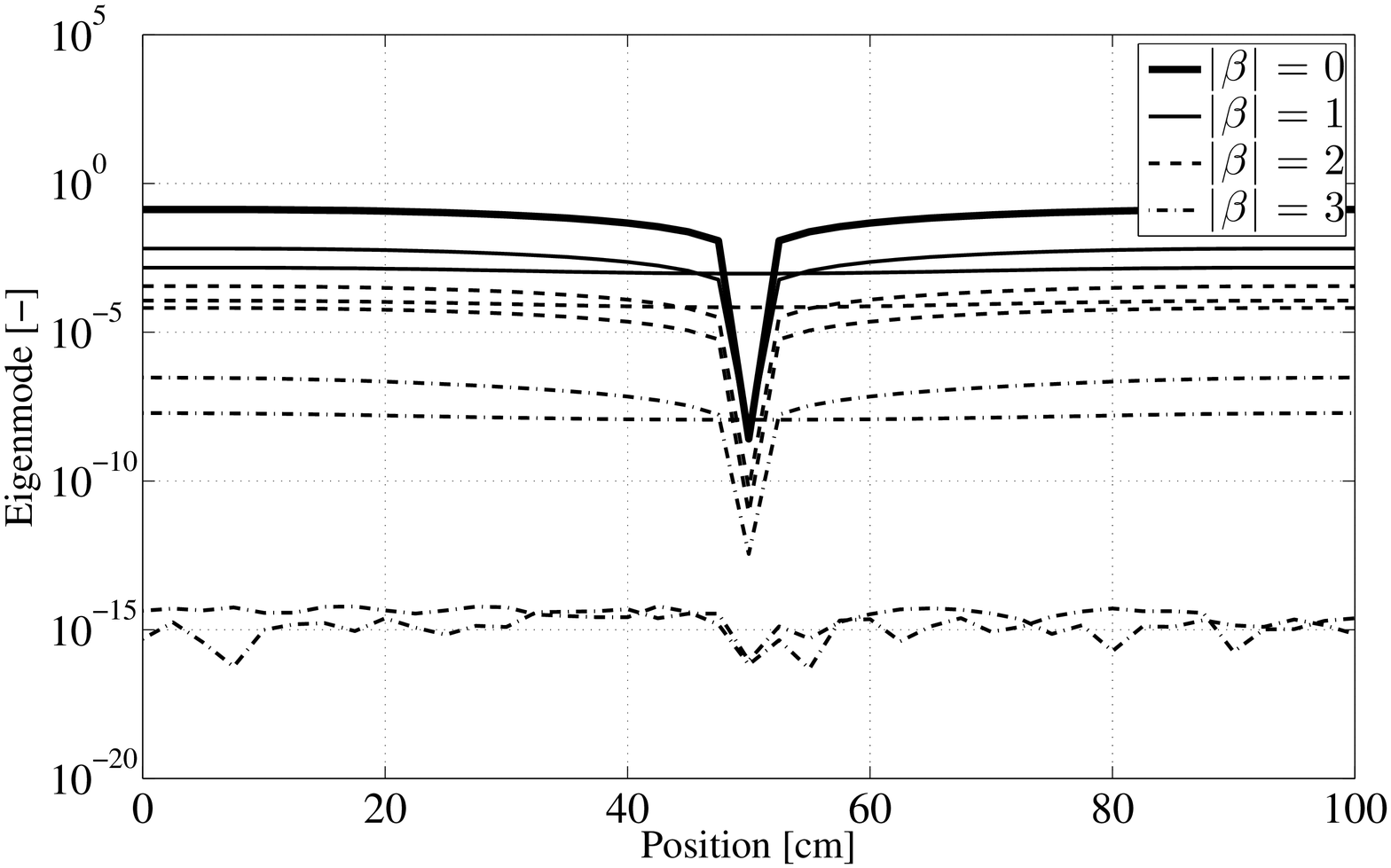}}
    \subfigure[Coefficients $\eta^{\infty}_{1,\boldsymbol{\alpha}}$ of the chaos expansion $\eta^{\infty,p}_{1}(\boldsymbol{\xi})=\sum_{|\boldsymbol{\alpha}|=1}^{p}\eta^{\infty}_{1,\boldsymbol{\alpha}}\varphi_{\boldsymbol{\alpha}}(\boldsymbol{\xi})$.]{\includegraphics[width= 0.47\textwidth]{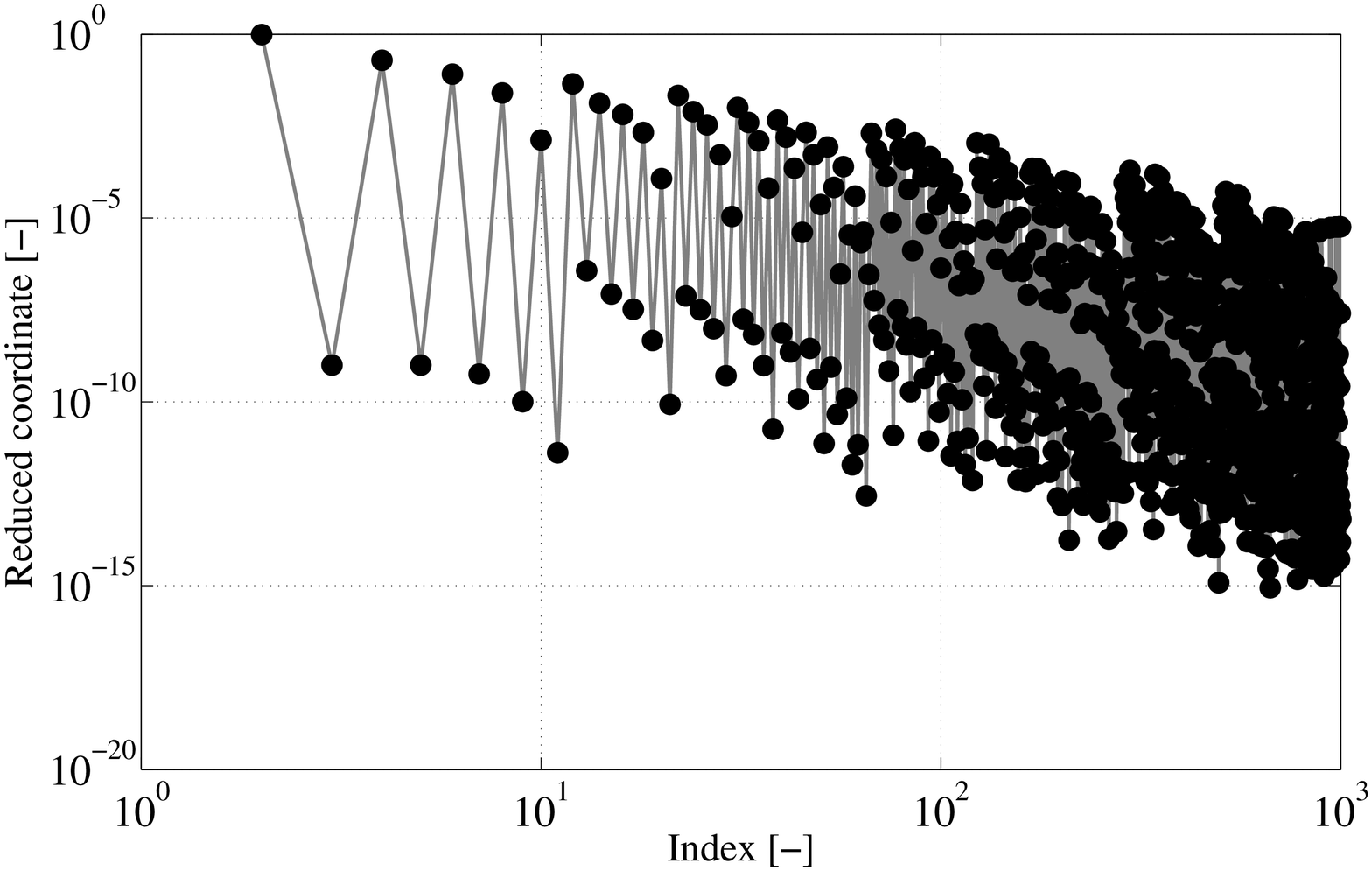}}
    \hfill
    \subfigure[Coefficients $\eta^{\infty}_{2,\boldsymbol{\alpha}}$ of the chaos expansion $\eta^{\infty,p}_{2}(\boldsymbol{\xi})=\sum_{|\boldsymbol{\alpha}|=1}^{p}\eta^{\infty}_{2,\boldsymbol{\alpha}}\varphi_{\boldsymbol{\alpha}}(\boldsymbol{\xi})$.]{\includegraphics[width= 0.47\textwidth]{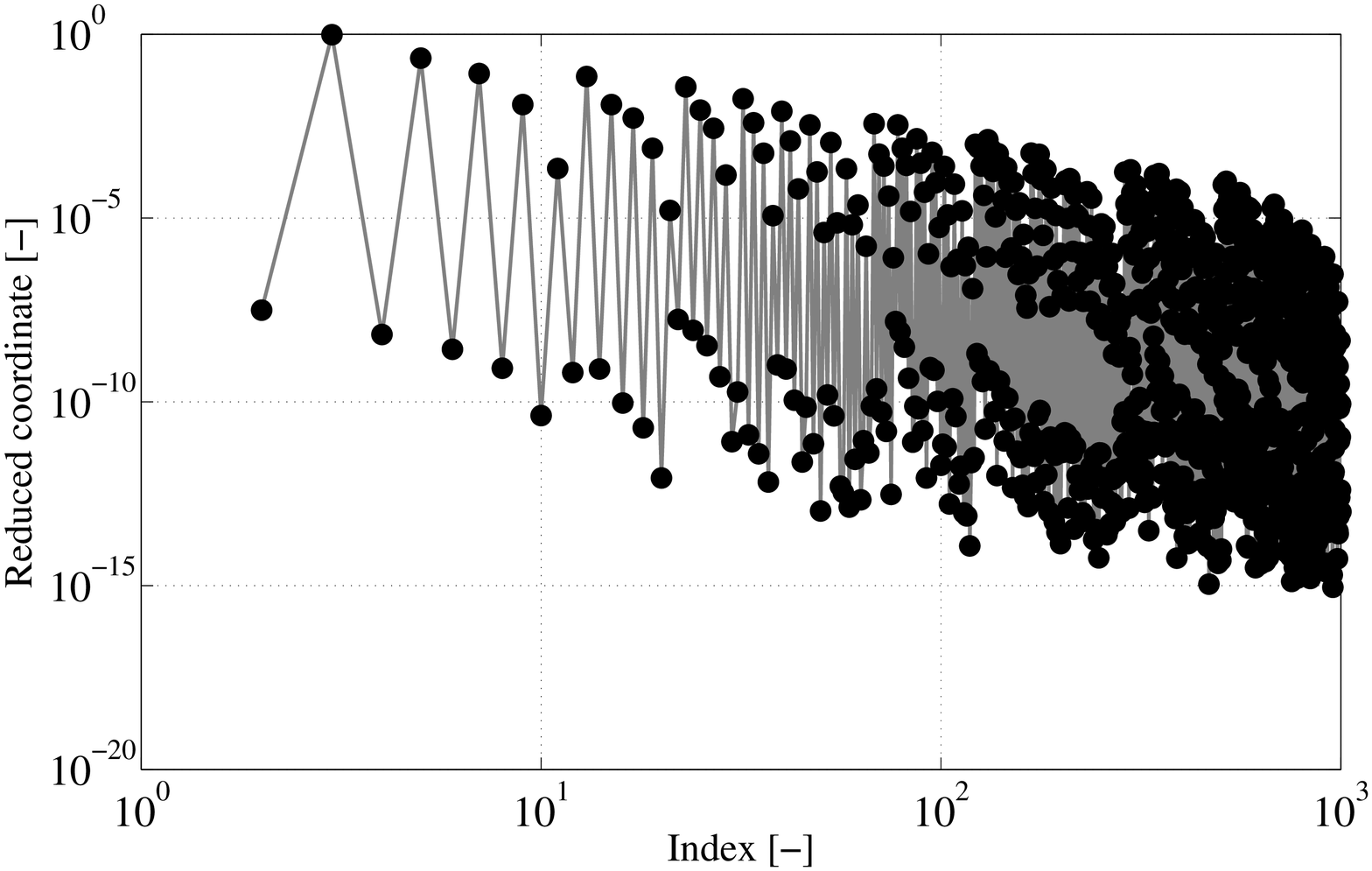}}\\
    \vspace{3mm}
    \caption{PC-based simulation: a few components of the reduced chaos expansion with random coefficients, $\widehat{\boldsymbol{T}}{}^{\infty,p,d}(\boldsymbol{\xi},\boldsymbol{\zeta})=\overline{\boldsymbol{T}}{}^{\infty,p}(\boldsymbol{\zeta})+\sum_{j=1}^{d}\sqrt{\lambda_{j}^{\infty}}\,\eta_{j}^{\infty,p}(\boldsymbol{\xi})\,\boldsymbol{\phi}^{j,\infty,p-1}(\boldsymbol{\zeta})$, of the random temperature.}\label{fig:figure4}
  \end{center}
\end{figure}

\begin{figure}[htp]
  \begin{center}
    \subfigure[Joint probability density function.]{\includegraphics[width=0.55\textwidth]{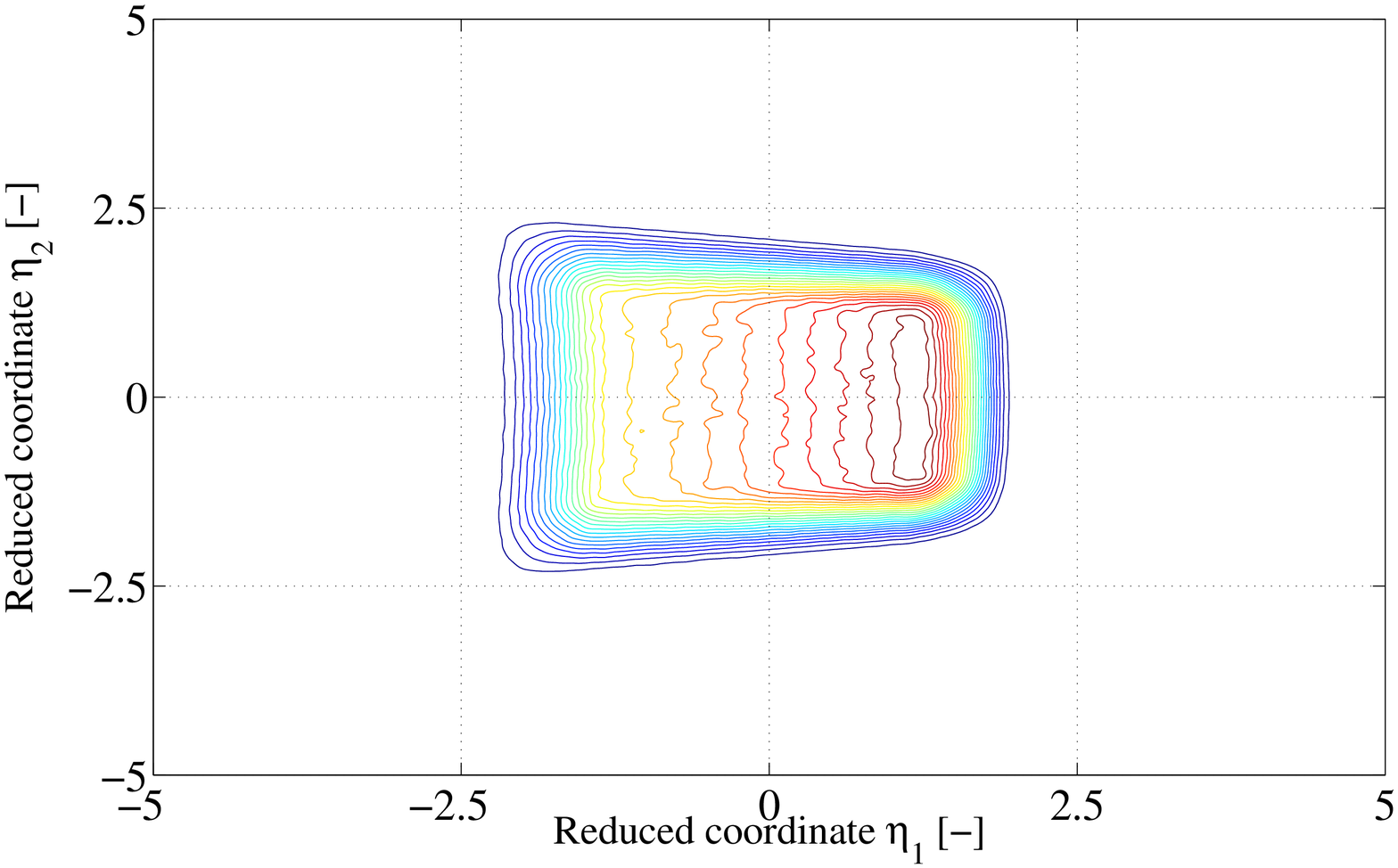}}
    \vfill
    \subfigure[First marginal probability density function.]{\includegraphics[width=0.55\textwidth]{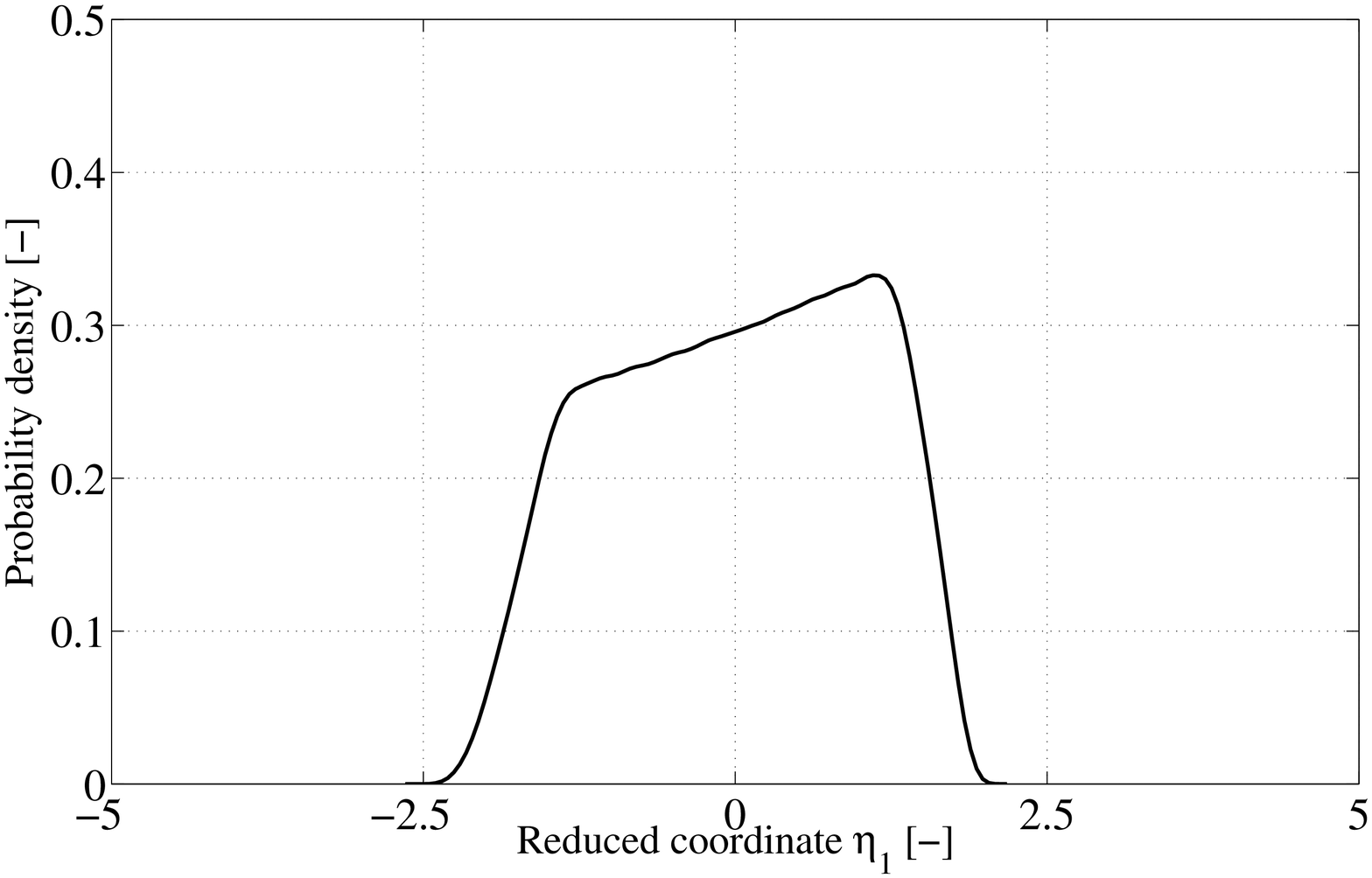}}
     \vfill
    \subfigure[Second marginal probability density function.]{\includegraphics[width=0.55\textwidth]{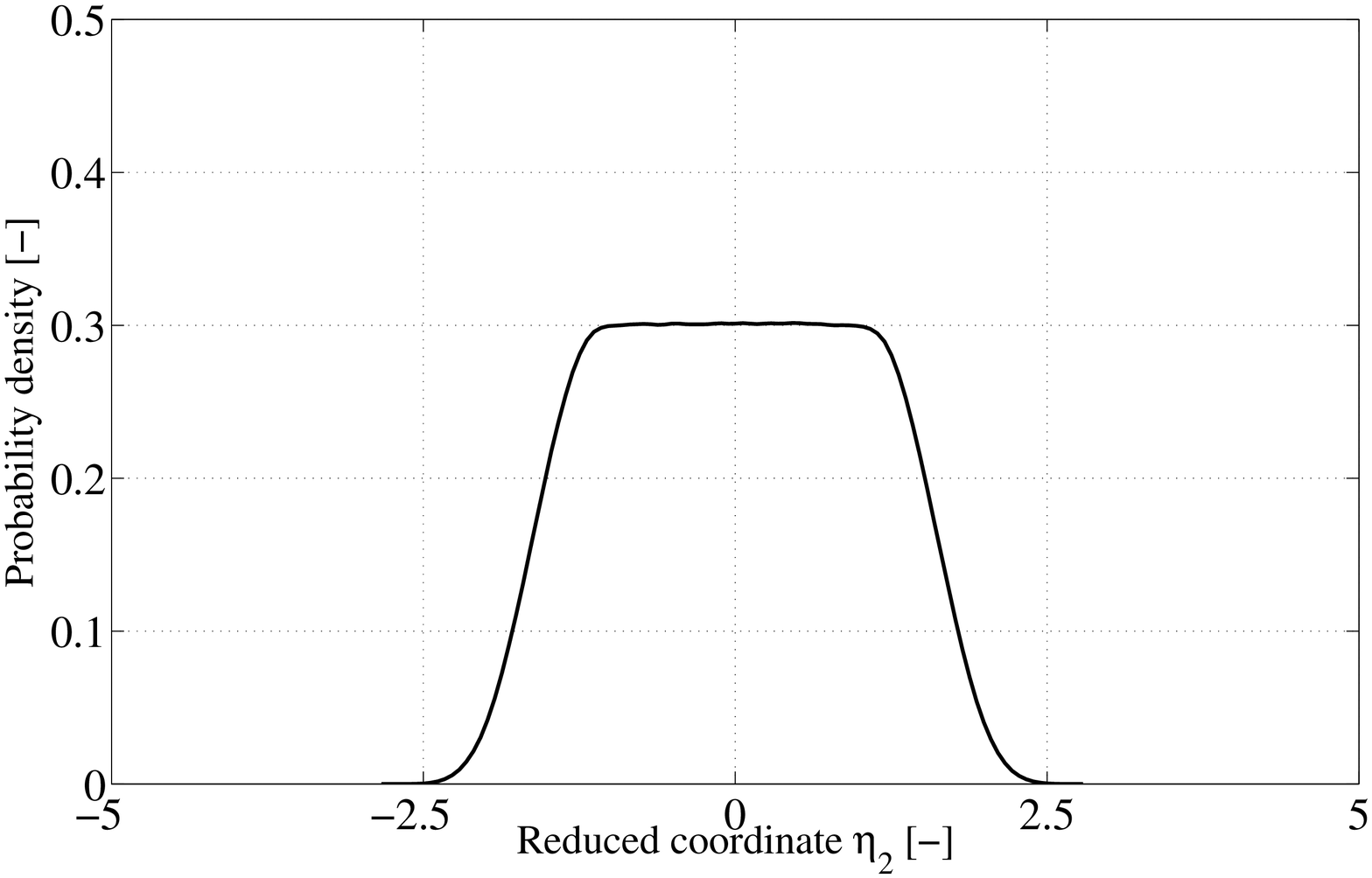}}
    \caption{PC-based simulation: probability distribution of the first and second reduced random variables of the reduced chaos expansion with random coefficients of the random temperature.}\label{fig:figure5}
  \end{center}
\end{figure}

Figure~\ref{fig:figure4} shows a few components of the reduced chaos expansion with random coefficients of the random temperature.
We can observe that the eigenvalues of the reduced chaos expansion with random coefficients of the random temperature (Fig.~\ref{fig:figure4}(b)) decay at a higher rate than those of the KL decomposition of the thermal transmittance random field (Fig.~\ref{fig:figure1ab}(b)).
This result is consistent with our earlier observation that the samples of the random temperature are smoother than those of the thermal transmittance random field.

Figure~\ref{fig:figure5} shows the joint and marginal probability density functions of the reduced random variables.
Clearly, the joint probability density function shows statistical dependence, and the marginal probability density functions are not ``labeled."

At iteration~$\ell=20$, a reduced chaos expansion with random coefficients obtained by retaining only $d=2$ terms was found to be sufficiently accurate to satisfy~(\ref{eq:criterion1}) for~$\epsilon_{1}=0.01$; thus, at this iteration, the measure transformation necessitated the construction of polynomial chaos and quadrature rules with respect to the joint probability distribution of the reduced random variables~$\boldsymbol{\eta}^{\ell,p}=(\eta^{\ell,p}_{1},\eta^{\ell,p}_{2})$ and the input random variables~$\boldsymbol{\zeta}=(\zeta_{1},\zeta_{2})$.

\begin{figure}[htp]
  \begin{center}
    \subfigure[$\Gamma_{00}^{\infty}$.]{\includegraphics[width=0.22\textwidth]{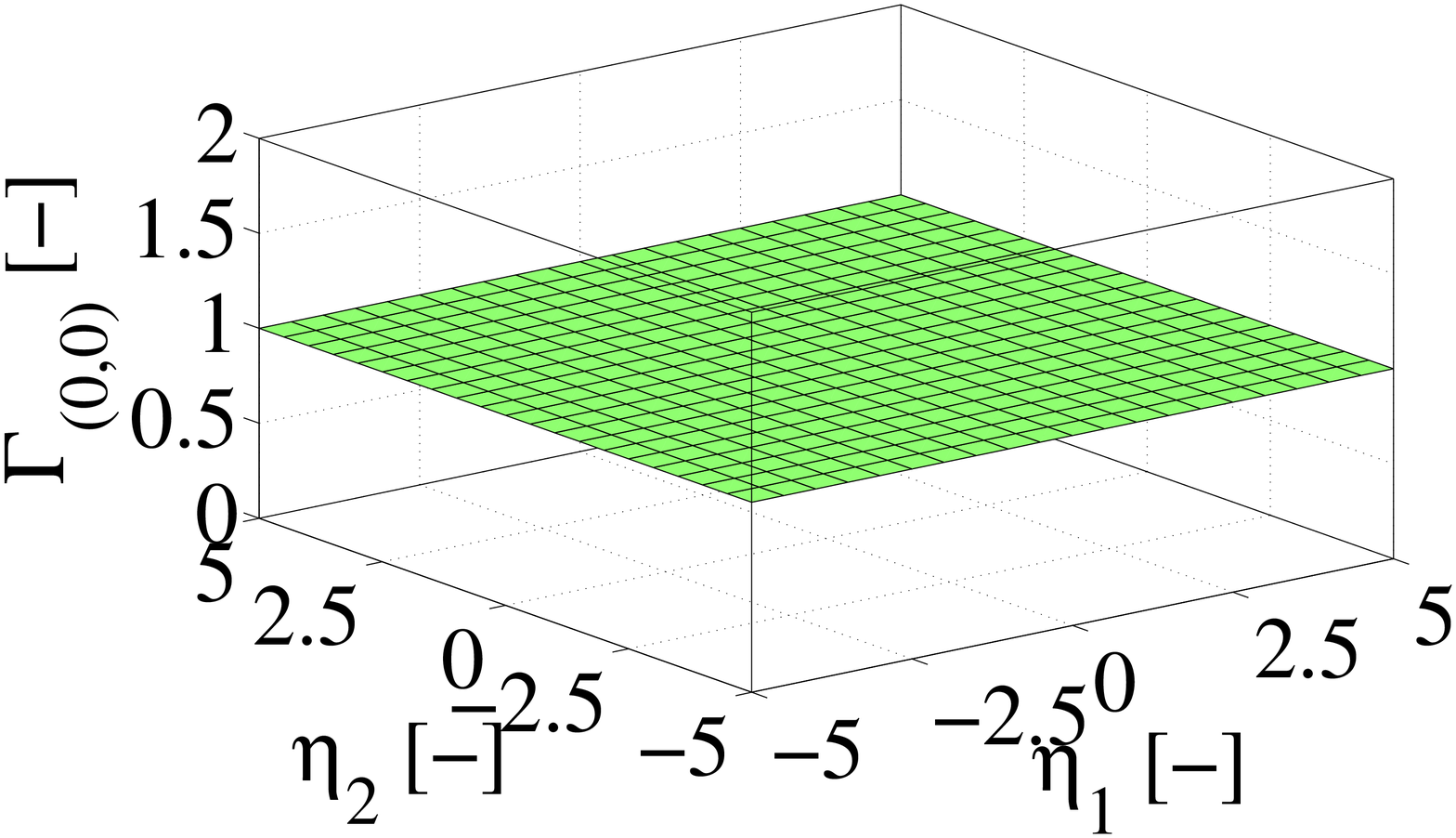}}
    \vfill
    \subfigure[$\Gamma_{10}^{\infty}$.]{\includegraphics[width=0.22\textwidth]{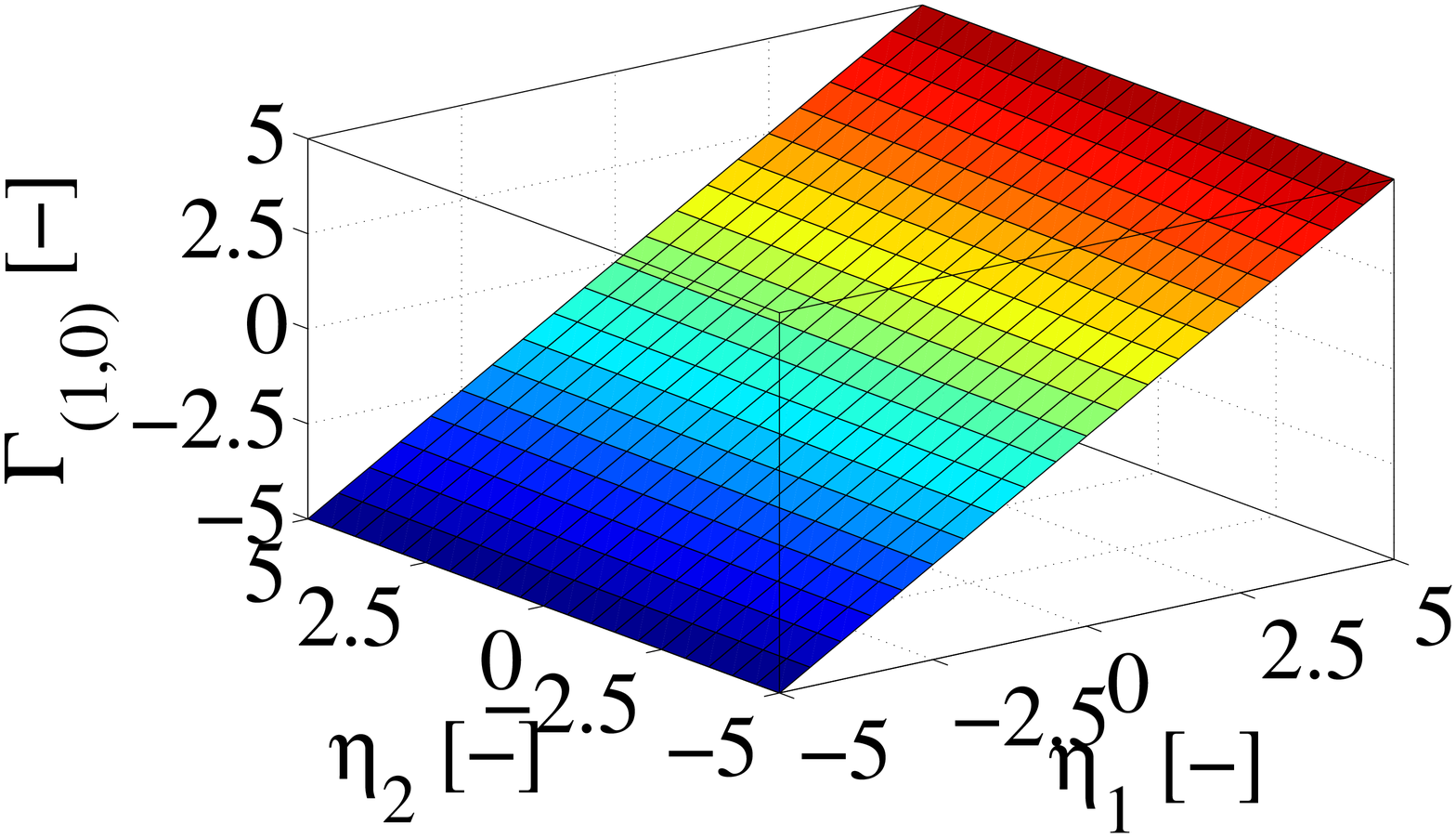}}
    \hspace{1mm}
    \subfigure[$\Gamma_{01}^{\infty}$.]{\includegraphics[width=0.22\textwidth]{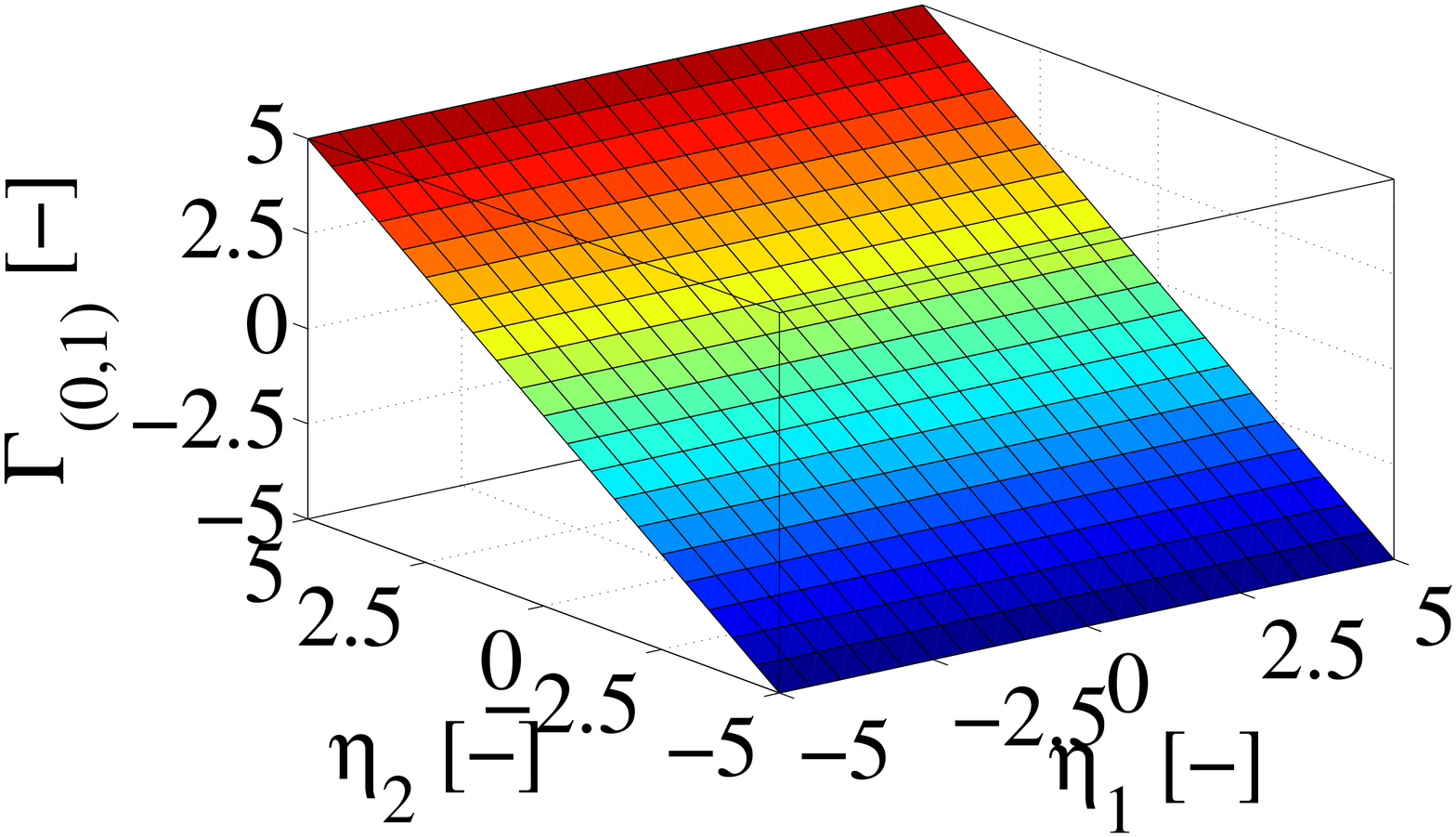}}
    \vfill
    \subfigure[$\Gamma_{20}^{\infty}$.]{\includegraphics[width=0.22\textwidth]{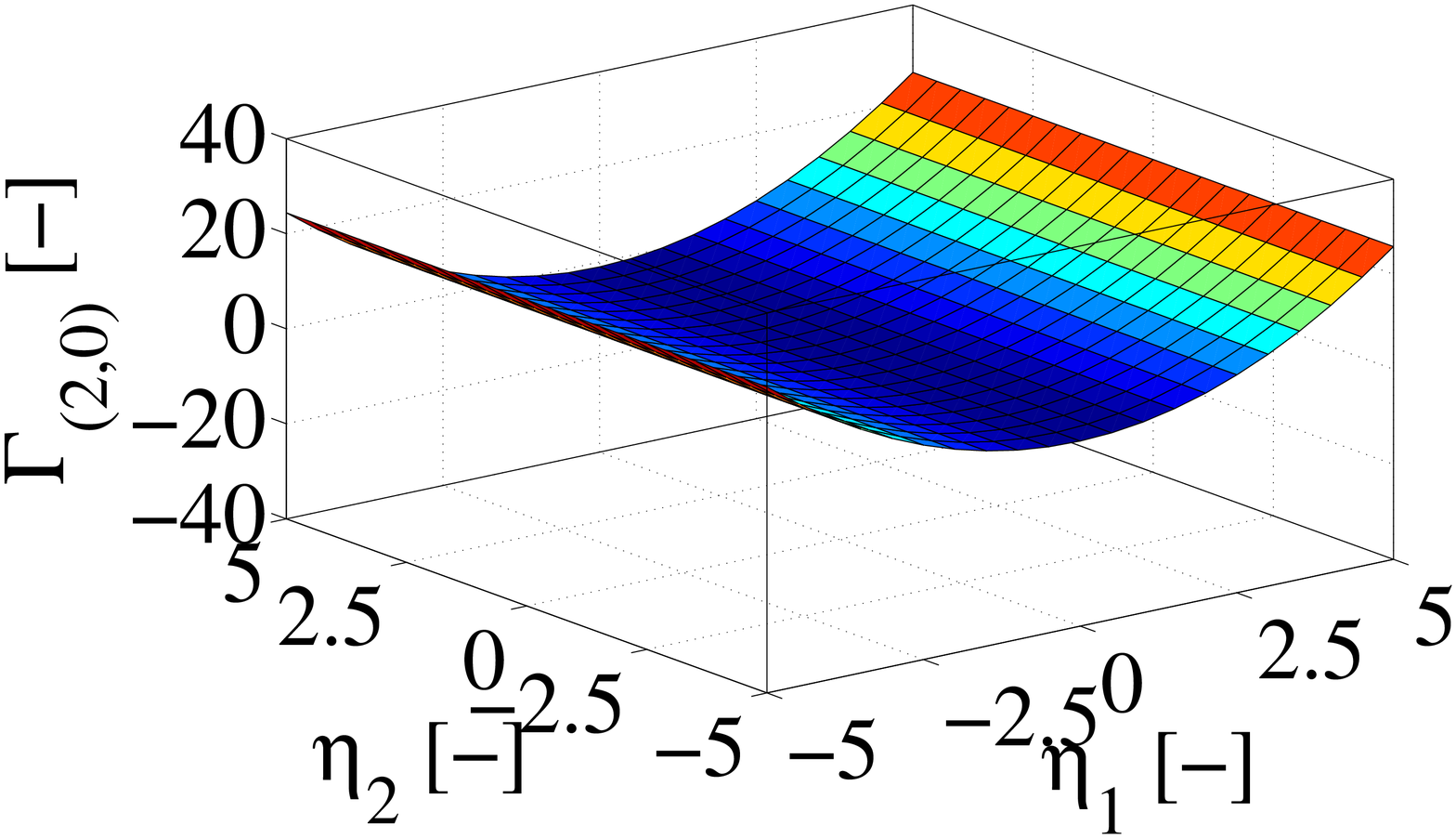}}
    \hspace{1mm}
    \subfigure[$\Gamma_{11}^{\infty}$.]{\includegraphics[width=0.22\textwidth]{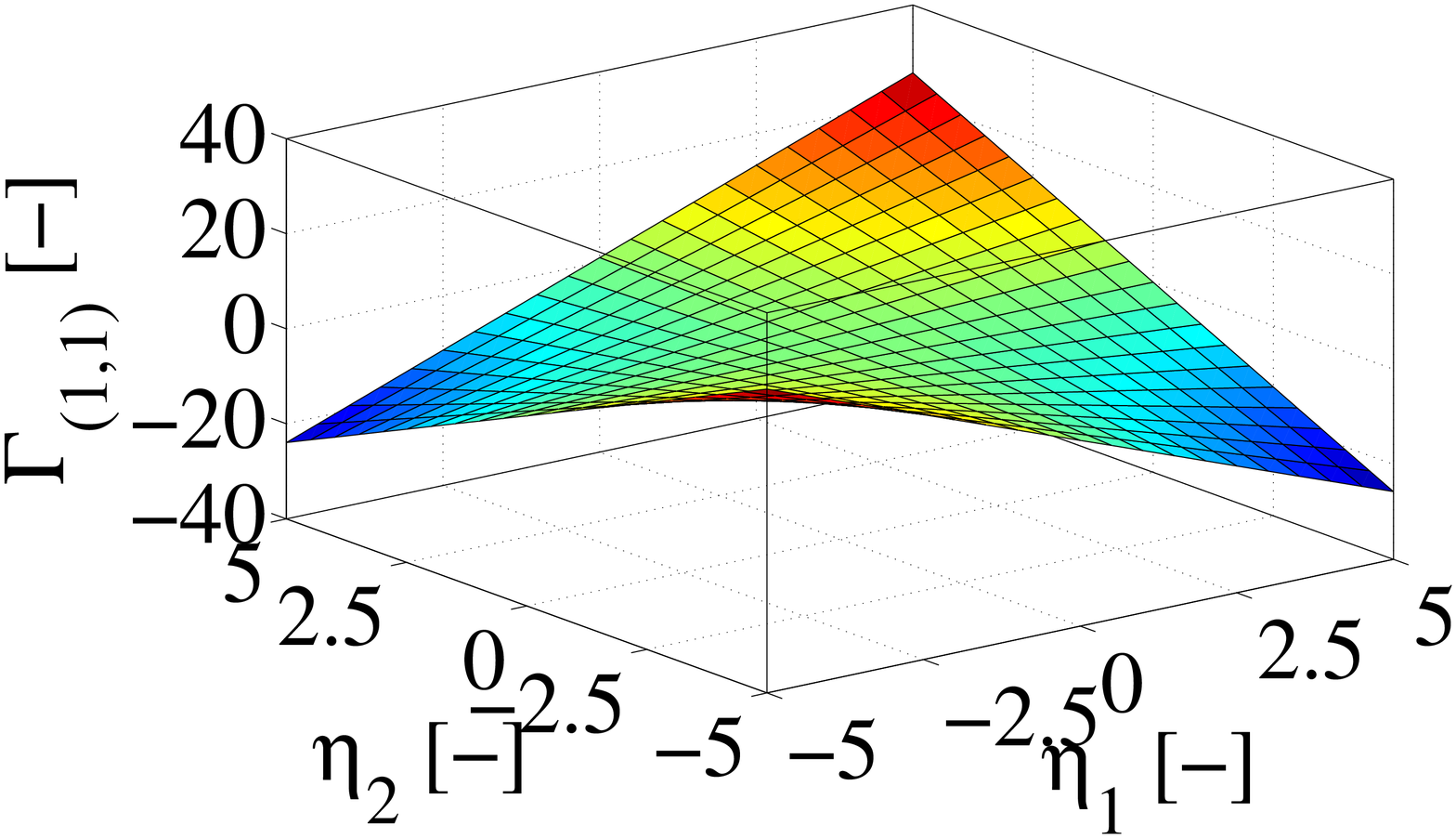}}
    \hspace{1mm}
    \subfigure[$\Gamma_{02}^{\infty}$.]{\includegraphics[width=0.22\textwidth]{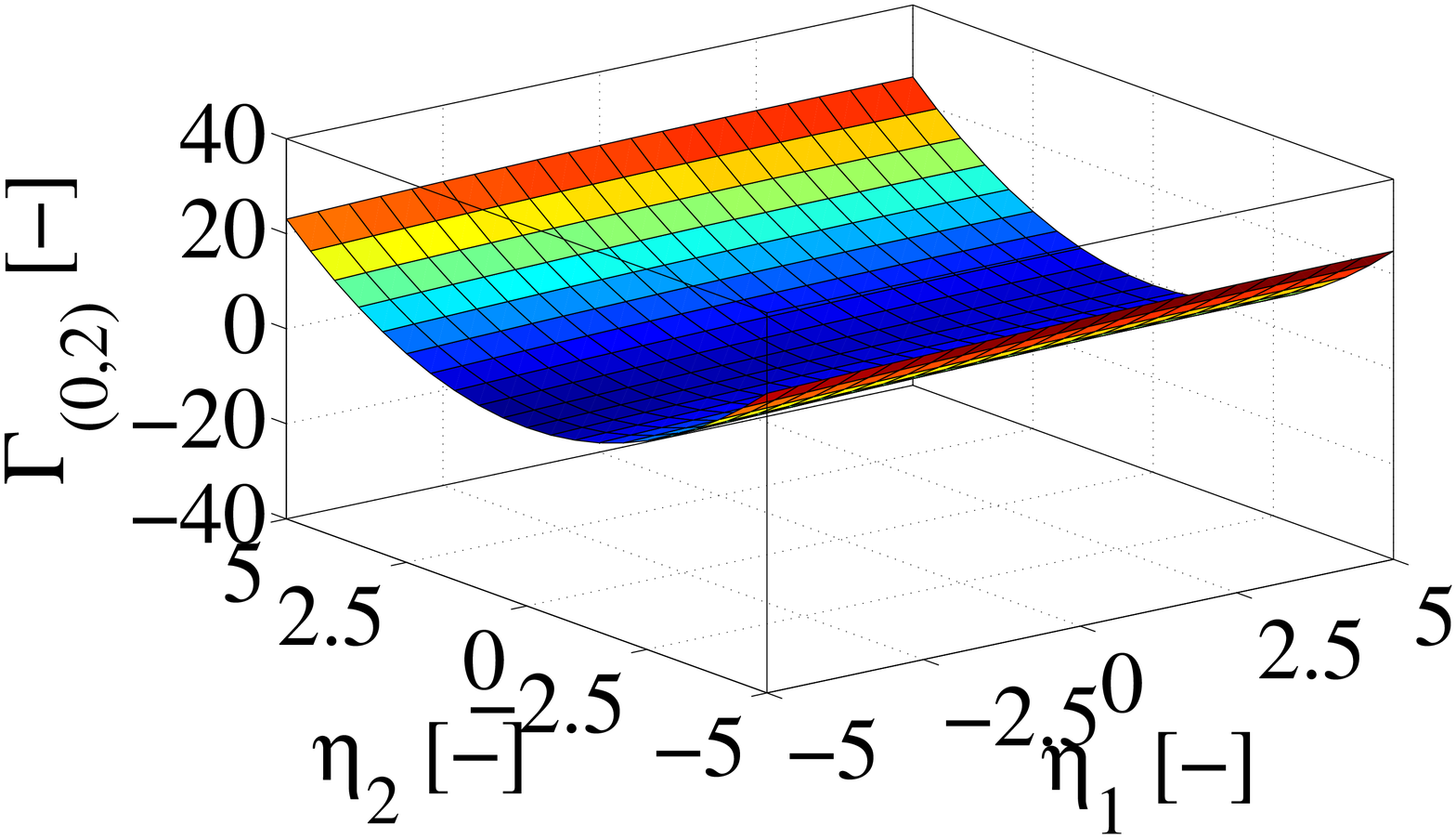}}
    \caption{PC-based simulation: computed polynomial chaos up to a total degree~$q$ of $2$.}\label{fig:figure6}
  \end{center}
\end{figure}

Figure~\ref{fig:figure6} illustrates the proposed computational construction of polynomial chaos with respect to~$P^{\ell}_{\boldsymbol{\eta}}$; specifically, the figure shows the polynomial chaos obtained up to a total degree~$q$ of $2$.
With reference to Sec.~\ref{sec:measuretransformationillustration}, the requisite polynomial chaos with respect to~$P_{\boldsymbol{\eta}}^{\ell}\times P_{\boldsymbol{\zeta}}$ are synthesized from the computed polynomial chaos thus obtained with respect to~$P^{\ell}_{\boldsymbol{\eta}}$ and the normalized Legendre polynomials with respect to~$P_{\boldsymbol{\zeta}}$ by tensorization.

\begin{figure}[htp]
  \begin{center}
    \subfigure[$\{\boldsymbol{\eta}_{k}^{\infty,\lambda},1\leq k\leq\nu^{\infty}_{\lambda}\}$ for $\lambda=1$.]{\includegraphics[width=0.39\textwidth]{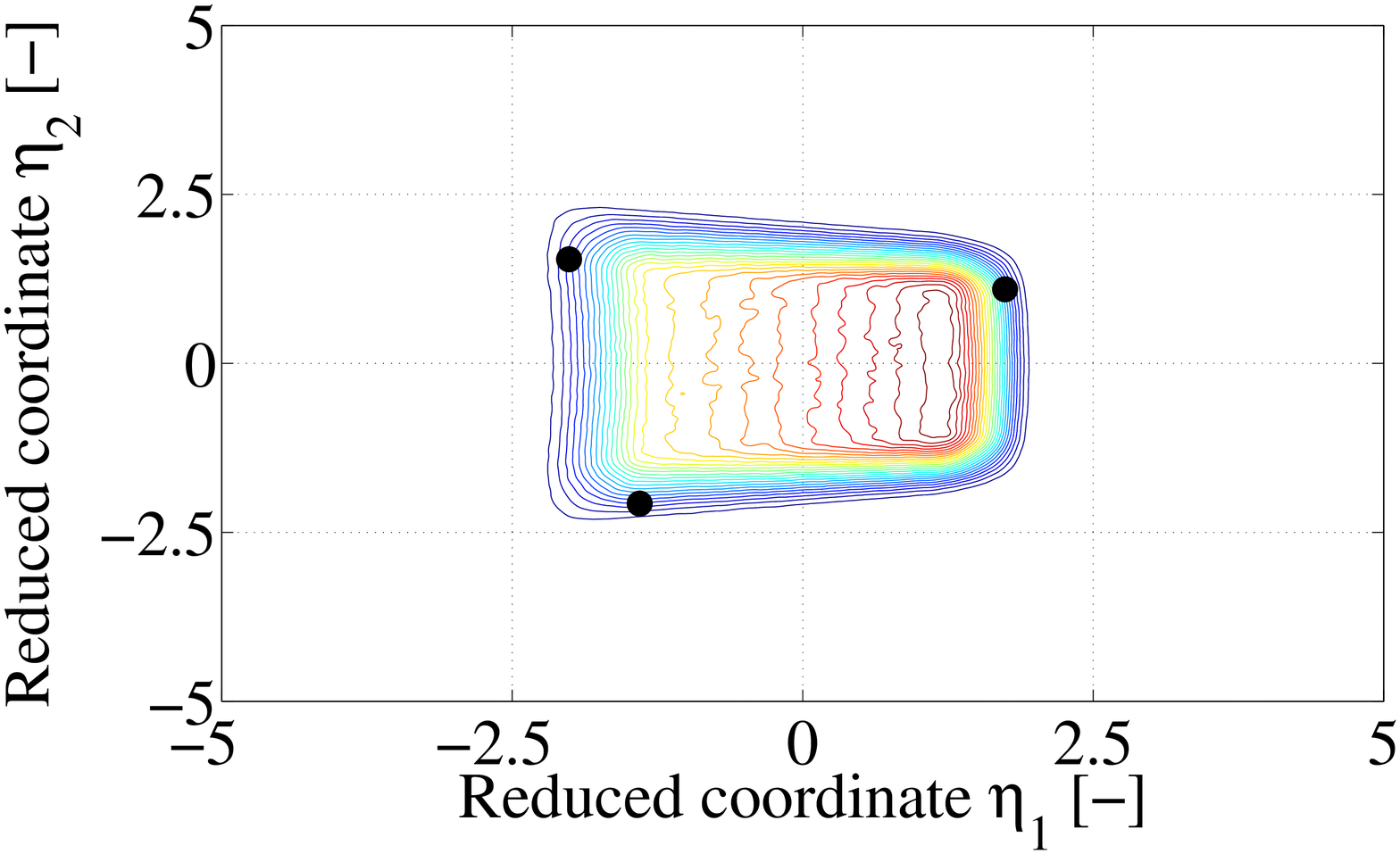}}
    \hfill
    \subfigure[$\{w_{k}^{\infty,\lambda},1\leq k\leq\nu^{\infty}_{\lambda}\}$ for $\lambda=1$.]{\includegraphics[width=0.39\textwidth]{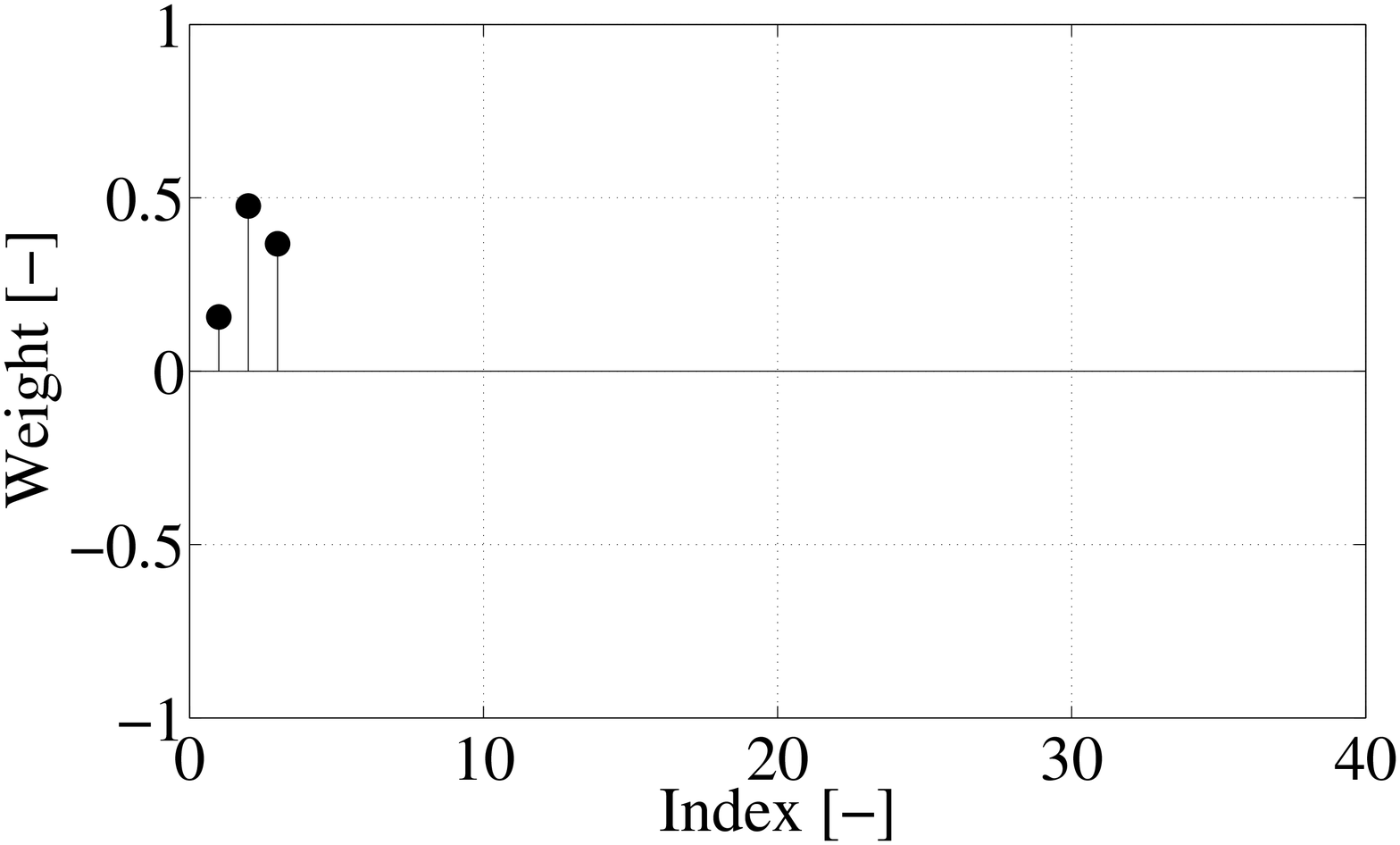}}
    \subfigure[$\{\boldsymbol{\eta}_{k}^{\infty,\lambda},1\leq k\leq\nu^{\infty}_{\lambda}\}$ for $\lambda=2$.]{\includegraphics[width=0.39\textwidth]{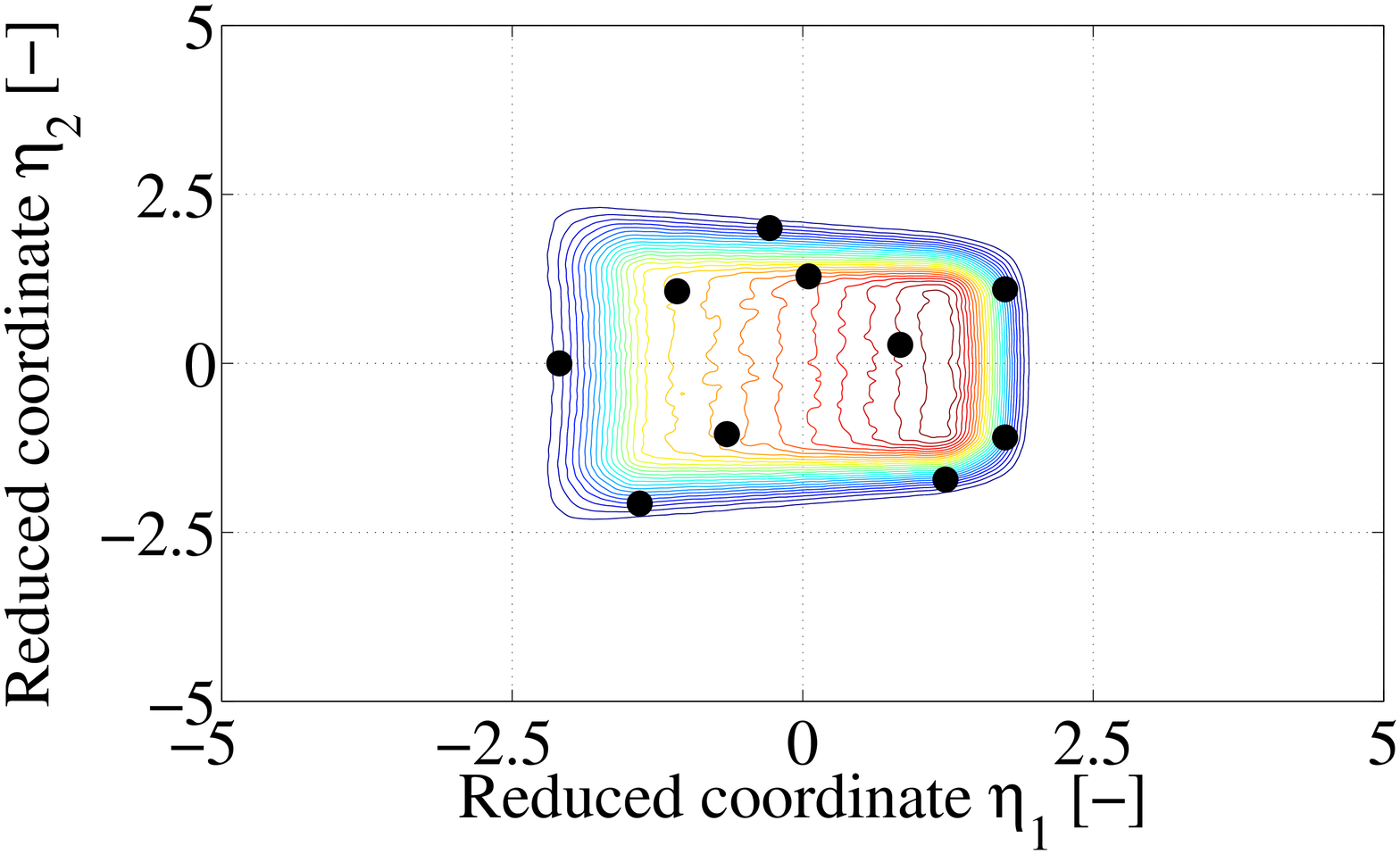}}
    \hfill
    \subfigure[$\{w_{k}^{\infty,\lambda},1\leq k\leq\nu^{\infty}_{\lambda}\}$ for $\lambda=2$.]{\includegraphics[width=0.39\textwidth]{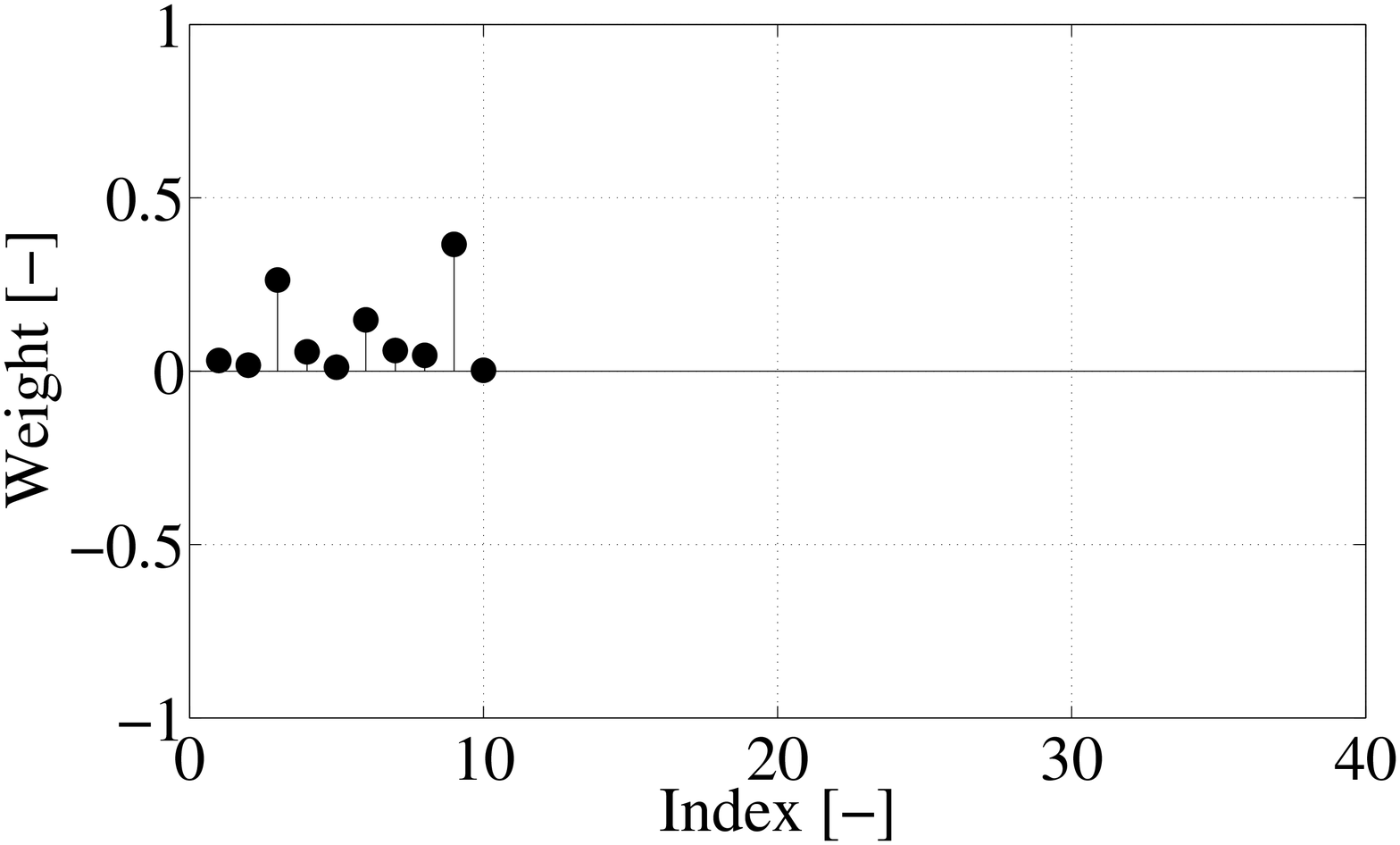}}
    \subfigure[$\{\boldsymbol{\eta}_{k}^{\infty,\lambda},1\leq k\leq\nu^{\infty}_{\lambda}\}$ for $\lambda=3$.]{\includegraphics[width=0.39\textwidth]{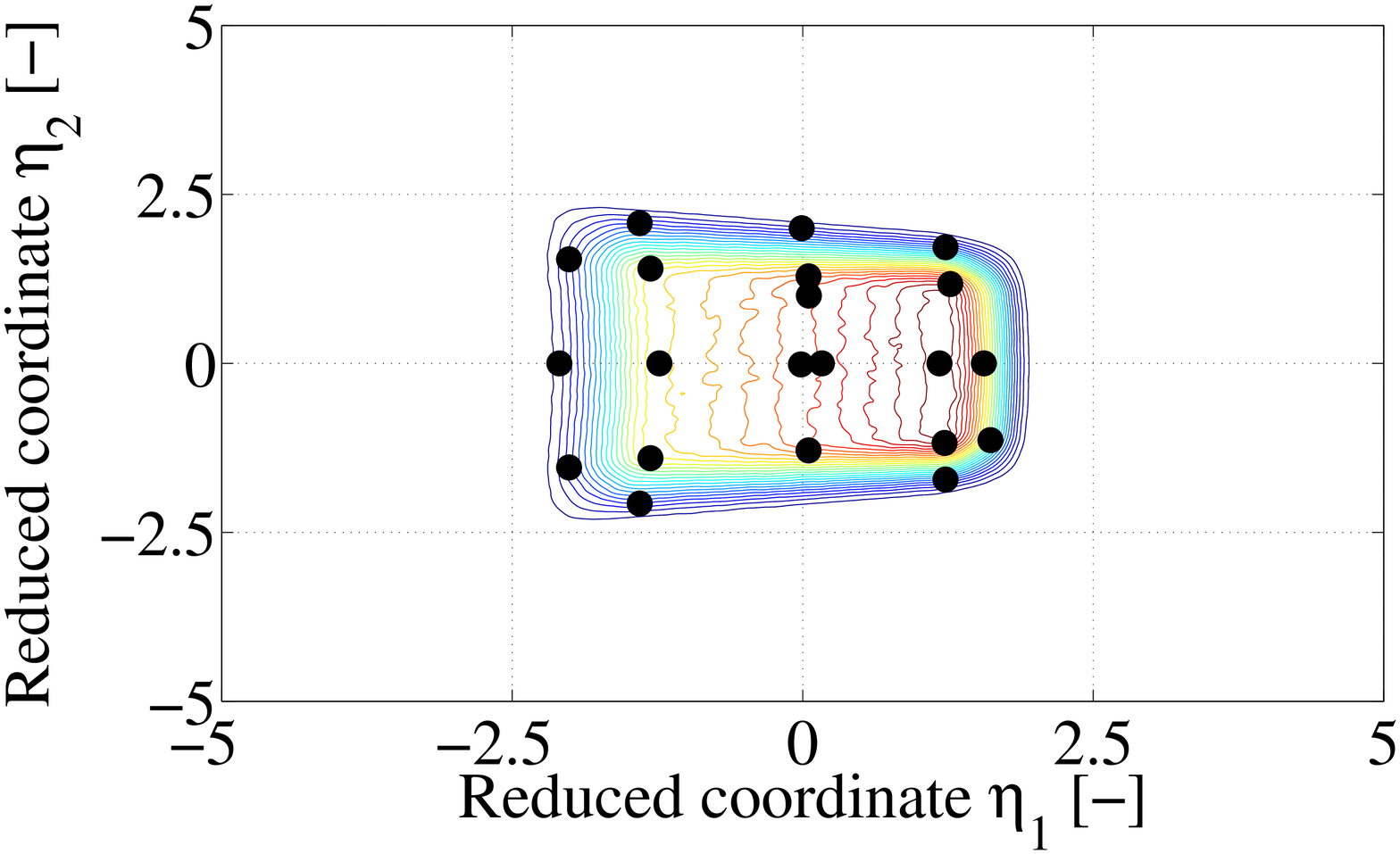}}
    \hfill
    \subfigure[$\{w_{k}^{\infty,\lambda},1\leq k\leq\nu^{\infty}_{\lambda}\}$ for $\lambda=3$.]{\includegraphics[width=0.39\textwidth]{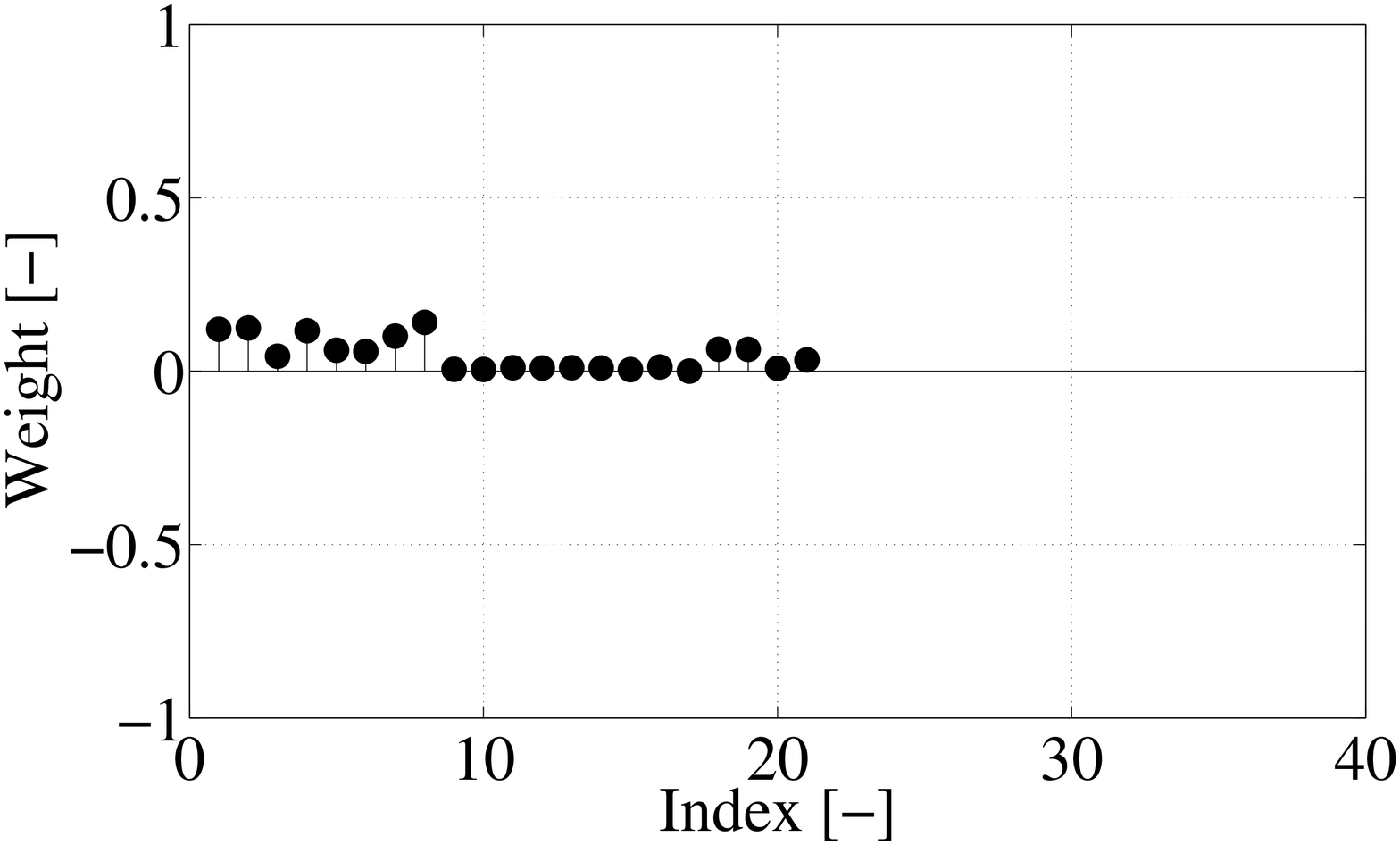}}
    \subfigure[$\{\boldsymbol{\eta}_{k}^{\infty,\lambda},1\leq k\leq\nu^{\infty}_{\lambda}\}$ for $\lambda=4$.]{\includegraphics[width=0.39\textwidth]{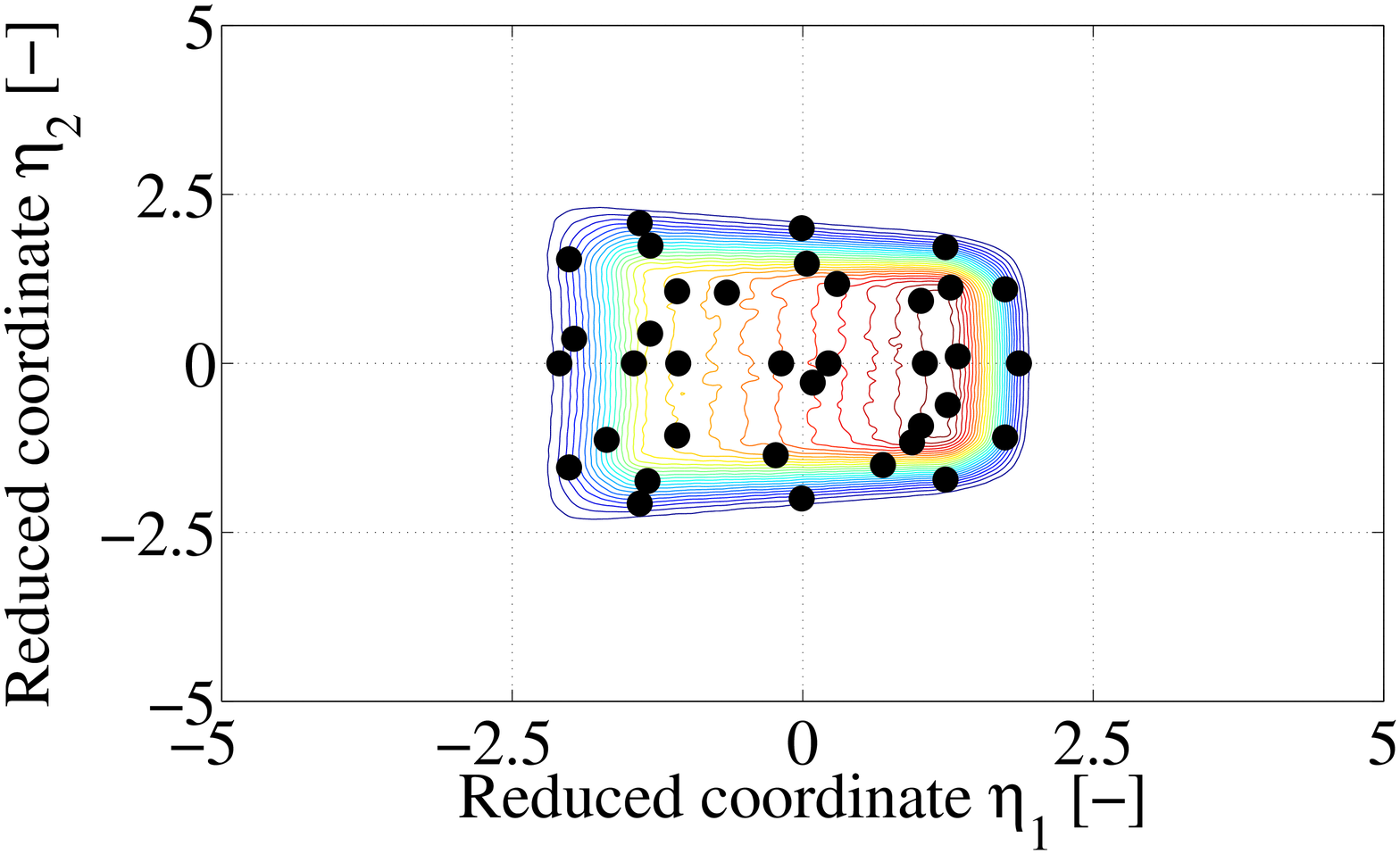}}
    \hfill
    \subfigure[$\{w_{k}^{\infty,\lambda},1\leq k\leq\nu^{\infty}_{\lambda}\}$ for $\lambda=4$.]{\includegraphics[width=0.39\textwidth]{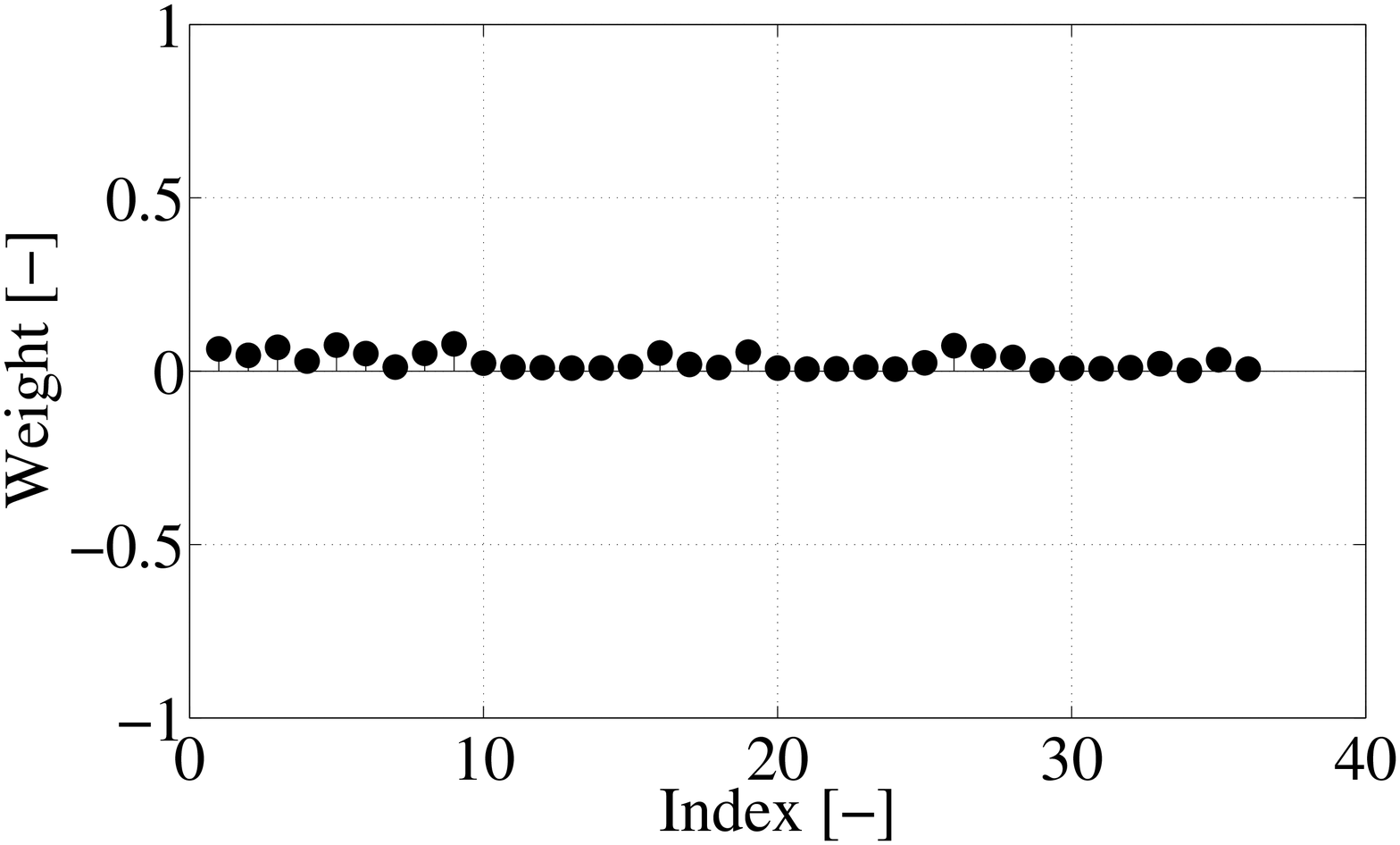}}
    \caption{PC-based simulation: computed quadrature rules up to a level~$\lambda$ of $4$.}\label{fig:figure7}
  \end{center}
\end{figure}

Figure~\ref{fig:figure7} illustrates the proposed computational construction of quadrature rules with respect to~$P^{\ell}_{\boldsymbol{\eta}}$; specifically, the figure shows the quadrature rules obtained up to a level~$\lambda$ of $4$.
This figure indicates that as $\lambda$ was increased, higher accuracy was required, and thus, a quadrature rule with more nodes and weights was systematically obtained.
With reference to Sec.~\ref{sec:measuretransformationillustration}, the requisite quadrature rules with respect to~$P_{\boldsymbol{\eta}}^{\ell}\times P_{\boldsymbol{\zeta}}$ are synthesized from the computed quadrature rules obtained with respect to~$P^{\ell}_{\boldsymbol{\eta}}$ and the fully tensorized Gauss-Legendre quadrature rules with respect to~$P_{\boldsymbol{\zeta}}$ by tensorization.

\begin{figure}[htp]
  \begin{center}
   \subfigure[Coefficients $\widehat{\boldsymbol{T}}{}^{\infty}_{\boldsymbol{\alpha}\boldsymbol{\beta}}$ indexed by multi-indices $(\boldsymbol{\alpha},\boldsymbol{\beta})$ with $|\boldsymbol{\beta}|=0$ (thick circles), $|\boldsymbol{\beta}|=1$ (thin circles), $|\boldsymbol{\beta}|=2$ (squares), $|\boldsymbol{\beta}|=3$ (diamonds), and $|\boldsymbol{\beta}|=4$ (pentagrams) of the chaos expansion $\widehat{\boldsymbol{T}}{}^{\infty,p}(\boldsymbol{\xi},\boldsymbol{\zeta})=\sum_{|\boldsymbol{\alpha}|+|\boldsymbol{\beta}|=0}^{p}\widehat{\boldsymbol{T}}{}^{\infty}_{\boldsymbol{\alpha}\boldsymbol{\beta}}\varphi_{\boldsymbol{\alpha}}(\boldsymbol{\xi})\psi_{\boldsymbol{\beta}}(\boldsymbol{\zeta})$ at $x=10\,\text{$[\text{cm}]$}$.]{\includegraphics[width=.8\textwidth]{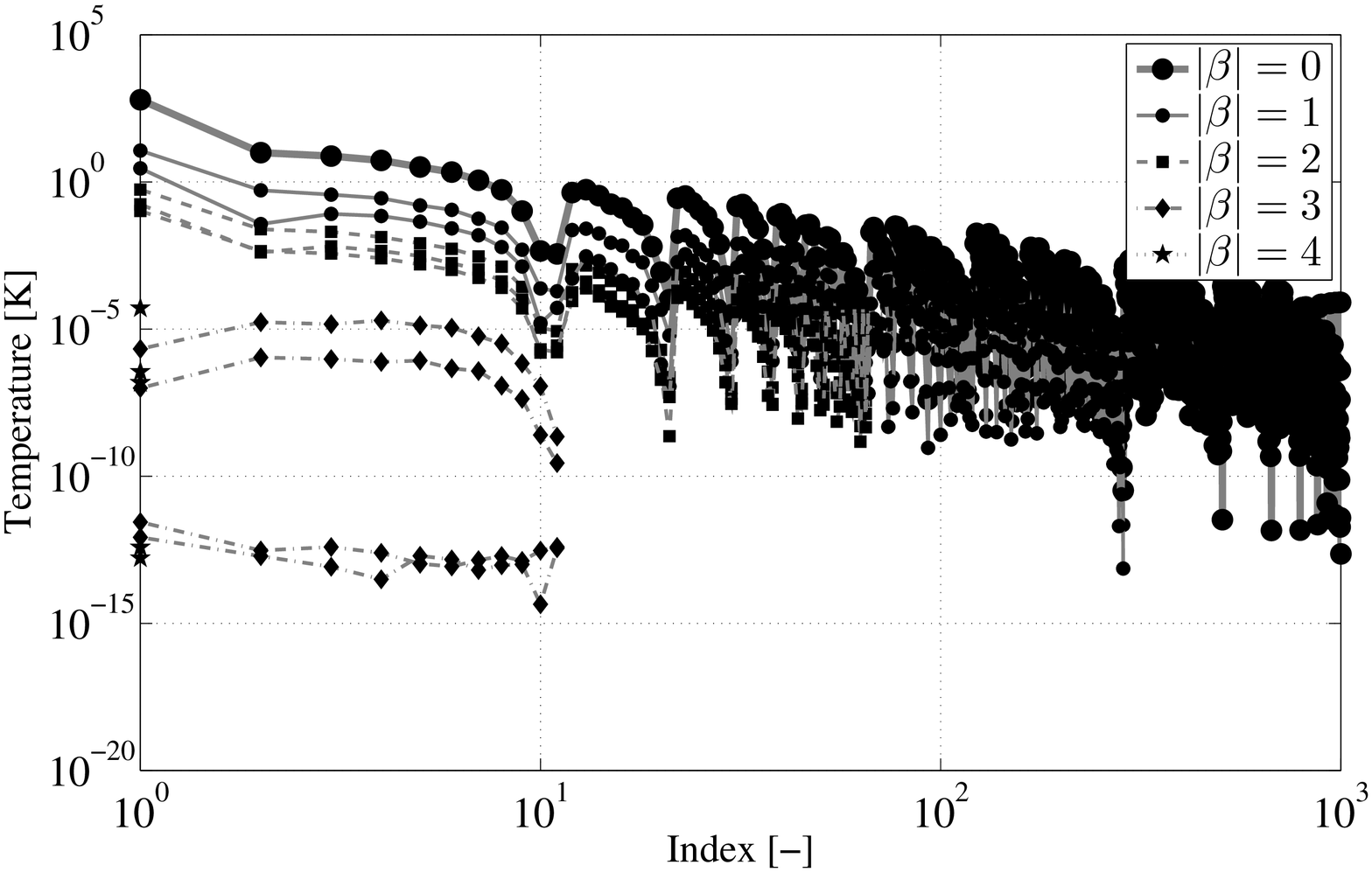}}
    \hfill
    \subfigure[Coefficients $\widehat{\boldsymbol{\Phi}}{}^{\infty}_{\boldsymbol{\beta}\boldsymbol{\gamma}}$ indexed by multi-indices $(\boldsymbol{\beta},\boldsymbol{\gamma})$ with $|\boldsymbol{\gamma}|=0$ (thick circles), $|\boldsymbol{\gamma}|=1$ (thin circles), and $|\boldsymbol{\gamma}|=2$ (squares) of the chaos expansion $\widehat{\boldsymbol{\Phi}}{}^{\infty,q}(\boldsymbol{\eta}^{\infty,p},\boldsymbol{\zeta})=\sum_{|\boldsymbol{\gamma}|+|\boldsymbol{\beta}|=0}^{q}\widehat{\boldsymbol{\Phi}}{}^{\infty}_{\boldsymbol{\gamma}\boldsymbol{\beta}}\Gamma_{\boldsymbol{\gamma}}^{\infty}(\boldsymbol{\eta}^{\infty,p})\psi_{\boldsymbol{\beta}}(\boldsymbol{\zeta})$ at $x=10\,\text{$[\text{cm}]$}$.]{\includegraphics[width=.8\textwidth]{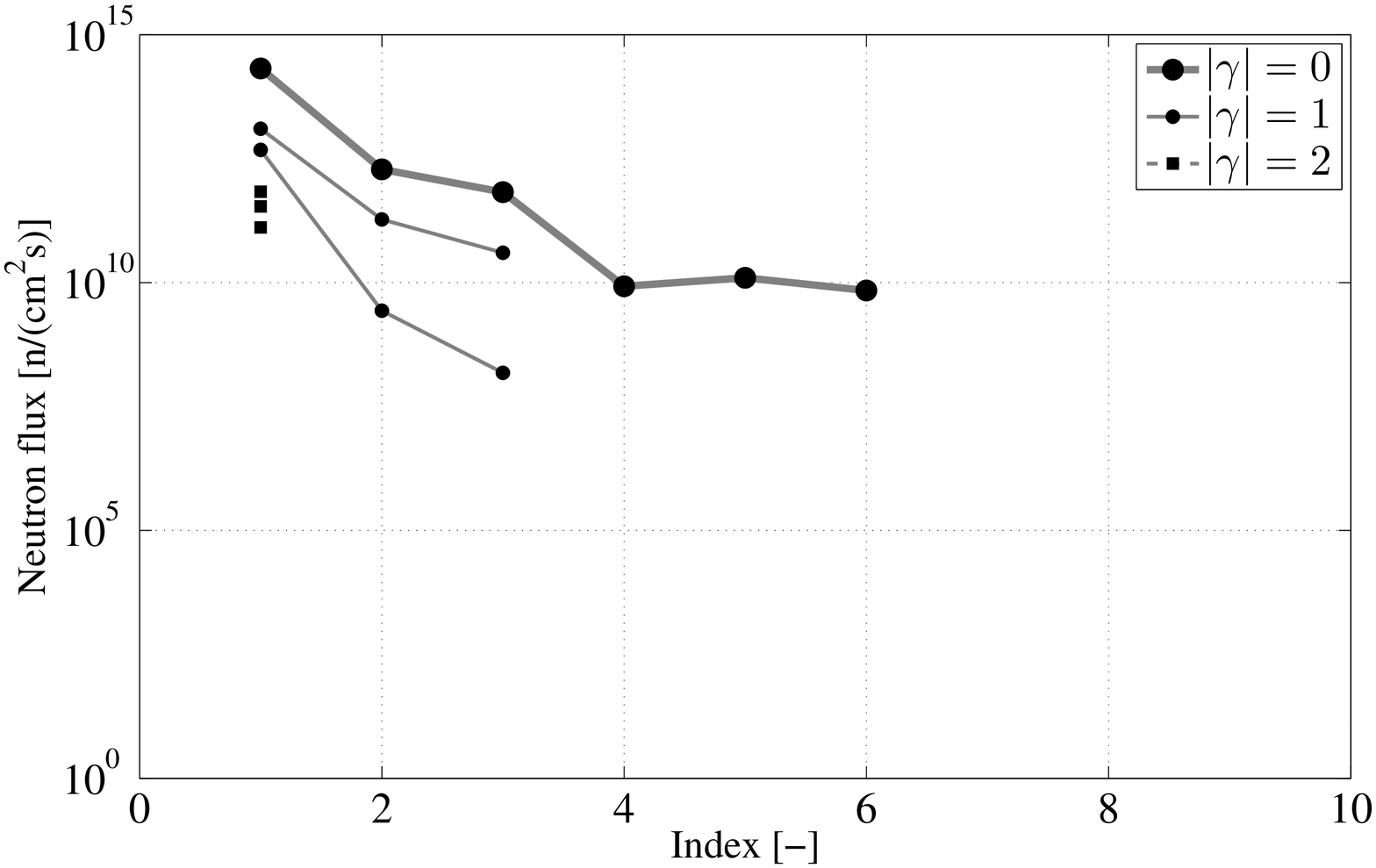}}
    \caption{PC-based simulation: coefficients of the solution at $x=10\,\text{$[\text{cm}]$}$.}\label{fig:figure8}
  \end{center}
\end{figure}

At iteration~$\ell=20$, a chaos expansion truncated at $q=2$ was found to be sufficiently accurate to satisfy~(\ref{eq:criterion2}) for~$\epsilon_{2}= 0.01$.
The representation of the random temperature by a chaos expansion of dimension $m+n=12$ and total degree $p=4$ requires~$1,820=16!/12!/4!$ terms; in contrast, the representation of the random neutron flux by a chaos expansion of dimension~$d+n=4$ and total degree~$q=2$ requires only~$15=6!/4!/2!$ terms.
The sparse-grid Gauss-Legendre quadrature rule of dimension $m+n=12$ and level $p+1=5$ used to compute the coefficients of the chaos expansion of the random temperature has $34,065$ nodes and weights; in contrast, the sparse-grid quadrature rule of dimension $d+n=4$ and level $q+2=4$ used to compute the coefficients of the chaos expansion of the random neutron flux has only $346$ nodes and weights.
Figure~\ref{fig:figure8} shows a few coefficients of the solution.

\begin{figure}[htp]
  \begin{center}
    \subfigure[Temperature.]{\includegraphics[width=0.8\textwidth]{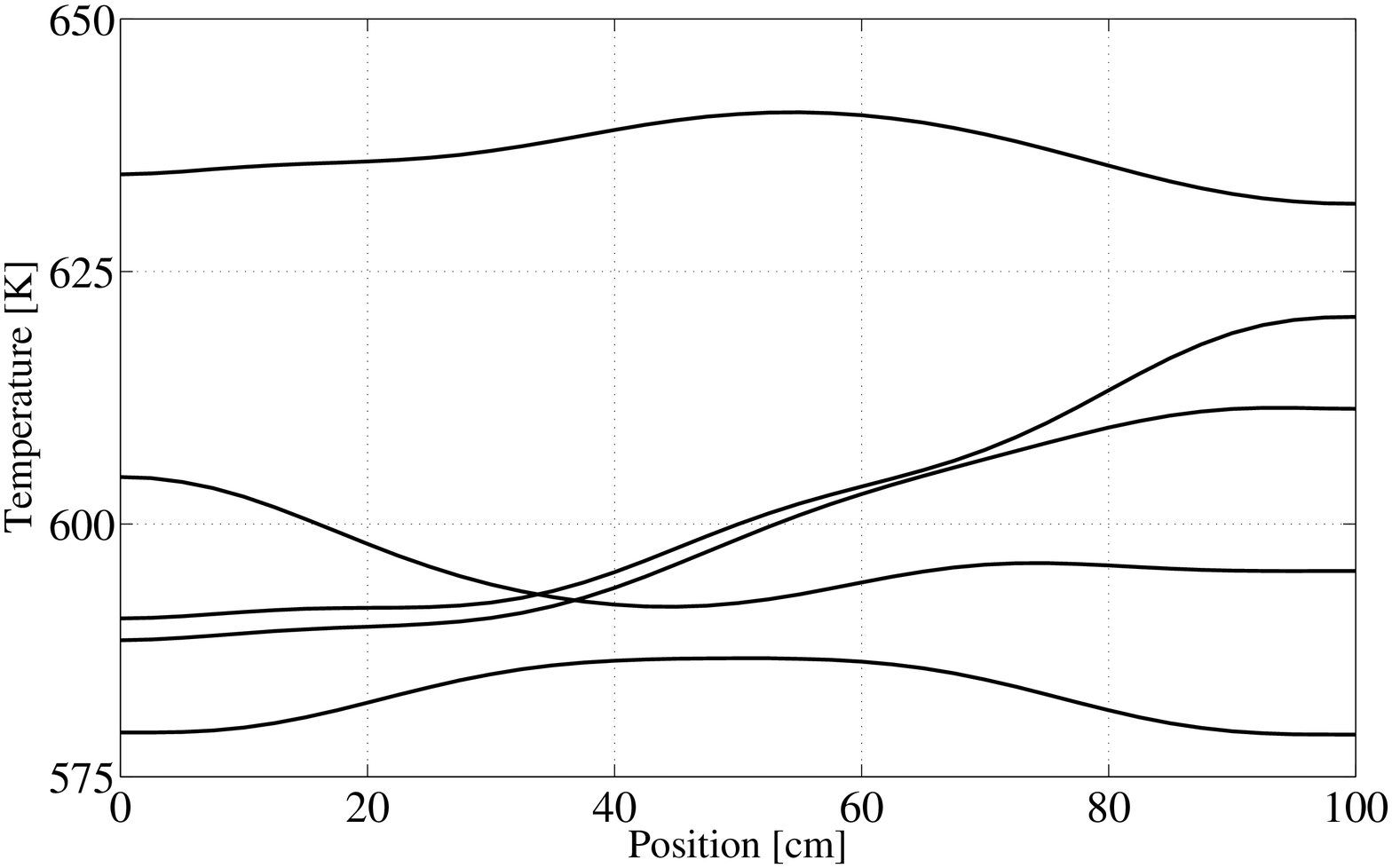}}
    \hfill
    \subfigure[Neutron flux.]{\includegraphics[width=0.8\textwidth]{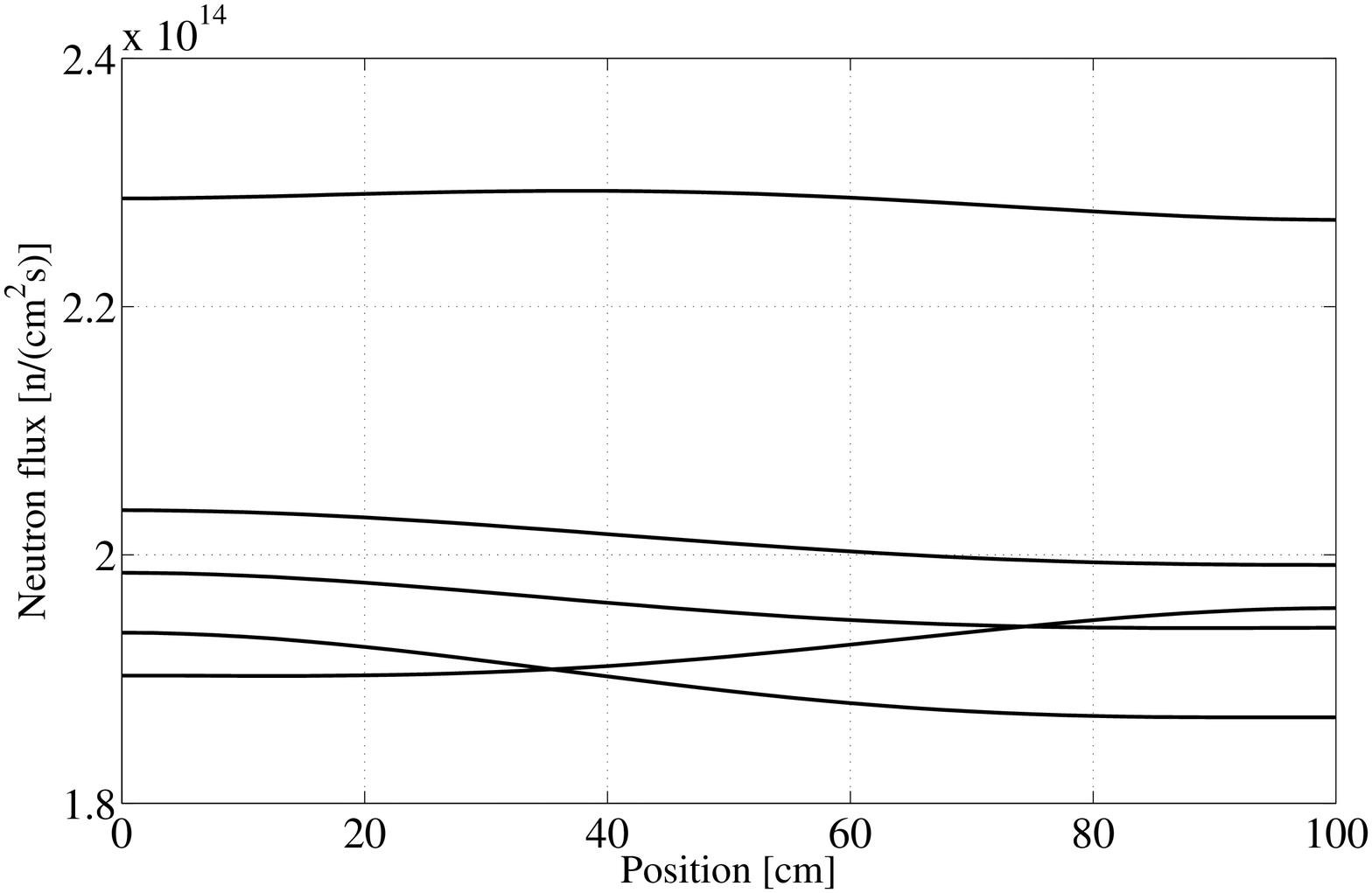}}
    \caption{PC-based simulation: five samples of the solution.}\label{fig:figure9}
  \end{center}
\end{figure}

Figure~\ref{fig:figure9} shows a few samples of the random temperature and random neutron flux deduced from the chaos expansions obtained as the output of the solution algorithm.
The samples of the input random variables used to synthesize the samples of the random temperature and random neutron flux shown in Fig.~\ref{fig:figure9} were identical to those used to generate the samples shown in Fig.~\ref{fig:figure2}.
The similarity of the samples in Figs.~\ref{fig:figure2} and~\ref{fig:figure9} indicates that the surrogate model based on polynomial chaos not only provides an accurate \textit{global} representation of the multiphysics model but is also capable of accurately reproducing a \textit{sample-wise} response.

\subsection{Sensitivity analysis of the random temperature and random neutron flux}
The chaos expansions in~(\ref{eq:PCTphi}) facilitate the following decompositions of the variances of the random temperature and random neutron flux:
\begin{align}
V^{T}&=V^{T}_{\boldsymbol{\xi}}+V^{T}_{\boldsymbol{\zeta}}+V^{T}_{(\boldsymbol{\xi},\boldsymbol{\zeta})},&V^{\Phi}&=V^{\Phi}_{\boldsymbol{\xi}}+V^{\Phi}_{\boldsymbol{\zeta}}+V^{\Phi}_{(\boldsymbol{\xi},\boldsymbol{\zeta})},\\
\notag V^{T}&=\sum_{|\boldsymbol{\alpha}|+|\boldsymbol{\beta}|=1}^{p}\|\widehat{T}{}^{\infty}_{\boldsymbol{\alpha}\boldsymbol{\beta}}\|_{\boldsymbol{W}}^{2},&V^{\Phi}&=\sum_{|\boldsymbol{\gamma}|+|\boldsymbol{\beta}|=1}^{q}\|\widehat{\Phi}{}^{\infty}_{\boldsymbol{\gamma}\boldsymbol{\beta}}\|_{\boldsymbol{W}}^{2},\\
\notag V_{\boldsymbol{\xi}}^{T}&=\sum_{\;\;\;\,|\boldsymbol{\alpha}|=1\;\;\;\,}^{p}\|\widehat{T}{}^{\infty}_{\boldsymbol{\alpha}\boldsymbol{0}}\|_{\boldsymbol{W}}^{2},&V_{\boldsymbol{\xi}}^{\Phi}&=\sum_{\;\;\;\,|\boldsymbol{\gamma}|=1\;\;\;\,}^{q}\|\widehat{\Phi}{}^{\infty}_{\boldsymbol{\gamma}\boldsymbol{0}}\|_{\boldsymbol{W}}^{2},\\
\notag V_{\boldsymbol{\zeta}}^{T}&=\sum_{\;\;\;\,|\boldsymbol{\beta}|=1\;\;\;\,}^{p}\|\widehat{T}{}^{\infty}_{\boldsymbol{0}\boldsymbol{\beta}}\|_{\boldsymbol{W}}^{2},&V_{\boldsymbol{\zeta}}^{\Phi}&=\sum_{\;\;\;\,|\boldsymbol{\beta}|=1\;\;\;\,}^{q}\|\widehat{\Phi}{}^{\infty}_{\boldsymbol{0}\boldsymbol{\beta}}\|_{\boldsymbol{W}}^{2},\\
\notag V_{(\boldsymbol{\xi},\boldsymbol{\zeta})}^{T}&=\sum_{\substack{|\boldsymbol{\alpha}|+|\boldsymbol{\beta}|=1\\\boldsymbol{\alpha}\neq\boldsymbol{0},\boldsymbol{\beta}\neq\boldsymbol{0}}}^{p}\|\widehat{T}{}^{\infty}_{\boldsymbol{\alpha}\boldsymbol{\beta}}\|_{\boldsymbol{W}}^{2},&V_{(\boldsymbol{\xi},\boldsymbol{\zeta})}^{\Phi}&=\sum_{\substack{|\boldsymbol{\gamma}|+|\boldsymbol{\beta}|=1\\\boldsymbol{\gamma}\neq\boldsymbol{0},\boldsymbol{\beta}\neq\boldsymbol{0}}}^{q}\|\widehat{\Phi}{}^{\infty}_{\boldsymbol{\gamma}\boldsymbol{\beta}}\|_{\boldsymbol{W}}^{2}.
\end{align}
Here, $V_{\boldsymbol{\xi}}^{T}$ and~$V_{\boldsymbol{\xi}}^{\Phi}$ are the sums of the variances of those terms in the chaos expansions of the random temperature and random neutron flux, respectively, which depend only on the input random variables $\boldsymbol{\xi}$ that describe the uncertainty in the parameters of the heat subproblem; hence, following the approach given in~\citep{oakley2004,sudret2008,crestaux2009}, $V_{\boldsymbol{\xi}}^{T}$ and~$V_{\boldsymbol{\xi}}^{\Phi}$ can be interpreted as the portions contributed by the uncertainty in the parameters of the heat subproblem to the variances of the random temperature and random neutron flux, respectively.  
Conversely, $V_{\boldsymbol{\zeta}}^{T}$ and~$V_{\boldsymbol{\zeta}}^{\Phi}$ are the sums of the variances of those terms in the chaos expansions of the random temperature and random neutron flux which depend only on the input random variables $\boldsymbol{\zeta}$ that describe the uncertainty in the parameters of the neutronics subproblem, and they can be interpreted as the portions contributed by the uncertainty in the parameters of the neutronics subproblem to the variances of the random temperature and random neutron flux, respectively. 
Lastly, $V_{(\boldsymbol{\xi},\boldsymbol{\zeta})}^{T}$ and~$V_{(\boldsymbol{\xi},\boldsymbol{\zeta})}^{\Phi}$ are the portions contributed by the interaction of the uncertainties in the parameters of the heat subproblem and those in the parameters of the neutronics subproblem to the variances of the random temperature and random neutron flux, respectively. 

\begin{table}
\begin{center}
\begin{tabular}{|c|c|c|}
\hline
& temperature  & neutron flux \\
\hline
\parbox{8cm}{\centering main effect\\ of uncertainty in data of heat subproblem} & $\frac{V^{T}_{\boldsymbol{\xi}}}{V^{T}}=55.45\%$ & $\frac{V^{\Phi}_{\boldsymbol{\xi}}}{V^{\Phi}}=2.15\%$ \\
\hline
\parbox{8cm}{\centering main effect\\ of uncertainty in data of neutronics subproblem} & $\frac{V^{T}_{\boldsymbol{\zeta}}}{V^{T}}=44.40\%$ & $\frac{V^{\Phi}_{\boldsymbol{\zeta}}}{V^{\Phi}}=97.83\%$\\
\hline
\parbox{8cm}{\centering interaction\\of uncertainty in data of heat subproblem \\and uncertainty in data of neutronics subproblem} & $\frac{V^{T}_{(\boldsymbol{\xi},\boldsymbol{\zeta})}}{V^{T}}=0.02\%$ & $\frac{V^{\Phi}_{(\boldsymbol{\xi},\boldsymbol{\zeta})}}{V^{\Phi}}=0.08\%$\\
\hline
\end{tabular}
\end{center}
\caption{Decomposition of the variances of the random temperature and random neutron flux.}\label{table:table1}
\end{table}

We obtained the values given in Table~\ref{table:table1}; these values indicate that the uncertainty in the parameters of the heat subproblem predominantly drives the uncertainty in the temperature ($55.45\%$) and that the uncertainty in the parameters of the neutronics subproblem predominantly drives the uncertainty in the neutron flux ($97.83\%$).
Further, we can observe that the uncertainty in the parameters of the heat subproblem is less important in driving the uncertainty in the neutron flux than the uncertainty in the parameters of the neutronics subproblem is in driving the uncertainty in the temperature ($2.15\%$ versus $44.40\%$).
This observation can be explained as a consequence of the coupling mechanism adopted in this work; with reference to~(\ref{eq:couplingmechanism}), this coupling mechanism involves the division of the temperature by a large reference temperature as it feeds into the coefficients of the neutronics subproblem, thus diminishing the impact that fluctuations in the temperature have on the neutron flux.

\subsection{Convergence analysis}

\begin{figure}[htp]
  \begin{center}
    \hspace{-13mm}\subfigure[Reduced dimension $d$.]{\includegraphics[width=0.55\textwidth]{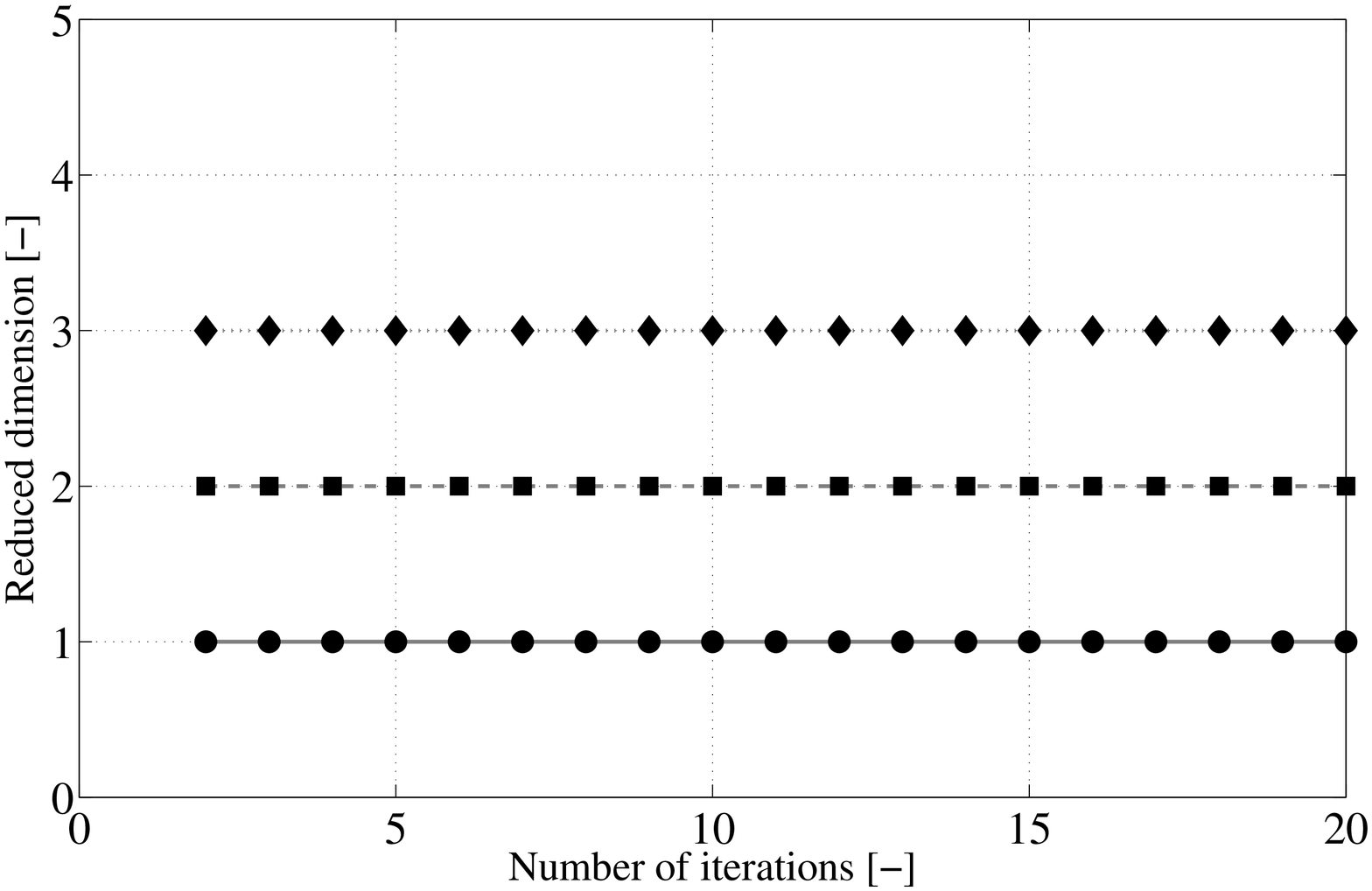}}
    \vfill
    \hspace{-13mm}\subfigure[$\ell\mapsto\sqrt{\frac{1}{MC}\sum_{k=1}^{MC}\|\boldsymbol{T}^{\ell}(\boldsymbol{\xi}_{k},\boldsymbol{\zeta}_{k})-\widehat{\boldsymbol{T}}{}^{\ell,p}(\boldsymbol{\xi}_{k},\boldsymbol{\zeta}_{k})\|_{\boldsymbol{W}}^{2}}\Big/\sqrt{\frac{1}{MC}\sum_{k=1}^{MC}\|\boldsymbol{T}^{\infty}(\boldsymbol{\xi}_{k},\boldsymbol{\zeta}_{k})\|_{\boldsymbol{W}}^{2}}$.]{\includegraphics[width=0.55\textwidth]{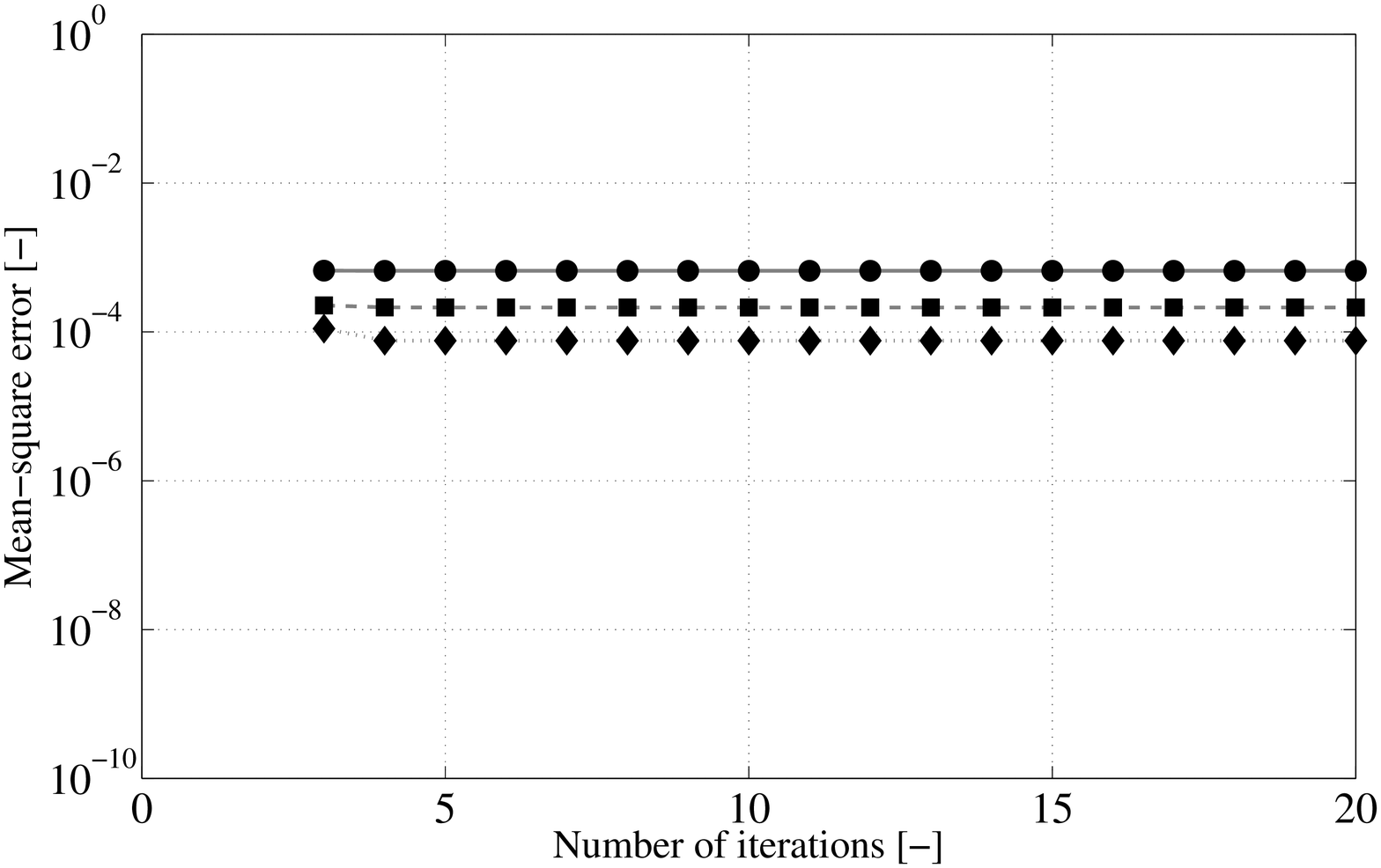}}
    \vfill
    \hspace{-13mm}\subfigure[$\ell\mapsto\sqrt{\frac{1}{MC}\sum_{k=1}^{MC}\|\boldsymbol{\Phi}^{\ell}(\boldsymbol{\xi}_{k},\boldsymbol{\zeta}_{k})-\widehat{\boldsymbol{\Phi}}{}^{\ell,q}(\boldsymbol{\eta}^{\ell,p}(\boldsymbol{\xi}_{k}),\boldsymbol{\zeta}_{k})\|_{\boldsymbol{W}}^{2}}\Big/\sqrt{\frac{1}{MC}\sum_{k=1}^{MC}\|\boldsymbol{\Phi}^{\infty}(\boldsymbol{\xi}_{k},\boldsymbol{\zeta}_{k})\|_{\boldsymbol{W}}^{2}}$.]{\includegraphics[width=0.55\textwidth]{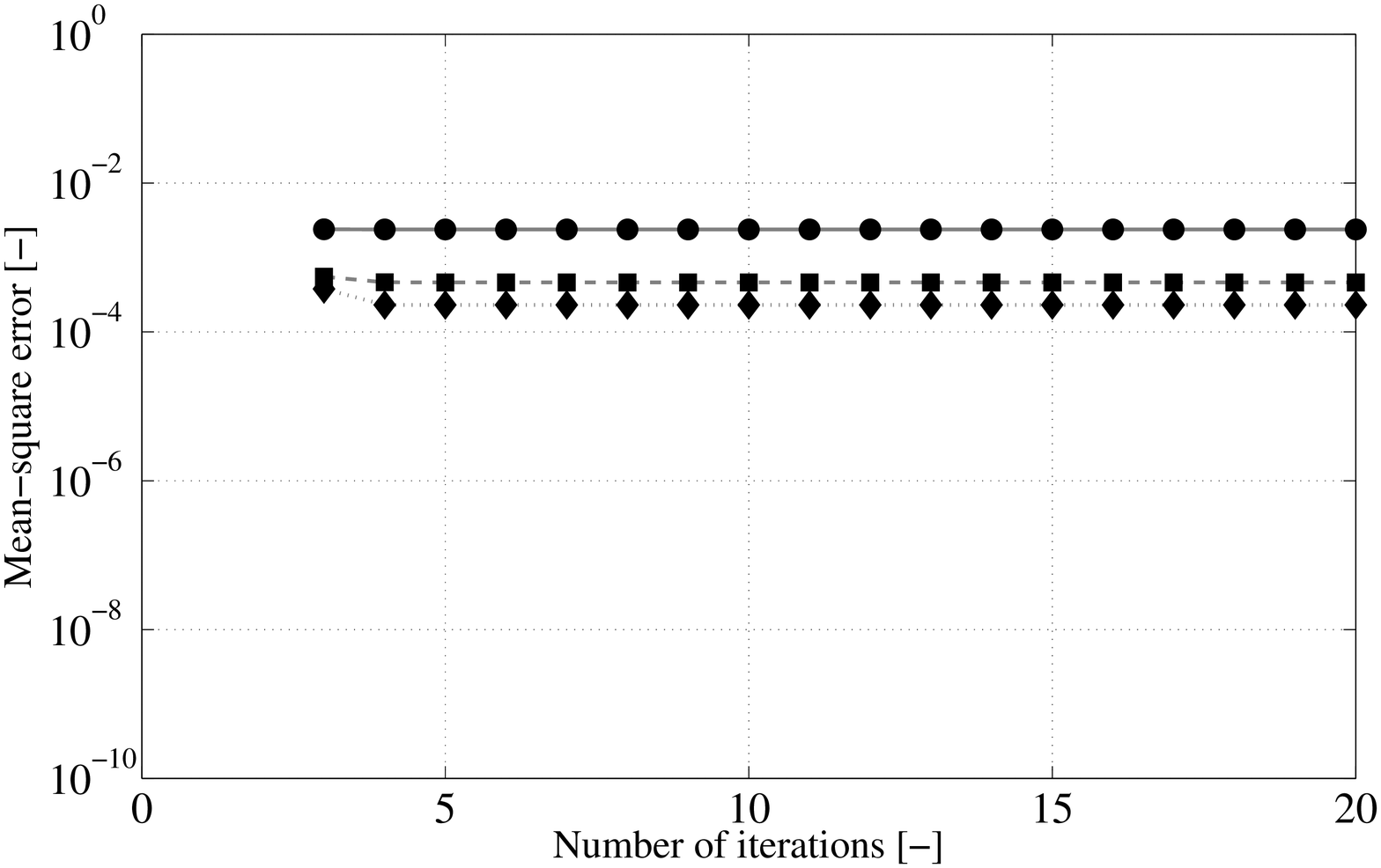}}
     \caption{Convergence analysis: (a) reduced dimension and (b, c) mean-square distance between the successive approximations determined by the simulation based on Monte Carlo sampling and the simulation based on polynomial chaos for $\epsilon_{1}=0.02$~(circles), $\epsilon_{1}=0.01$~(squares), and $\epsilon_{1}=0.005$~(diamonds) and $\epsilon_{2}=0.01$ as a function of the number of iterations.}\label{fig:figure10}
  \end{center}
\end{figure}

\begin{figure}[htp]
  \begin{center}
    \hspace{-13mm}\subfigure[Total degree $q$.]{\includegraphics[width=0.55\textwidth]{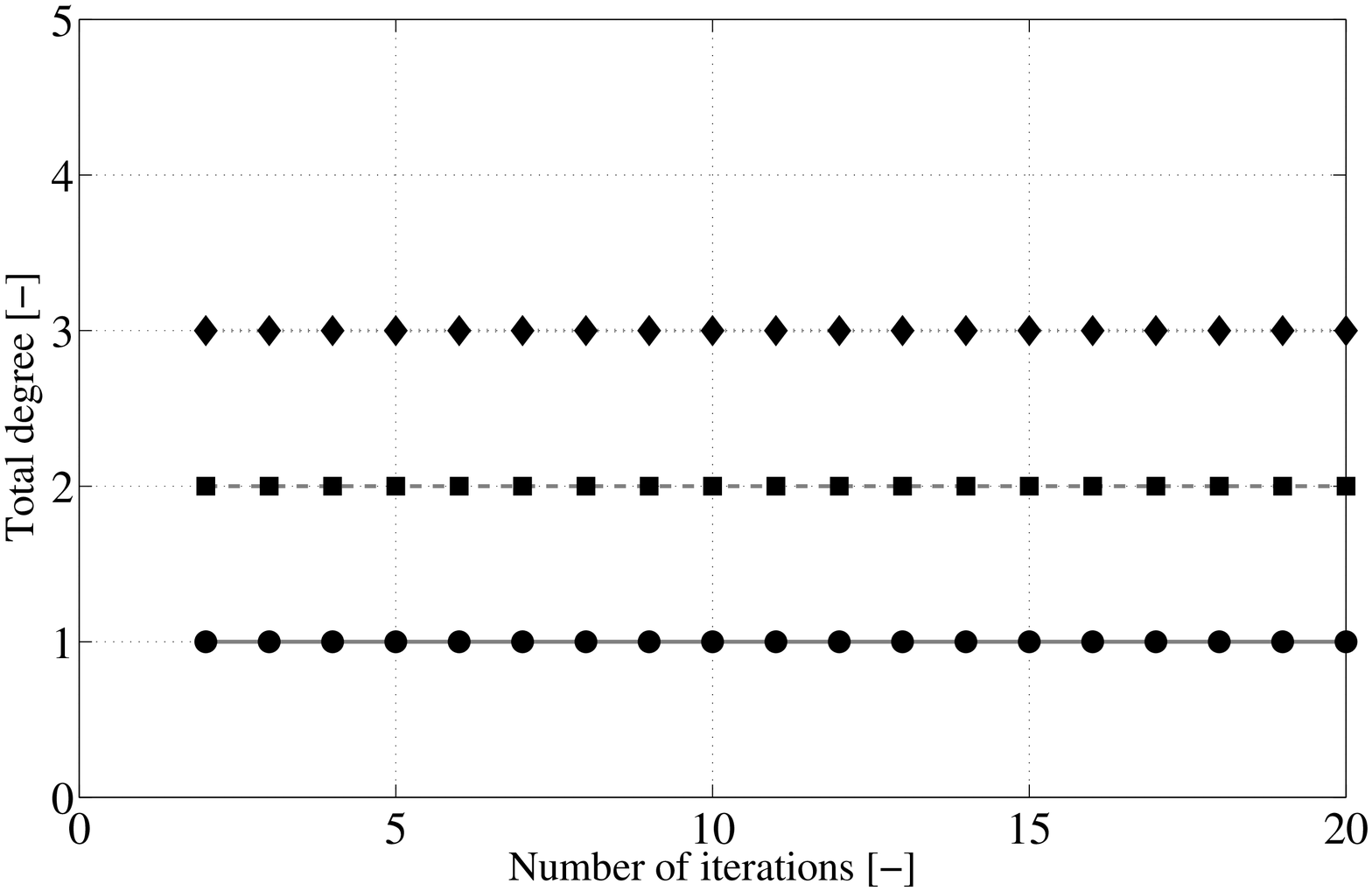}}
    \vfill
    \hspace{-13mm}\subfigure[$\ell\mapsto\sqrt{\frac{1}{MC}\sum_{k=1}^{MC}\|\boldsymbol{T}^{\ell}(\boldsymbol{\xi}_{k},\boldsymbol{\zeta}_{k})-\widehat{\boldsymbol{T}}{}^{\ell,p}(\boldsymbol{\xi}_{k},\boldsymbol{\zeta}_{k})\|_{\boldsymbol{W}}^{2}}\Big/\sqrt{\frac{1}{MC}\sum_{k=1}^{MC}\|\boldsymbol{T}^{\infty}(\boldsymbol{\xi}_{k},\boldsymbol{\zeta}_{k})\|_{\boldsymbol{W}}^{2}}$.]{\includegraphics[width=0.55\textwidth]{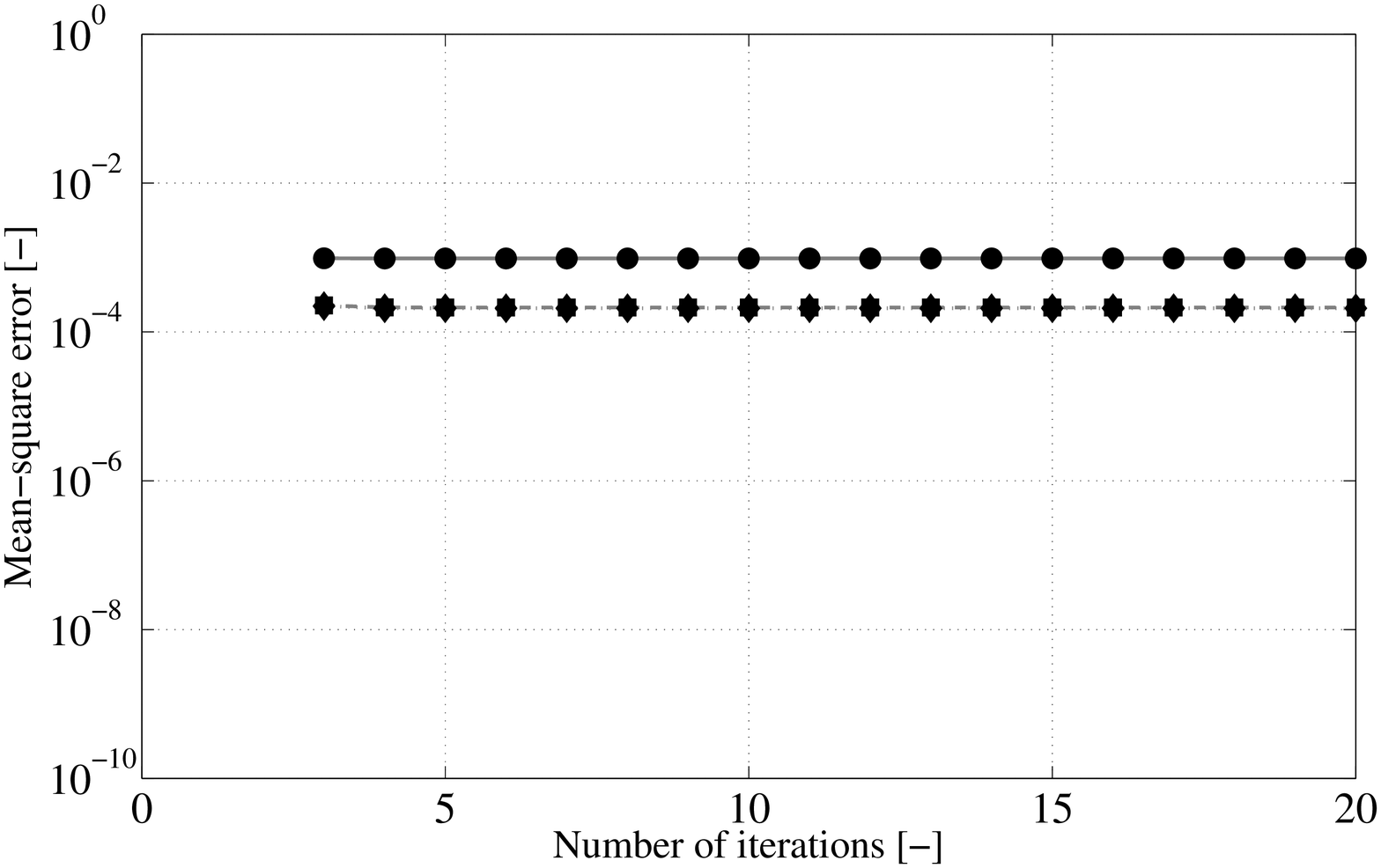}}
    \vfill
    \hspace{-13mm}\subfigure[$\ell\mapsto\sqrt{\frac{1}{MC}\sum_{k=1}^{MC}\|\boldsymbol{\Phi}^{\ell}(\boldsymbol{\xi}_{k},\boldsymbol{\zeta}_{k})-\widehat{\boldsymbol{\Phi}}{}^{\ell,q}(\boldsymbol{\eta}^{\ell,p}(\boldsymbol{\xi}_{k}),\boldsymbol{\zeta}_{k})\|_{\boldsymbol{W}}^{2}}\Big/\sqrt{\frac{1}{MC}\sum_{k=1}^{MC}\|\boldsymbol{\Phi}^{\infty}(\boldsymbol{\xi}_{k},\boldsymbol{\zeta}_{k})\|_{\boldsymbol{W}}^{2}}$.]{\includegraphics[width=0.55\textwidth]{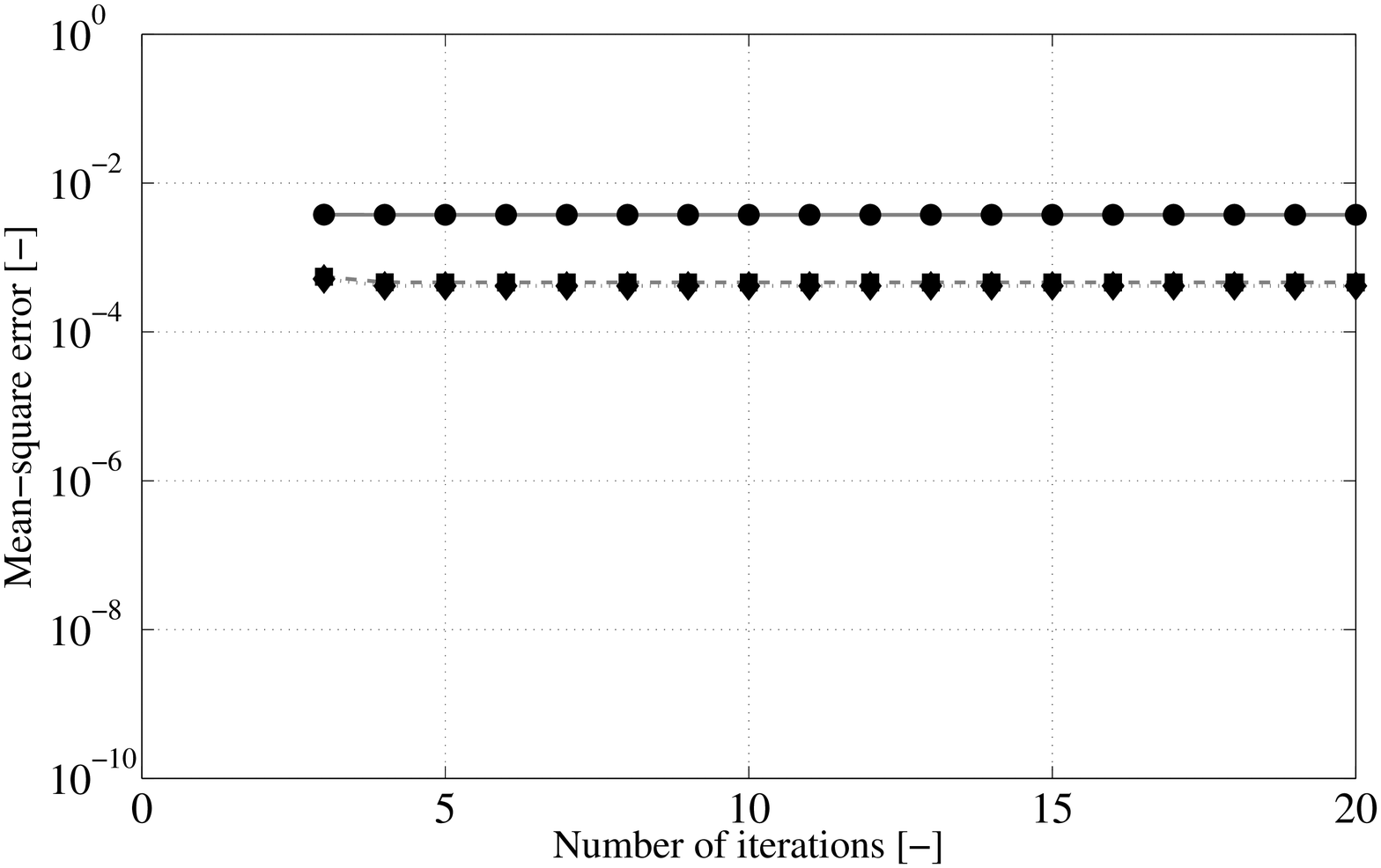}}
     \caption{Convergence analysis: (a) total degree and (b,c) mean-square distance between the successive approximations determined by the simulation based on Monte Carlo sampling and the simulation based on polynomial chaos for $\epsilon_{1}=0.01$ and $\epsilon_{2}=0.1$~(circles), $\epsilon_{2}=0.01$~(squares), and~$\epsilon_{2}=0.001$~(diamonds) as a function of the number of iterations.}\label{fig:figure11}
  \end{center}
\end{figure}

We repeated the simulation based on polynomial chaos for several values of the error tolerance levels.
Each error tolerance level corresponded to a specific accuracy that the reduced chaos expansion with random coefficients of the random temperature and the chaos expansion of the random neutron flux were required to maintain at each iteration. 
Figures~\ref{fig:figure10}(a) and~\ref{fig:figure11}(a) indicate that more terms were retained in these expansions when higher accuracy was required.

Further, Figs.~\ref{fig:figure10} ((b) and~(c)) and~\ref{fig:figure11} ((b) and~(c)) indicate that the distance between the successive approximations determined by the simulation based on Monte Carlo sampling and the simulation based on polynomial chaos remained bounded as the iterations progressed and that this distance can be reduced systematically by improving the accuracy of the reduced chaos expansion with random coefficients of the random temperature and the chaos expansion of the random neutron flux by decreasing the respective error tolerance levels.

\subsection{Concluding remarks}
The proposed methodology provided the solution of the neutronics subproblem in a reduced-dimensional space because the reduced chaos expansion with random coefficients could facilitate a low-dimensional representation of the random temperature as it passed from the heat subproblem to the neutronics subproblem.
While accuracy was maintained, the solution in a reduced-dimensional space resulted in computational gains because of the following two factors.
First, the solution in a reduced-dimensional space facilitated the accurate representation of the random neutron flux by a chaos expansion that contained only a few terms.
Second, the coefficients in the chaos expansion of the random neutron flux could be computed by using a quadrature rule that had only a few nodes, and thus, the solution of only a few samples of the neutronics subproblem was required at each iteration. 

\section{Conclusion}\label{sec:sec9}

We presented a characterization of information by a reduced chaos expansion with random coefficients as this information passes from a subproblem of a stochastic coupled problem to another and from iteration to iteration. 
This expansion provides a reduced-dimensional representation of the exchanged information, while maintaining segregation between statistically independent sources of uncertainty that stem from different subproblems.
Further, we presented a measure transformation that allows stochastic expansion methods to exploit this dimension reduction to obtain an efficient solution of subproblems in a reduced-dimensional space.
We showed that owing to the uncertainty source segregation, polynomial chaos and quadrature rules required for the implementation of this measure transformation can be readily obtained by tensorization.
The effectiveness of the methodology was demonstrated by considering a multiphysics problem in nuclear engineering.

\section*{Acknowledgements}
This work was supported by the Department of Energy (DOE) through an Applied Scientific Computing Research (ASCR) grant. 
The authors would also like to thank Professor Christian Soize for relevant discussions during the final stages of the preparation of this paper. 

\bibliography{multiphysics4}

\end{document}